\documentclass[10pt,letterpaper,twoside]{amsart}

\usepackage{multicol}
\usepackage{amssymb}
\usepackage{amsmath}
\usepackage{amsbsy}
\usepackage{graphicx}
\usepackage{amsthm}
\usepackage[latin1]{inputenc}

\newtheorem{thm}{Theorem}
\newtheorem{cor}[]{Corollary}
\newtheorem{lem}[]{Lemma}
\newtheorem{defn}[]{Definition}
\newtheorem{prop}[]{Proposition}
\theoremstyle{remark}
\newtheorem{rem}[]{Remark}

\newcommand{\sign}[1]{{\rm sign\/}(#1)}
\newcommand{\supp}[1]{{\rm supp\/}(#1)}

\newcommand{\lcm}{\text{lcm}}
\newcommand{\sop}{\text{\rm supp}}

\newcommand{\R}{\mathbb{R}}

\newcommand{\C}{\mathbb{C}}

\newcommand{\Co}{\text{\rm Co}}
\newcommand{\Z}{\mathbb{Z}}

\newcommand{\lnj}{\mathcal{A}_{{\bf n},j}}
\newcommand{\lnlj}{\mathcal{A}_{{\bf n}^{l},j}}
\newcommand{\lnmj}{\mathcal{A}_{{\bf n},-j}}
\newcommand{\lnlmj}{\mathcal{A}_{{\bf n}^{l},-j}}
\newcommand{\gnj}{{\mathcal{G}}_{{\bf n},j}}
\newcommand{\gnju}{{\mathcal{G}}_{{\bf n},j+1}}
\newcommand{\gnmj}{{\mathcal{G}}_{{\bf n},-j}}
\newcommand{\gnmju}{{\mathcal{G}}_{{\bf n},-j-1}}

\newcommand{\fp}{\hfill $\Box$}

\renewcommand{\O}{\mathcal{O}}

\begin{document}

\title[Mixed type multiple orthogonal polynomials]{Mixed type multiple orthogonal polynomials for two Nikishin
Systems}
\thanks{
The first three authors in alphabetical order were supported by
grants MTM 2006-13000-C03-02 of Ministerio de Ciencia y
Tecnolog\'{\i}a and CCG07-UC3M/ESP-3339 of Comunidad Aut\'onoma de
Madrid-Universidad Carlos III de Madrid. V. N. Sorokin received
support from grants RFBR-08-01-00317 and NSh-3906.2008.1.}

\author[U. Fidalgo]{U. Fidalgo Prieto}
\address{U. Fidalgo Prieto \\
Dpto. de Matem\'aticas\\
Universidad Carlos III de Madrid \\
c/ Universidad 30, 28911  Legan\'es, Spain.}
\email{ufidalgo@math.uc3m.es}

\author[A. L\'opez ]{A. L\'{o}pez Garc\'{\i}a}
\address{A. L\'opez Garc\'{\i}a\\
Department of Mathematics\\
Vanderbilt University\\
Nashville, TN 37240, USA.} \email{abey.lopez@Vanderbilt.edu}

\author[G. L\'opez ]{G. L\'opez Lagomasino}
\address{G. L\'opez Lagomasino\\
Dpto. de Matem\'aticas\\
Universidad Carlos III de Madrid \\
c/ Universidad 30, 28911, Legan\'es, Spain.}
\email{lago@math.uc3m.es}

\author[V. N. Sorokin]{V. N. Sorokin}
\address{V. N. Sorokin\\
Dept. of Function Theory and Functional Analysis\\
Moscow State University\\
119992 Leninskye Gory, Moscow, Russia.}
\email{vnsormm@mech.math.msu.su}

\begin{abstract}

We study the logarithmic and ratio asymptotic of linear forms
constructed from a Nikishin system which satisfy orthogonality
conditions with respect to a system of measures generated from
another Nikishin system. This construction combines type I and
type II multiple orthogonal polynomials. The logarithmic
asymptotic of the linear forms is expressed in terms of the
extremal solution of an associated vector valued equilibrium
problem for the logarithmic potential. The ratio asymptotic is
described by means of a conformal representation of an appropriate
Riemann surface of genus zero onto the extended complex plane.
\end{abstract}

\maketitle \vspace{1cm}


{\it A.M.S. Subject Classification.} Primary: 30E10, 42C05;
Secondary: 41A20.

\section{Introduction}\label{int}

Let $s$ be a finite positive  Borel measure  supported on a
compact subset $\sop (s)$ of the real line, and
$(w^1_0,\ldots,w^1_{m_1}), (w^2_0,\ldots,w^2_{m_2})$ be two
systems of continuous functions on $\sop (s)$. Fix ${\bf
n}_1=(n_{1,0},\,n_{1,1},\ldots,\,n_{1,m_1})\in\Z_+^{m_1+1}$ and
${\bf n}_2=(n_{2,0},\,n_{2,1},\ldots,\,n_{2,m_2})\in\Z_+^{m_2+1}$.
Set $|{\bf n}_1|=n_{1,0}+n_{1,1}+\cdots+n_{1,m_1}$, $|{\bf n}_2|=
n_{2,0}+\cdots+n_{2,m_2}$, and ${\bf n} = ({\bf n}_1;{\bf n}_2)$.
In the sequel, we suppose that $|{\bf n}_2| +1 = |{\bf n}_1|$.

Let $|{\bf n}_{1}| \geq 1.$ It is easy to see that there exist
polynomials $a_{{\bf n},0},\,a_{{\bf n},1},\ldots,\,a_{{\bf
n},m_1}$ such that:
\begin{itemize}
\item [i)] $\deg(a_{{\bf n},j})\le n_{1,j}-1,\,j=0,\ldots,\,m_1$,
not all identically equal to zero. \item [ii)] For
$k=0,\ldots,m_2$
\[
 \int x^{\nu}  \sum_{j=0}^{m_1} a_{{\bf
n},j}(x)w^1_j(x) w_k^2(x) d s (x)=0, \qquad \nu =0,\ldots,n_{2,k}
-1.
\]
\end{itemize}
($\deg (a_{{\bf n},j}) \leq -1$ means that $a_{{\bf n},j} \equiv
0.$)

When $m_2=0$ the polynomials $(a_{{\bf n},0},\ldots,a_{{\bf
n},m_1})$ are called type I multiple orthogonal polynomials. If
$m_1 =0, a_{{\bf n},0}$ is called a type II multiple orthogonal
polynomial. The case $m_1=m_2=0$ reduces to the usual definition
of orthogonal polynomial. When $m_1,m_2 \geq 1$ these multiple
orthogonal polynomials are called of mixed type.

Multiple orthogonal polynomials appear in problems connected with
the algebraic independence of functions and numbers (type I) and
in questions related with simultaneous rational approximation
(type II). Those of type II are formed by polynomials which share
orthogonality conditions with a system of measures which may be
written in  the form of orthogonality relations with respect to a
family of generalized polynomials (that in the sequel we call
linear forms). In type I,  the linear forms are defined through
full orthogonality relations with respect to a single measure.
Mixed type multiple orthogonal polynomials occur in stochastic
models connected with random matrices and non intersecting random
paths, see \cite{DK}. Mixed type multiple orthogonal polynomials
as presented above were considered in \cite{S1} and their
algebraic properties studied in \cite{SV}.

We will restrict our attention to mixed type multiple orthogonal
polynomials in which the linear forms are generated by two (not
necessarily distinct) Nikishin systems of measures. Nikishin
systems of measures were introduced in \cite{nik}. Before going
into details let us mention some papers which constitute our
starting point.

E. M. Nikishin studied the asymptotic behavior of the linear forms
generated by a Nikishin system of measures in \cite{nik2} (see
also \cite{lomi} and the last section in \cite{niso}). He
described the logarithmic asymptotic of type I multiple orthogonal
polynomials in terms of the solution of a vector equilibrium
problem for the logarithmic potential. Later,
Gonchar-Rakhmanov-Sorokin studied in \cite{GRS} the rate of
convergence of Hermite-Pad\'{e} approximation of generalized
Nikishin systems of functions and the logarithmic asymptotic of
their associated type II multiple orthogonal polynomials. The
solution is also characterized by a similar vector equilibrium
problem.  In \cite{S2}, V. N. Sorokin defines mixed type multiple
orthogonal polynomials for two Nikishin systems and gives their
logarithmic asymptotic.

Let $s$ be a finite positive Borel measure supported on a bounded
interval $\Delta$ of the real line ${\mathbb{R}}$ such that
$s^{\prime} > 0$ almost everywhere on $\Delta$ and let $\{Q_n\},
n\in {\mathbb{Z}}_+,$ be the corresponding sequence of monic
orthogonal polynomials; that is, with leading coefficients equal
to one. In a series of two papers (see \cite{kn:Rak1} and
\cite{kn:Rak2}), E. A. Rakhmanov proved that under these
conditions
\begin{equation} \label{Rakteo}
\lim_{n \in {\mathbb{Z}}_+} \frac{Q_{n+1}(z)}{Q_n(z)} =
\frac{\varphi(z)}{\varphi^{\prime}(\infty)}, \qquad \mathcal{K}
\subset {\mathbb{C}} \setminus \Delta
\end{equation}
(uniformly on each compact subset of ${\mathbb{C}} \setminus
\Delta$), where $\varphi(z)$ denotes the conformal representation
of $\overline{\mathbb{C}} \setminus \Delta$ onto $\{w: |w| > 1\}$
such that $\varphi(\infty) = \infty$ and $\varphi^{\prime}(\infty)
>0$. This result attracted great attention because of its
theoretical interest within the general theory of orthogonal
polynomials and its applications to the theory of rational
approximation of analytic functions. Simplified proofs of
Rakhmanov's theorem may be found in \cite{kn:Rak3} and
\cite{kn:Nev1}.

This result has been extended in several directions. Orthogonal
polynomials with respect to varying measures (depending on the
degree of the polynomial) arise in the study of multipoint Pad\'e
approximation of Markov functions. In this context, in
\cite{kn:Gui1} and \cite{kn:Gui2}, an analogue of Rakhmanov's
theorem for such sequences of orthogonal polynomials was proved.
Recently, S. A. Denisov \cite{Denisov} (see also \cite{NT})
extended Rakhmanov's result to the case when  $\sop (s) =
\widetilde{\Delta} \cup e \subset \mathbb{R}$, where
$\widetilde{\Delta}$ is a bounded interval, $e$ is a set without
accumulation points  in $\overline{\mathbb{R}} \setminus
\widetilde{\Delta} $, and $s^{\prime}
> 0$ a.e. on $\widetilde{\Delta}$. A version for orthogonal
polynomials with respect to varying Denisov type measures was
given in \cite{BCL}.

For multiple orthogonal polynomials associated with Nikishin
systems of measures an analogue of Rakhmanov's theorem was proved
in \cite{AptLopRoc} and extended in \cite{LL} to the case when the
measures in the Nikishin system are as those considered by
Denisov.

Let us define the notion of  Nikishin system of measures. Let
$\sigma_{\alpha}$, $\sigma_{\beta}$ be two measures with constant
sign supported on $\R$ and let $\Delta_{\alpha}$, $\Delta_{\beta}$
denote the smallest intervals containing their supports,
$\sop(\sigma_{\alpha})$ and $\sop(\sigma_{\beta})$, respectively.
We write $\Co(\sop(\sigma_{\alpha}))=\Delta_{\alpha}$. Assume that
$\Delta_{\alpha}\cap\Delta_{\beta}=\emptyset$ and define
\[
\langle\sigma_{\alpha},\,\sigma_{\beta}\rangle(x) :=\int
\frac{d\sigma_{\beta}(t)}{x-t}d\sigma_{\alpha}(x)=\widehat{\sigma}_{\beta}(x)d\sigma_{\alpha}(x).
\]
Therefore, $\langle\sigma_{\alpha},\,\sigma_{\beta}\rangle$ is a
measure with constant sign and support equal to that of
$\sigma_{\alpha}$.

For a system of intervals $\Delta_0,\ldots,\,\Delta_{m}$ contained
in $\R$ satisfying $\Delta_j\cap\Delta_{j+1}=\emptyset$,
$j=0,\ldots,\,m-1$, and finite Borel measures
$\sigma_0,\ldots,\,\sigma_{m}$ with constant sign in
$\Co(\sop(\sigma_j))=\Delta_j,$ such that each one has infinitely
many points in its support, we define recursively
\[
\langle\sigma_0,\,\sigma_1,\ldots,\,\sigma_j\rangle=\langle\sigma_0,\,
\langle\sigma_1,\ldots,\,\sigma_j\rangle\rangle,\quad
j=1,\ldots,\,m.
\]
We say that $
(s_0,\ldots,\,s_{m})=\mathcal{N}(\sigma_0,\ldots,\,\sigma_{m})$,
where
\[
s_0=\langle\sigma_0\rangle=\sigma_0,\quad
s_1=\langle\sigma_0,\,\sigma_1\rangle,\ldots,\quad
\,s_{m}=\langle\sigma_0,\ldots,\,\sigma_{m}\rangle
\]
is the Nikishin system of measures generated by
$(\sigma_0,\ldots,\sigma_{m})$. In the sequel, when referring to a
Nikishin system the condition  $\Delta_j \cap \Delta_{j+1} =
\emptyset, j =0,\ldots,m-1,$ is always assumed to hold.  Notice
that all the measures in a Nikishin system have the same support,
namely $\sop(\sigma_0)$.  We will denote $(s_{j,j} = \sigma_j)$
\[ s_{j,k} = \langle \sigma_j , \ldots, \sigma_k \rangle, \qquad 0 \leq j
\leq k \leq m.
\]

Take two  systems $S^1 = (s^1_0,\ldots,s^1_{m_1}) =
\mathcal{N}(\sigma^1_0,\ldots,\sigma^1_{m_1}), S^2 =
(s^2_0,\ldots,s^2_{m_2}) =
\mathcal{N}(\sigma^2_0,\ldots,\sigma^2_{m_2})$ generated by
$m_1+1$ and $m_2+1$ measures, respectively. The two systems need
not coincide, but we will assume that $\sigma^1_0 = \sigma^2_0$;
that is, both systems stem from the same basis measure.  The
smallest interval containing $\sop(\sigma^i_j)$  will be denoted
$\Co(\sop(\sigma^i_j))=\Delta^i_j$.

Fix ${\bf
n}_1=(n_{1,0},\,n_{1,1},\ldots,\,n_{1,m_1})\in\Z_+^{m_1+1}$ and
${\bf n}_2=(n_{2,0},\,n_{2,1},\ldots,\,n_{2,m_2})\in\Z_+^{m_2+1}$.
Set $|{\bf n}_1|=n_{1,0}+n_{1,1}+\cdots+n_{1,m_1}$, $|{\bf n}_2|=
n_{2,0}+\cdots+n_{2,m_2}$, and ${\bf n} = ({\bf n}_1;{\bf n}_2)$.
In the sequel, we suppose that $|{\bf n}_2| +1 = |{\bf n}_1|$.

Let $|{\bf n}_1| \geq 1.$ The system of polynomials $a_{{\bf
n},0},\,a_{{\bf n},1},\ldots,\,a_{{\bf n},m_1}$  satisfying:
\begin{itemize}
\item [i')] $\deg(a_{{\bf n},j})\le n_{1,j}-1,\,j=0,\ldots,\,m_1$,
not all identically equal to zero. \item [ii')] For
$k=0,\ldots,m_2$
\begin{equation}\label{definition}
 \int x^{\nu}\left(a_{{\bf n},0}(x)+\sum_{j=1}^{m_1} a_{{\bf
n},j}(x)\widehat{s}^1_{1,j}(x)\right) d s^2_{0,k}(x)=0, \qquad \nu
=0,\ldots,n_{2,k}  -1,
\end{equation}
\end{itemize}
($\deg (a_{{\bf n},j}) \leq -1$ means that $a_{{\bf n},j} \equiv
0$) is called a system of mixed type multiple orthogonal
polynomials relative to the multi-index ${\bf n}= ({\bf n}_1;{\bf
n}_2)$ and the pair $(S^1,S^2)$ of Nikishin systems. This concept
was first introduced in \cite{S2}.

Finding $a_{{\bf n},0},\ldots,a_{{\bf n},m_1}$ reduces to solving
a homogeneous linear system of  $|{\bf n}_2|$ equations on $|{\bf
n}_1|$ unknowns. Since $|{\bf n}_2| = |{\bf n}_1| -1$  a
non-trivial solution is guaranteed.

A multi-index ${\bf n}= ({\bf n}_1; {\bf n}_2)$ is said to be
normal if every solution to i')-ii') satisfies $\deg a_{{\bf n},j}
= n_{1,j} -1, j=0,\ldots,m$. If ${\bf n}$ is normal, it is easy to
verify that the vector $(a_{{\bf n},0},\ldots,a_{{\bf n},m_1})$ is
uniquely determined except for a constant factor, and in that case
we normalize it to be ``monic'' meaning by this that its last
entry different from zero has leading coefficient equal to $1$.
Set
\[ \mathbb{Z}_+^{m_1+1}(\bullet) = \{ {\bf n}_1 \in
\mathbb{Z}_+^{m_1 + 1}: n_{1,0}\ge \cdots \ge n_{1,m_1} \}.
\]
In Proposition \ref{normality}, we  prove that  all  ${\bf n} =
({\bf n}_1;{\bf n}_2) \in\mathbb{Z}_+^{m_1+1}(\bullet) \times
\mathbb{Z}_+^{m_2 +1}(\bullet)$ are normal. For the sequences of
multi-indices we shall consider, for almost all ${\bf n}$ we will
have that $n_{1,m_1} \geq 1$ and a ``monic''  $(a_{{\bf
n},0},\,a_{{\bf n},1},\ldots,\, a_{{\bf n},m_1})$ will have
$a_{{\bf n},m_1}$ monic.

Theorem 1 gives the rate of convergence of the $|{\bf n}_1|$-th
root of the linear forms
\[
\mathcal{A}_{{\bf n},0}(z)= a_{{\bf n},0}(z)+\sum_{k=1}^{m_1}
a_{{\bf n},k}(z)\widehat{s}^1_{1,k}(z),
\]
under mild conditions on the sequence of multi-indices  and the
measures generating both Nikishin systems.   A measure $\sigma$ is
said to be regular if
\[ \lim_{n \to \infty} \kappa_n^{1/n} = 1/\mbox{cap}(\sop(\sigma)),
\]
where $\mbox{cap}(\cdot)$ denotes the logarithmic capacity of the
Borel set $(\cdot)$ and $\kappa_n$ denotes the leading coefficient
of the $n$th orthonormal polynomial with respect to $\sigma$. For
different equivalent forms of defining regular measures see
sections 3.1 to 3.3 in \cite{stto} (in particular Theorem 3.1.1).
For short, we write $ (S^1,S^2) \in \mbox{\bf Reg}$ to mean that
all the measures which generate both Nikishin systems are regular
and their supports are regular compact sets. Recall that a compact
set is regular when the Green's function with singularity at
$\infty$ of the unbounded connected component of the complement of
the compact set can be extended continuously to all $\mathbb{C}$.
Before stating Theorem 1, we need to introduce some notation and
results from potential theory.

Let $E_k,\, k=-m_2,\ldots,m_1,$ be (not necessarily distinct)
compact subsets of the real line and $\mathcal{C} = (c_{j,k}),
-m_2\leq j,k \leq m_1,$ a real, positive definite, symmetric
matrix of order $m_1+m_2+1$. $\mathcal{C}$ will be called the
interaction matrix. By $\mathcal{M}(E_k)$ we denote the class of
all finite, positive, Borel measures with compact support
consisting of an infinite set of points contained in $E_k$ and
$\mathcal{M}_1(E_k) $ is the subclass of probability measures of
$\mathcal{M}(E_k).$ Set
\[\mathcal{M}_1= \mathcal{M}_1(E_{-m_2})
\times \cdots \times \mathcal{M}_1(E_{m_1})  \,.
\]

Given a vector measure $\mu=(\mu_{-m_2},\ldots,\,\mu_{m_1}) \in
\mathcal{M}_1$ and $j= -m_2,\ldots,m_1,$ we define the combined
potential
\begin{equation}
\label{combpot} W^{\mu}_j(x) = \sum_{k=-m_2}^{m_1} c_{j,k}
V^{\mu_k}(x) \,,
\end{equation}
where
\[ V^{\mu_k}(x) = \int \log \frac{1}{|x-t|} \,d\mu_k(t)\,,
\]
denotes the standard logarithmic potential of $\mu_k$. We denote
\[ \omega_j^{\mu} = \inf \{W_j^{\mu}(x): x \in E_j\} \,, \quad
j=-m_2,\ldots,m_1\,.
\]

In Chapter 5 of \cite{niso}  the authors prove (we state the
result in a form convenient for our purpose).

\begin{lem} \label{niksor} Assume that the compact sets
$E_k,k=-m_2,\ldots,m_1,$ are regular with respect to the Dirichlet
problem. Let $\mathcal{C}$ be a real, positive definite, symmetric
matrix of order $m_1+m_2+1$. If there exists $\overline{\mu} =
(\overline{\mu}_{-m_2},\ldots,\overline{\mu}_{m_1})\in
\mathcal{M}_1$ such that for each $j=-m_2,\ldots,m_1$
\[
W_j^{\overline{\mu}} (x) = \omega_j^{\overline{\mu}}\,, \qquad x
\in \supp{\overline{\mu}_j}\,,
\]
then $\overline{\mu}$ is unique. Moreover, if $c_{j,k} \geq 0$
when $E_j \cap E_k \neq \emptyset$, then $\overline{\mu}$ exists.
\end{lem}

For details on how this lemma is derived from  \cite[Chapter
5]{niso} see  \cite[Section 4]{bel}. The vector measure
$\overline{\mu} \in \mathcal{M}_1$ is called the equilibrium
solution for the vector potential problem determined by the
interaction matrix $\mathcal{C}$ on the system of compact sets
$E_j\,, j = -m_2,\ldots,m_1\,.$

Let $\Lambda =
\Lambda(p_{1,0},\ldots,p_{1,m_1};p_{2,0},\ldots,p_{2,m_2}) \subset
\mathbb{Z}_+^{m_1+1}(\bullet) \times \mathbb{Z}_+^{m_2
+1}(\bullet)$ be an infinite sequence of distinct multi-indices
such that
\[
\lim_{ {\bf n} \in \Lambda} \frac{n_{1,j}}{|{\bf n}_1|} = p_{1,j}
\in (0,1), \quad j=0,\ldots,m_1, \quad \lim_{{\bf n} \in \Lambda}
\frac{n_{2,j}}{|{\bf n}_1|} = p_{2,j} \in (0,1), \quad
j=0,\ldots,m_2.
\]
Obviously, $p_{1,0}\ge \cdots\ge p_{1,m_1}, p_{2,0} \geq \cdots
\geq p_{2,m_2}$, and $\sum_{j=0}^{m_1} p_{1,j}= \sum_{j=0}^{m_2}
p_{2,j} = 1$. Set
\[
P_j=\sum_{k=j}^{m_1} p_{1,k},\,\, j=0,\ldots,m_1, \qquad P_{-j} =
\sum_{k=j}^{m_2} p_{2,k},\,\, j=0,\ldots,m_2.
\]

Let us define the interaction matrix $\mathcal{C}$ which is
relevant for the rest of the paper. Take the tri-diagonal matrix
\begin{equation}\label{matriz}
\mathcal{C}=
\begin{pmatrix}
P_{-m_2}^2 & -\frac{P_{-m_2}P_{-m_2+1}}{2} & 0 &  \cdots  & 0\\
-\frac{P_{-m_2}P_{-m_2+1}}{2} & P_{-m_2+1}^2 &
-\frac{P_{-m_2+1}P_{-m_2+2}}{2} & \cdots  & 0\\
0 & -\frac{P_{-m_2+1}P_{-m_2+2}}{2} & P_{-m_2+2}^2
& \cdots  & 0 \\
\vdots &\vdots&\vdots&\ddots &\vdots\\
0 & 0& 0 &\cdots  & P_{m_1}^2
\end{pmatrix}.
\end{equation}
This matrix satisfies all the assumptions of Lemma \ref{niksor} on
the compact sets $E_j = \sop(\sigma^1_j), j=0,1,\ldots,m_1, E_j =
\sop(\sigma^2_{-j}), j=0,-1,\ldots,-m_2,$ including $c_{j,k} \geq
0$ when $E_j \cap E_k \neq \emptyset$ (recall that $\sigma^1_0 =
\sigma^2_0$), and it is positive definite because the principal
section $\mathcal{C}_r, r=1,\ldots,m_1+m_2+1,$ of $\mathcal{C}$
satisfies
\begin{displaymath}
\det(\mathcal{C}_r)=P_{-m_2}^2\cdots P_{-m_2+r-1}^2 \det
\begin{pmatrix}
1 & -\frac{1}{2} & 0 & \cdots & 0 & 0 \\
-\frac{1}{2} & 1 & -\frac{1}{2}  & \cdots & 0 & 0 \\
0 & -\frac{1}{2} &
1  & \cdots & 0 & 0 \\
\vdots & \vdots & \vdots   & \ddots & \vdots & \vdots\\
0 & 0 & 0 & \cdots & 1 & -\frac{1}{2}\\
0 & 0 & 0  & \cdots & -\frac{1}{2} & 1
\end{pmatrix}_{r\times r} > 0.
\end{displaymath}
Let $\overline{\mu}(\mathcal{C})$ be the equilibrium solution for
the corresponding vector potential problem.  We have

\begin{thm}\label{principal}
Let $\Lambda =
\Lambda(p_{1,0},\ldots,p_{1,m_1};p_{2,0},\ldots,p_{2,m_2}) \subset
\mathbb{Z}_+^{m_1+1}(\bullet) \times
\mathbb{Z}_+^{m_2+1}(\bullet), (S^1,S^2)\in\mbox{\bf {Reg}},$ $
S^1= \mathcal{N}(\sigma^1_0,\ldots,\sigma^1_{m_1}),$ and $S^2=
\mathcal{N}(\sigma^2_0,\ldots,\sigma^2_{m_2})$ be given.  Then
\begin{equation}\label{conver1a}
\lim_{{\bf n} \in \Lambda}\left|\mathcal{A}_{{\bf
n},0}(z)\right|^{1/|{\bf n}_1|}=G_0(z), \qquad {\mathcal{K}}
\subset {\C}\setminus (\Delta^1_0\cup\Delta^1_1),
\end{equation}
uniformly on each compact subset ${\mathcal{K}} \subset
{\C}\setminus (\Delta^1_0\cup\Delta^1_1)$, where
\[
G_0(z)=\exp\left(P_1V^{\overline{\mu}_1}(z)-V^{\overline{\mu}_{0}}(z)-2\sum_{k=1}^{m_1}
\frac{\omega_k^{\overline{\mu}}}{P_k}\right).
\]
$\overline{\mu}= \overline{\mu}(\mathcal{C})=
(\overline{\mu}_{-m_2},\ldots,\overline{\mu}_{m_1})$ is  the
equilibrium vector measure   and
$(\omega_{-m_2}^{\overline{\mu}},\ldots,\omega_{m_1}^{\overline{\mu}})$
is the system of equilibrium constants for the vector potential
problem determined by the interaction matrix $\mathcal{C}$ defined
in $(\ref{matriz})$ on the system of compact sets  $E_j =
\sop(\sigma^1_j),j=0,\ldots,m_1, E_j = \sop(\sigma_{-j}^2),
j=-m_2,\ldots,0$.
\end{thm}

Throughout the paper, the notation
\[ \lim_{n \in \Lambda} g_n(z) = g(z), \qquad
{\mathcal{K}} \subset \Omega, \] stands for uniform convergence of
the sequence $\{g_n\}, n\in \Lambda,$ to $g$ on each compact
subset $\mathcal{K}$ contained in the indicated region (in this
case $\Omega$).

For the next result, we assume that $\sop (\sigma^{i}_j) =
\widetilde{\Delta}^{i}_j \cup e^{i}_j, i=1,2$, where
$\widetilde{\Delta}^{i}_j$ is a bounded interval of the real line,
$|(\sigma^{i}_j)^{\prime}| > 0$ a.e. on
$\widetilde{\Delta}_j^{i}$, and $e^{i}_j$ is at most a denumerable
set without accumulation points in $\mathbb{R} \setminus
\widetilde{\Delta}_j^{i}$. We denote this writing $S^1=
\mathcal{N}'(\sigma^1_0,\ldots,\sigma^1_{m_1}), S^2=
\mathcal{N}'(\sigma^2_0,\ldots,\sigma^2_{m_2})$. Fix a vector
$l:=(l_{1};l_{2})$ where $0\leq l_{1}\leq m_{1}$ and $0\leq
l_{2}\leq m_{2}$. We define the multi-index ${\bf n}^{l}:=({\bf
n}_{1}+{\bf e}^{l_{1}};{\bf n}_{2}+{\bf e}^{l_{2}})=({\bf
n}_{1}^{l_{1}};{\bf n}_{2}^{l_{2}})$, where ${\bf e}^{l_{i}}$
denotes the unit vector of length $m_{i}+1$ with all components
equal to zero except the component $(l_{i}+1)$ which equals $1$.
It is always assumed that both $\bf n$ and ${\bf n}^{l}$ belong to
$\mathbb{Z}_{+}^{m_{1}+1}(\bullet)\times\mathbb{Z}_{+}^{m_{2}+1}(\bullet)$.

\begin{thm} \label{teofund1}
Let $S^1= \mathcal{N}'(\sigma^1_0,\ldots,\sigma^1_{m_1}), S^2=
\mathcal{N}'(\sigma^2_0,\ldots,\sigma^2_{m_2})$, and $\Lambda
\subset \mathbb{Z}_+^{m_1+1}(\bullet)\times
\mathbb{Z}_+^{m_2+1}(\bullet)$ is an infinite sequence of distinct
multi-indices such that
\[
\displaystyle{\sup_{{\bf n} \in \Lambda }((m_2 +1) n_{2,0} - |{\bf
n}_2| ) }< \infty, \qquad \displaystyle{\sup_{{\bf n} \in \Lambda
}((m_1+1)n_{1,0} -  |{\bf n}_1| ) }< \infty.
\]
Assume that there exists $l =(l_{1};l_{2}), 0\leq l_{1}\leq m_{1},
0\leq l_{2}\leq m_{2},$ such that for all ${\bf n} \in \Lambda$ we
have that ${\bf n}^l = ({\bf n}_1^{l_1};{\bf n}_2^{l_2}) \in
\mathbb{Z}_{+}^{m_1+1}(\bullet)\times
\mathbb{Z}_{+}^{m_2+1}(\bullet)$. Then
\begin{equation}\label{eq30*}
\lim_{{\bf n}\in\Lambda}\,\frac{{\mathcal{A}}_{{\bf
n}^l,0}(z)}{{\mathcal{A}}_{{\bf
n},0}(z)}={\mathcal{A}}^{(l)}_{0}(z),\qquad  \mathcal{K}
\subset\mathbb{C}\setminus(\supp{\sigma_{0}^{1}}\cup\supp{\sigma_{1}^{1}})\,,
\end{equation}
where ${\mathcal{A}}^{(l)}_{0}(z)$ is a one to one analytic
function in $\mathbb{C} \setminus
(\widetilde{\Delta}_{0}^{1}\cup\widetilde{\Delta}_{1}^{1})$.
\end{thm}

An expression for ${\mathcal{A}}^{(l)}_{0}(z)$ will be given in
Theorem \ref{teoLratio}. The answer depends of the conformal
representation of an associated Riemann surface with $m_1+m_2+2$
sheets and genus zero onto the extended complex plane. The
previous result is already new when $m_2 = 0$ and $m_1 \geq 2.$

Besides normality, in Section 2 we obtain the orthogonality
relations satisfied by the linear forms involved in the
construction. Section \ref{ext-pro} is devoted to the study of the
asymptotic distribution of zeros of a system of linear forms
associated with $\mathcal{A}_{{\bf n},0}$ that allows to prove
Theorem \ref{general} in Section \ref{thm1} of which Theorem
\ref{principal} is a corollary. Theorem 5 was first stated in
\cite{S1} under the stronger assumptions $S^1=
\mathcal{N}'(\sigma^1_0,\ldots,\sigma^1_{m_1}), S^2=
\mathcal{N}'(\sigma^2_0,\ldots,\sigma^2_{m_2})$, and  $\sop
(\sigma^{i}_j) = \widetilde{\Delta}^{i}_j, i=1,2$; and the proof
was only carried out for $m_1 =m_2 =1.$

In Section \ref{intprop} we study the interlacing properties of
the zeros of the linear forms which is needed for the proof of
Theorem \ref{teoLratio} in Section \ref{ratio} from which Theorem
\ref{teofund1} follows.  Section 7 contains a Markov type theorem
for mixed type Hermite-Pad\'{e} approximation and some
reinterpretation of the theory developed in the context of systems
of bi-orthogonal linear forms.

\section{Normality and orthogonality relations}\label{normal}

Recall that
\[
s_{j,k}=\langle \sigma_j,\ldots,\sigma_k\rangle,\quad 0\le j \le k
\le m,\qquad s_{j,j}=\sigma_j \,.
\]
We denote $(\widehat{s}^1_{j+1,j}(z)\equiv 1, \mathcal{A}_{{\bf
n},m_1} \equiv a_{{\bf n},m_1})$
\[
\mathcal{A}_{{\bf n},j}(z):=\sum_{k=j}^{m_1} a_{{\bf
n},k}(z)\widehat{s}^1_{j+1,k}(z),\qquad j=0,\ldots,\,m_1.
\]

In \cite{nik}, E. M. Nikishin introduced the following definition.

\begin{defn}
A set of real continuous functions $u_0(x),\ldots,\,u_{m_1}(x)$
defined on an interval $\Delta$, is called an AT-system for  ${\bf
n}_1 = (n_{1,0},\ldots,n_{1,m_1}) \in\Z_+^{m_1+1}$, if for any
polynomials $h_0\ldots,\,h_{m_1}$ such that $\deg(h_i)\le
n_{1,i}-1$, $i=0,\ldots,m_1$, not simultaneously identically equal
to zero, the function
\[
 h_0(x)u_0(x)+\cdots+h_{m_1}(x)u_{m_1}(x),
\]
has at most $|{\bf n}_1|-1$ zeros on $\Delta$ $(\deg (h_j)\le -1$
means that $h_j\equiv 0)$.
\end{defn}

Let $\Z_+^{m_1+1}(*)$ be the set of multi-indices given by
\[
\Z_+^{m_1+1}(*)=\{{\bf n}_1\in\Z^{m_1+1}_{+}:\, \not\exists\,\,
i<k<j\,\, \mbox{such that}\,\, n_{i} < n_{j} < n_{k}\}.
\]
In connection with AT-systems, in \cite{uli} U. Fidalgo  and G.
L\'opez  proved

\begin{lem}\label{ulises} Let ${\bf
n}_1\in\Z_+^{m_1+1}(*)$ and
$(s_1,\ldots,\,s_{m_1})=\mathcal{N}(\sigma_1,\ldots,\,\sigma_{m_1})$,
then the system of functions
$(1,\,\widehat{s}_1,\ldots,\,\widehat{s}_{m_1})$ defines an
AT-system with respect to ${\bf n}_1 = (n_{1,0},
\ldots,n_{1,m_1})$ on any interval disjoint from
$\Co(\supp{\sigma_1})$.
\end{lem}

For each $j = 0,\ldots,m_1-1$, we have that
$(s^1_{j+1,j+1},\ldots,s^1_{j+1,m_1}) =
\mathcal{N}(\sigma^1_{j+1},\ldots,\sigma^1_{m_1})$. Using  Lemma
\ref{ulises} it follows that for ${\bf n}_1 \in
\mathbb{Z}_+^{m_1+1}(\bullet) \subset \mathbb{Z}_+^{m_1+1}(*)$ the
linear form $\mathcal{A}_{{\bf n},j}$ cannot have more that
$N_{1,j} -1$ zeros on any interval disjoint from $\Delta^1_{j+1}$,
where \[ N_{1,j} = N_{1,j}({\bf n}) = n_{1,j}+\cdots
+n_{1,m_{1}}\,.\] Obviously, the same is true for the polynomial
$\mathcal{A}_{{\bf n},m_1}\equiv a_{{\bf n},m_1}$. Below, we also
use the previous lemma for linear forms generated by the second
Nikishin system.

Notice that $d s^2_{0,k}(x) = \widehat{s}^2_{1,k}(x)  d
\sigma_{0}^2(x)$. On the other hand, we can replace $x^{\nu}$ by
any polynomial of degree $\leq n_{2,k}-1$ inside the integral in
(\ref{definition}). Set
\[
{\mathcal{B}}_{{\bf n}_2}(z)= \sum_{k=0}^{m_2} b_{{\bf
n}_2,k}(z)\widehat{s}^2_{1,k}(z), \quad \deg b_{{\bf n}_2,k} \leq
n_{2,k} -1, \quad k =0,\ldots,m_2.
\]
$(\widehat{s}^2_{1,0}(z)\equiv 1).$ Then  (\ref{definition}) is
equivalent to
\begin{equation}\label{definition2}
 \int {\mathcal{B}}_{{\bf n}_2}(x)\mathcal{A}_{{\bf n},0}(x)d \sigma_{0}^2(x)=0,
\end{equation}
for all ${\mathcal{B}}_{{\bf n}_2}$ as indicated.

Suppose that $\mathcal{A}_{{\bf n},0}$ has at most $N < |{\bf
n}_1| -1 = |{\bf n}_2|$ sign changes on the interval $\Delta^2_0$.
Choose polynomials $b_{{\bf n}_2,k}$ conveniently so that
${\mathcal{B}}_{{\bf n}_2}$ changes sign exactly at those points
where $\mathcal{A}_{{\bf n},0}$ changes sign on $\Delta^2_0$ and
has a zero of order $|{\bf n}_2| -1 -N$ at one of the extreme
points of $\Delta^1_0=\Delta_{0}^{2}$. By Lemma \ref{ulises}, the
linear form ${\mathcal{B}}_{{\bf n}_2}$ has on $\Delta^2_0$ at
most $|{\bf n}_2| -1$ zeros, thus it can only have those zeros
which we have assigned to it. The continuous function
${\mathcal{B}}_{{\bf n}_2}\mathcal{A}_{{\bf n},0}$ has constant
sign on $\Delta^2_0$. This contradicts (\ref{definition2}).

We have proved that   $\mathcal{A}_{{\bf n},0}$ has $|{\bf n}_1|
-1$ zeros with odd multiplicity in the interior of $\Delta^2_0 =
\Delta^1_0$. In connection with intervals of the real line, the
interior refers to the Euclidean topology of $\mathbb{R}$. In
short, we shall see that $\mathcal{A}_{{\bf n},0}$ has no other
zeros in ${\mathbb{C}} \setminus \Delta^1_1$ and that they are all
simple. Before proving this, let us turn to the question of
normality.

\begin{prop}\label{normality} Let ${\bf n}
\in\mathbb{Z}_+^{m_1+1}(\bullet)\times
\mathbb{Z}_+^{m_2+1}(\bullet),$   $S^1=
\mathcal{N}(\sigma^1_0,\ldots,\,\sigma^1_{m_1}),$ and $S^2=
\mathcal{N}(\sigma^2_0,\ldots,\,\sigma^2_{m_2}),$ be given. Then,
${\bf n}$ is normal and $(a_{{\bf n},0}, \ldots, a_{{\bf n},m_1})$
is uniquely determined except for a constant factor.
\end{prop}

\begin{proof} Assume that there exists $j \in \{0,\ldots,m_1\}$ such
that $\deg a_{{\bf n},j}\leq n_{1,j} -2$. Then ${\bf n}_1 - {\bf
e}^j \in \mathbb{Z}_+^{m_1+1}(*)$, where ${\bf e}^j$ denotes the
$m_1+1$ dimensional unit vector with all components equal to zero
except the component $j+1$ which equals $1$. According to Lemma
\ref{ulises} applied to ${\bf n}_1 - {\bf e}^j$, the linear form
$\mathcal{A}_{{\bf n},0}$ has at most $|{\bf n}_1| -2$ zeros on
$\Delta^1_0$, but we pointed out before that it has at least
$|{\bf n}_1|-1$ sign changes on this interval. This contradiction
yields that for all $j \in \{0,\ldots,m_1\},$ $\deg a_{{\bf n},j}=
n_{1,j} -1,$ which implies normality.

Now, let us assume that $(a_{{\bf n},0},\ldots, a_{{\bf n},m_1})$
and $(a_{{\bf n},0}^{*}, \ldots, a_{{\bf n},m_1}^{*})$ solve
i')-ii') and these vectors are not collinear. According to what we
just proved, for all $j \in \{0,\ldots,m_1\},$ $\deg a_{{\bf n},j}
= \deg a_{{\bf n},j}^{*} = n_{1,j}-1.$ Take $\lambda \in
\mathbb{C} \setminus \{0\}$ such that $\deg (a_{{\bf n},0} -
\lambda a_{{\bf n},0}^{*}) \leq n_{1,0} - 2$. Obviously, the
vector $(a_{{\bf n},0}-\lambda a_{{\bf n},0}^{*},\ldots, a_{{\bf
n},m_1}-\lambda a_{{\bf n},m_1}^{*})$ is not identically equal to
zero and also solves i')-ii') which is not possible since all non
trivial solutions must have all components of maximal degree.
\end{proof}

Proposition \ref{normality} allows us to determine the ``monic''
$(a_{{\bf n},0},\,a_{{\bf n},1},\ldots,\, a_{{\bf n},m_1})$
uniquely and we impose this normalization .  The next lemma will
be used on several occasions.

\begin{lem} \label{reduc}
Let $s_k,k=1,\ldots,m,$ be finite signed Borel measures, $
\mbox{\rm Co}(\sop(s_k)) = \Delta \subset {\mathbb{R}}$. Let
$F(z)= f_0(z) + \sum_{k=1}^m f_k(z)\widehat{ s_k}(z)  \in
{\mathcal{H}}(\overline{\mathbb{C}} \setminus \Delta),$ where $f_k
\in {\mathcal{H}}(V), k=0,\ldots,m,$ and $V$ is a neighborhood of
$\Delta$. If $F(z) = {\mathcal{O}}(1/z^2), z \to \infty,$ then
\begin{equation} \label{eq:f}
\sum_{k=1}^m \int  f_k(x)d { s_k}(x) =0
\end{equation}
and $F(z) = {\mathcal{O}}(1/z), z \to \infty,$ implies that
\begin{equation} \label{eq:g}
F(z) = \sum_{k=1}^m \int \frac{f_k(x)d { s_k}(x)}{z-x}.
\end{equation}
\end{lem}

\begin{proof} Let $\Gamma \subset V$ be a closed smooth Jordan curve
that surrounds $\Delta$. If $F(z) = {\mathcal{O}}(1/z^2), z \to
\infty,$ from Cauchy's theorem, Fubini's theorem and Cauchy's
integral formula, it follows that
\[ 0 = \int_{\Gamma} F(z) dz = \sum_{k=1}^m \int_{\Gamma} f_k(z)\widehat{ s_k}(z) dz = \sum_{k=1}^m \int
\int_{\Gamma} \frac{f_k(z)dz}{z-x}ds_k(x) = 2\pi i \sum_{k=1}^m
\int f_k(x) ds_k(x),
\]
and we obtain (\ref{eq:f}). On the other hand, if $F(z) =
{\mathcal{O}}(1/z), z \to \infty,$  and we assume that $z$ is in
the unbouded connected component of the complement of $\Gamma$,
Cauchy's integral formula and Fubini's theorem render
\[ F(z) = \frac{1}{2\pi i} \int_{\Gamma} \frac{F(\zeta)d\zeta}{z-\zeta} = \frac{1}{2\pi i} \sum_{k=1}^m \int_{\Gamma}
\frac{f_k(\zeta)\widehat{s}_k(\zeta)  d\zeta}{z-\zeta} = \] \[
\sum_{k=1}^m \int \frac{1}{2\pi i} \int_{\Gamma} \frac{f_k(\zeta)
d\zeta}{(z-\zeta)(\zeta -x)} ds_k(x) = \sum_{k=1}^m \int
\frac{f_k(x)ds_k(x)}{z-x}
\]
which is (\ref{eq:g}). \end{proof}

In the sequel,  ${\bf n} = ({\bf n}_1;{\bf
n}_2)\in\mathbb{Z}_+^{m_1+1}(\bullet)\times
\mathbb{Z}_+^{m_2+1}(\bullet)$ and $n_{1,m_1} \geq 1$. For
$j=0,\ldots,\,m_1$, let $Q_{{\bf n},j}$ be the monic polynomial
whose zeros are those of the linear form $\mathcal{A}_{{\bf n},j}$
in the region $ {\mathbb{C}} \setminus \Delta^1_{j+1}$, counting
multiplicities $(\Delta^1_{m_1+1} = \emptyset)$. In particular,
$\mathcal{A}_{{\bf n},m_1} = a_{{\bf n},m_1}= Q_{{\bf n},m_1}$.
From the previous proposition, if $n_{1,m_1} \geq 1,$ $\infty$ is
not a zero of any one  of these linear forms; thus, $\infty$
cannot be an accumulation point of such zeros. Though it is not
the case, in principle, some of these linear forms may have an
infinite number of zeros which accumulate on the boundary of the
corresponding region of meromorphy. In that case, for the time
being, $Q_{{\bf n},j}$ denotes a formal infinite product.

The next proposition is adapted from \cite{lomi}.

\begin{prop}\label{ortogonalidad}  Let ${\bf n} \in\mathbb{Z}_+^{m_1+1}(\bullet)\times
\mathbb{Z}_+^{m_2+1}(\bullet), n_{1,m_1}\geq 1,$ $S^1=
\mathcal{N}(\sigma_0^1,\ldots,\,\sigma^1_{m_1}),$ and $ S^2
=\mathcal{N}(\sigma^2_0,\ldots,\,\sigma^2_{m_2})\,\, $ be given
$($recall that $\sigma^1_0 = \sigma^2_0)$. Then, $\deg Q_{{\bf
n},j} = N_{1,j} -1, j=0,1\ldots,m_1,$ all its zeros are simple and
lie in the interior of $\Delta^1_j$. If $I$ denotes the closure of
a connected component of $\Delta^1_j \setminus \supp{\sigma^1_j}$
then $Q_{{\bf n},j}$ has at most one zero in $I$. Moreover,
\begin{equation} \label{orto1}
\int x^{\nu} \mathcal{A}_{{\bf n},j}(x)
\frac{d\sigma^1_j(x)}{Q_{{\bf n},j-1}(x)} = 0\,, \qquad
\nu=0,\ldots,N_{1,j}-2, \qquad j=1,\ldots,m_1,
\end{equation}
and for any polynomial $q, \deg q \leq N_{1,j+1} -1,$
\begin{equation} \label{formint}
\frac{q(z)\mathcal{A}_{{\bf n},j}(z)}{Q_{{\bf n},j}(z)} = \int
\frac{q(x)\mathcal{A}_{{\bf n},j+1}(x)}{Q_{{\bf n},j}(x)}
\frac{d\sigma^1_{j+1}(x)}{z-x}, \qquad j=0,\ldots,m_1-1.
\end{equation}
\end{prop}

\begin{proof} Using induction on $j$, we will prove simultaneously the general
statement concerning the zeros and (\ref{orto1}). Then, we prove
that on any interval $I$ there is at most one zero of $Q_{{\bf
n},j}$. Finally, we obtain (\ref{formint}). For $j=0$, we already
proved that $\mathcal{A}_{{\bf n},0}$ has $N_{1,0} -1 = |{\bf
n}_1|-1$ sign changes in the interior of $\Delta^1_0 =
\Delta^2_0$. Therefore, $\deg Q_{{\bf n},0} \geq N_{1,0} -1$. If
$\deg Q_{{\bf n},0} = N_{1,0} -1$ we conclude with the initial
step.

Suppose that $\deg Q_{{\bf n},0} \geq N_{1,0}$ (including the
possible case that $\deg Q_{{\bf n},0} = \infty$). It is easy to
see that $\mathcal{A}_{{\bf n},0}(\overline{z})=
\overline{\mathcal{A}_{{\bf n},0}({z})},$ so the zeros of $Q_{{\bf
n},0}$ come in conjugate pairs. Therefore, we can choose $N_{1,0}$
(or $N_{1,0}+1$ if necessary) zeros of $Q_{{\bf n},0}$ in such a
way that the monic polynomial $Q_{{\bf n},0}^*$ with this set of
zeros has constant sign on $\Delta^1_1\,\, (\Delta^1_1 \cap
\Delta^1_0 = \emptyset)$. Notice that
\[ \frac{\mathcal{A}_{{\bf n},0}}{Q_{{\bf n},0}^*} \in
\mathcal{H}(\overline{\mathbb{C}}\setminus \Delta^1_1)
\]
is analytic in the indicated region and
\[ \frac{z^{\nu}\mathcal{A}_{{\bf n},0}}{Q_{{\bf n},0}^*} =
\O\left(\frac{1}{z^2}\right)\,, \qquad \nu = 0,\ldots,N_{1,1}
-1\,.
\]
In  this paper the symbol $\mathcal{O}(\cdot)$ always refers to $z
\to \infty$. From (\ref{eq:f}),  we get
\[
0= \int x^{\nu}\mathcal{A}_{{\bf
n},1}(x)\frac{d\sigma^1_1(x)}{Q_{{\bf n},0}^*(x)}, \qquad
\nu=0,\ldots,N_{1,1}-1.
\]
This implies that $\mathcal{A}_{{\bf n},1}$ has at least $N_{1,1}$
zeros on $\Delta^1_1$. According to Lemma \ref{ulises} this linear
form can only have $N_{1,1} -1$ zeros on this interval.
Consequently, our initial assumption is false and $\deg Q_{{\bf
n},0} = N_{1,0} -1$.

Suppose that we have proved that for some $j  \in
\{0,\ldots,m_1-1\},$ $\deg Q_{{\bf n},j} = N_{1,j} -1$, all its
zeros are simple and lie in the interior of $\Delta^1_j$. Let us
show that then, (\ref{orto1}) and the statement concerning the
zeros are valid for $j+1$.

Indeed, the induction hypothesis implies that
\[ \frac{\mathcal{A}_{{\bf n},j}}{Q_{{\bf n},j} } \in
\mathcal{H}(\overline{\mathbb{C}}\setminus \Delta_{j+1}^1) \,,
\qquad \frac{z^{\nu}\mathcal{A}_{{\bf n},j}}{Q_{{\bf n},j}} =
\O\left(\frac{1}{z^2}\right)\,, \qquad \nu = 0,\ldots,N_{1,j+1}
-2\,.
\]
From (\ref{eq:f}), it follows that
\[
0= \int x^{\nu}\mathcal{A}_{{\bf
n},j+1}(x)\frac{d\sigma^1_{j+1}(x)}{Q_{{\bf n},j}(x)}, \qquad
\nu=0,\ldots,N_{1,j+1}-2.
\]
We have obtained (\ref{orto1}) for $j+1$.

Formula (\ref{orto1}) for $j+1$ implies that $Q_{{\bf n},j+1}$ has
at least $N_{1,j+1}-1$ sign changes in the interior of
$\Delta^1_{j+1}.$ If $\deg Q_{{\bf n},j+1} = N_{1,j+1}-1,$ we have
finished the proof (for example, this is the case when $j+1 = m_1$
because $\mathcal{A}_{{\bf n},m_1} \equiv a_{{\bf n},m_1}$). Let
us suppose that $\deg Q_{{\bf n},j+1} \geq N_{1,j+1}$ (including
the possible case that $\deg Q_{{\bf n},j+1} = \infty$, and of
course $j \leq m_1-2$). Since $\mathcal{A}_{{\bf
n},j+1}(\overline{z})= \overline{\mathcal{A}_{{\bf n},j+1}({z})}$,
we can choose $N_{1,j+1}$ (or $N_{1,j+1}+1$ if necessary) zeros of
$Q_{{\bf n},j+1}$ so that the monic polynomial $Q_{{\bf n},j+1}^*$
with this set of zeros has constant sign on $\Delta^1_{j+2}$. Then
\[ \frac{\mathcal{A}_{{\bf n},j+1}}{Q_{{\bf n},j+1}^*} \in
\mathcal{H}(\overline{\mathbb{C}}\setminus \Delta^1_{j+2}) \,,
\qquad \frac{z^{\nu}\mathcal{A}_{{\bf n},j+1}}{Q_{{\bf n},j+1}^*}
= \O\left(\frac{1}{z^2}\right)\,, \qquad \nu = 0,\ldots,N_{1,j+2}
-1\,.
\]
Using (\ref{eq:f}), it follows that
\[
0=\int x^{\nu}\mathcal{A}_{{\bf
n},j+2}(x)\frac{d\sigma^1_{j+2}(x)}{Q_{{\bf n},j+1}^*(x)}, \qquad
\nu=0,\ldots,N_{1,j+2}-1\,.
\]
This implies that $\mathcal{A}_{{\bf n},j+2}$ has at least
$N_{1,j+2}$ zeros on $\Delta^1_{j+2}$. According to Lemma
\ref{ulises} this linear form can only have $N_{1,j+2} -1$ zeros
on this interval. This implies that our initial assumption is
false; therefore, $\deg Q_{{\bf n},j+1} = N_{1,j+1} -1$ as stated.

Suppose that the interval $I$ contains two zeros $x_1,x_2$ of
$Q_{{\bf n},j}$; that is, of $\lnj$. According to (\ref{orto1})
\[
\int x^{\nu} \frac{\mathcal{A}_{{\bf n},j}(x)}{(x-x_1)(x-x_2)}
\frac{(x-x_1)(x-x_2)d\sigma^1_j(x)}{Q_{{\bf n},j-1}(x)} = 0\,,
\qquad \nu=0,\ldots,N_{1,j}-2.
\]
The function ${\mathcal{A}_{{\bf n},j}(x)}/{(x-x_1)(x-x_2)}$ has
$N_{1,j}-3$ sign changes on $\sop(\sigma^1_j)$ while the measure
${(x-x_1)(x-x_2)d\sigma^1_j(x)}/{Q_{{\bf n},j-1}(x)}$ has constant
sign on $\sop(\sigma^1_j)$. This is impossible because of the
number of orthogonality relations.

Formula (\ref{formint}) follows from (\ref{eq:g}) since  for any
$q, \deg q \leq N_{1,j+1}-1,$
\[ \frac{q \mathcal{A}_{{\bf n},j}}{Q_{{\bf n},j}} \in
\mathcal{H}(\overline{\mathbb{C}} \setminus \Delta^1_{j+1})\,,
\qquad \frac{q \mathcal{A}_{{\bf n},j}}{Q_{{\bf n},j}} =
\mathcal{O}\left(\frac{1}{z}\right).
\]
With this we conclude the proof. \end{proof}

We need to produce additional orthogonality relations. In the
second part of this section, we make use of some transformations
employed in \cite{GRS}. Let us define recursively the following
functions
\begin{equation}\label{defnlnjnegative}
\mathcal{A}_{{\bf n},-j-1}(z) = \int \frac{\mathcal{A}_{{\bf
n},-j}(x)}{z-x} d\sigma^2_j(x), \qquad j =0,\ldots,m_2.
\end{equation}

\begin{prop}\label{ortoga}  Let ${\bf n}  \in\mathbb{Z}_+^{m_1+1}(\bullet)\times
\mathbb{Z}_+^{m_2+1}(\bullet), n_{1,m_1}\geq 1,$ $S^1=
\mathcal{N}(\sigma^1_0,\ldots,\sigma^1_{m_1}),$ and $
S^2=\mathcal{N}(\sigma^2_0,\ldots,\,\sigma^2_{m_2})$ be given.
Then, for each $j=0,\ldots,m_2$
\begin{equation} \label{eq:10}
\int x^{\nu} \mathcal{A}_{{\bf n},-j}(x)d s^2_{j,k}(x)=0, \qquad
k=j,\ldots,m_2, \qquad \nu =0,\ldots,n_{2,k} -1.
\end{equation}
\end{prop}

\begin{proof}
When $j=0$ the statement reduces to the relations ii) which define
$\mathcal{A}_{{\bf n},0}.$ If $m_2=0$ we are done. Therefore, let
us assume that $m_2 \geq 1,$ that (\ref{eq:10}) holds for some $j
\in \{0,\ldots,m_2-1\},$ and prove that it is also satisfied for
$j+1$.

Fix $j \in \{0,\ldots,m_2-1\}, k \in \{j+1,\ldots,m_2\},$ and $\nu
\in \{0,\ldots,n_{2,k} -1\}$. Using the definition of
$\mathcal{A}_{{\bf n},-j-1}$, Fubini's theorem, and the induction
hypothesis, we obtain
\[ \int x^{\nu} \mathcal{A}_{{\bf
n},-j-1}(x)d s^2_{j+1,k}(x)
 = \int x^{\nu} \int \frac{\mathcal{A}_{{\bf
n},-j}(t)}{x-t}d\sigma^2_{j}(t) d s^2_{j+1,k}(x) =
\]
\[ \int
\mathcal{A}_{{\bf n},-j}(t)\int \frac{x^{\nu} \mp t^{\nu} }{x-t} d
s^2_{j+1,k}(x) d\sigma^2_{j}(t) =
\]
\[\int p_{\nu}(t) \mathcal{A}_{{\bf
n},-j}(t)d \sigma^2_{j}(t) - \int t^{\nu} \mathcal{A}_{{\bf
n},-j}(t)d s^2_{j,k}(t) = 0
\]
since $p_{\nu}$ is a polynomial of degree $\leq n_{2,k} -2,$ and
$n_{2,j+1} \geq n_{2,k}.$
\end{proof}

For $j= 1,\ldots,m_2+1$, let $Q_{{\bf n},-j}$ be the monic
polynomial whose zeros are those of $\mathcal{A}_{{\bf n},-j}$ in
the region $ {\mathbb{C}} \setminus \Delta^2_{j-1}$ counting
multiplicities.   As we did before, in the hypothetical case that
$\mathcal{A}_{{\bf n},-j}$ had infinitely many zeros in the
specified region, then $Q_{{\bf n},-j}$ denotes a formal infinite
product.

Taking linear combinations of the relations (\ref{eq:10}), we
obtain
\[
\int {\mathcal{B}}_{{\bf n}_2,j}(x) \mathcal{A}_{{\bf n},-j}(x)d
\sigma^2_{j}(x)=0, \qquad j=0,\ldots,m_2,
\]
where ${\mathcal{B}}_{{\bf n}_2,j}$ is an arbitrary linear form of
type
\[ {\mathcal{B}}_{{\bf n}_2,j}(x) = \sum_{k = j}^{m_2} b_{k}(x)
\widehat{s}^2_{j+1,k}(x), \qquad \deg b_k \leq n_{2,k} -1.
\]
Using Lemma \ref{ulises}, it follows that $\mathcal{A}_{{\bf
n},-j}$ has at least $N_{2,j}$ sign changes on $\Delta^2_{j}$,
where
\[ N_{2,j} = N_{2,j}({\bf n}) = n_{2,j}+\cdots +n_{2,m_2}, \qquad j=1,\ldots,m_2.
\]
Consequently, $\deg Q_{n,-j} \geq N_{2,j}, j=0,\ldots,m_2.$ Recall
that for $j=0$ we proved in Proposition \ref{ortogonalidad} that
$\deg Q_{n,0} = N_{2,0} = |{\bf n}_2| = |{\bf n}_1|-1,$ that its
zeros are simple, and lie in the interior of $\Delta^2_0=
\Delta^1_0$.

\begin{prop}\label{ortogb}  Let ${\bf n}  \in\mathbb{Z}_+^{m_1+1}(\bullet)\times
\mathbb{Z}_+^{m_2+1}(\bullet), n_{1,m_1}\geq 1,$ $S^1=
\mathcal{N}(\sigma^1_0,\ldots,\sigma^1_{m_1}),$ and $
S^2=\mathcal{N}(\sigma^2_0,\ldots,\,\sigma^2_{m_2})$ be given.
Then, $\deg Q_{{\bf n},-j} = N_{2,j}, j=0,\ldots,m_2,$ all its
zeros are simple and lie in the interior of $\Delta^2_{j}$, and
$Q_{{\bf n},-m_2-1} \equiv 1$.  If $I$ denotes the closure of a
connected component of $\Delta^2_j \setminus \supp{\sigma^2_j}$
then $Q_{{\bf n},-j}$ has at most one zero in $I$. Moreover,
\begin{equation} \label{orto2}
\int x^{\nu} \mathcal{A}_{{\bf n},-j}(x) \frac{d
\sigma^2_{j}(x)}{Q_{{\bf n},-j-1}(x)}=0,  \qquad
\nu=0,\ldots,N_{2,j}-1, \qquad j=0,\ldots,m_2,
\end{equation}
and for any polynomial $q, \deg q \leq N_{2,j-1},$
\begin{equation} \label{formint2}
\frac{q(z)\mathcal{A}_{{\bf n},-j}(z)}{Q_{{\bf n},-j}(z)} = \int
\frac{q(x)\mathcal{A}_{{\bf n},-j+1}(x)}{Q_{{\bf n},-j}(x)}
\frac{d\sigma^2_{j-1}(x)}{z-x}, \qquad j=1,\ldots,m_2+1.
\end{equation}
\end{prop}

\begin{proof} Fix $j \in \{0,\ldots,m_2\}$.  From (\ref{eq:10}) we have that for each
 $q, \deg q \leq n_{2,j},$
\[
\int \frac{q(z)-q(x)}{z-x}\mathcal{A}_{{\bf n},-j}(x)d
\sigma^2_{j}(x)=0.
\]
It follows that
\[\mathcal{A}_{{\bf
n},-j-1}(z) = \frac{1}{q(z)}\int \frac{q(x)}{z-x}\mathcal{A}_{{\bf
n},-j}(x)d \sigma^2_{j}(x) = \mathcal{O}\left(1/z^{n_{2,j}
+1}\right), \qquad z \to \infty.
\]

We have shown that $\deg Q_{{\bf n},-j-1} \geq N_{2,j+1}
(N_{2,m_2+1} =0)$. The zeros of $Q_{{\bf n},-j-1}$ come in
conjugate pair since $\mathcal{A}_{{\bf n},-j-1}$ is also
symmetric with respect to the real line. If $\deg Q_{{\bf n},-j-1}
> N_{2,j+1}$ take $N_{2,j+1} +1$ (or $N_{2,j+1} +2$ if necessary)
zeros from $Q_{{\bf n},-j-1}$ so that the monic polynomial
$Q_{{\bf n},-j-1}^*$ with these zeros has constant sign on
$\Delta^2_j$. If $\deg Q_{{\bf n},-j-1} = N_{2,j+1}$ take $Q_{{\bf
n},-j-1}^* = Q_{{\bf n},-j-1}.$

Therefore,
\[ \frac{\mathcal{A}_{{\bf
n},-j-1}}{Q_{{\bf n},-j-1}^*} = \mathcal{O}\left( 1/z^{n_{2,j} +
\deg Q_{{\bf n},-j-1}^* + 1}\right) \in
\mathcal{H}(\overline{\mathbb{C}} \setminus \Delta^2_{j}),
\]
and
\[ \frac{z^{\nu}\mathcal{A}_{{\bf
n},-j-1}}{Q_{{\bf n},-j-1}^*} = \mathcal{O}\left( 1/z^2\right) \in
\mathcal{H}(\overline{\mathbb{C}} \setminus \Delta^2_{j}), \qquad
\nu = 0,\ldots, n_{2,j} + \deg Q_{{\bf n},-j-1}^* -1.
\]
Using (\ref{eq:f}), we obtain
\[ 0 =  \int x^{\nu} \mathcal{A}_{{\bf n},-j}(x) \frac{d \sigma^2_{j}(x)}{Q_{{\bf
n},-j-1}^*(x)}, \qquad \nu = 0,\ldots, n_{2,j} + \deg Q_{{\bf
n},-j-1}^* -1.
\]

This formula implies that $\mathcal{A}_{{\bf n},-j}$ has at least
$n_{2,j} + \deg Q_{{\bf n},-j-1}^* \geq N_{2,j} $ sign changes on
$\Delta^2_{j}$. In particular, we have proved that if for some $j
\in \{0,\ldots,m_2\},$ $ \deg {Q}_{{\bf n},-j-1} >  N_{2,j+1}$
then $\deg {Q}_{{\bf n},-j}> N_{2,j}$. Going downwards on the
index $j$ we would obtain that $\deg {Q}_{{\bf n},0} > N_{2,0} =
|{\bf n}_2|=|{\bf n}_1| -1$ which is false according to
Proposition \ref{ortogonalidad}. Consequently, for all $j \in
\{0,\ldots,m_2\}, \deg {Q}_{{\bf n},-j-1} =  N_{2,j+1}$ (in
particular, $Q_{{\bf n},-m_2-1} \equiv 1$). Hence, $Q_{{\bf
n},-j-1}^* = Q_{{\bf n},-j-1}$ and (\ref{orto2}) follows. The
proof that $I$ contains at most one zero of $Q_{{\bf n},-j}$ is
the same as in Proposition \ref{ortogonalidad}.

Now, fix $j \in \{1,\ldots,m_2+1\}$. Notice that for any $q, \deg
q \leq N_{2,j-1},$
\[ \frac{q \mathcal{A}_{{\bf n},-j}}{Q_{{\bf n},-j}} \in
\mathcal{H}(\overline{\mathbb{C}} \setminus \Delta^2_{j-1})\,,
\quad \frac{q \mathcal{A}_{{\bf n},-j}}{Q_{{\bf n},-j}} =
\mathcal{O}\left(\frac{1}{z}\right)\,, z \to \infty\,.
\]
Using (\ref{eq:g}),  (\ref{formint2}) readily follows.
\end{proof}

\section{Interlacing properties} \label{intprop}
Fix a vector $l:=(l_{1};l_{2})$ where $0\leq l_{1}\leq m_{1}$ and
$0\leq l_{2}\leq m_{2}$. We define the multi-index ${\bf
n}^{l}:=({\bf n}_{1}+{\bf e}^{l_{1}};{\bf n}_{2}+{\bf
e}^{l_{2}})=({\bf n}_{1}^{l_{1}};{\bf n}_{2}^{l_{2}})$, where
${\bf e}^{l_{i}}$ denotes the unit vector of length $m_{i}+1$ with
all components equal to zero except the component $(l_{i}+1)$
which equals $1$. In this section it is always assumed that both
$\bf n$ and ${\bf n}^{l}$ belong to
$\mathbb{Z}_{+}^{m_{1}+1}(\bullet)\times\mathbb{Z}_{+}^{m_{2}+1}(\bullet)$.

Fix real constants $A,B$ such that $|A|+|B|>0$ and define
\[ {\mathcal{G}}_{{\bf n},j}:=A\lnj+B\lnlj, \qquad j=0,\ldots,m_1. \]
Since $\deg a_{{\bf n}^{l},l_1} = \deg a_{{\bf n},l_1} +1$ it is
obvious that ${\mathcal{G}}_{{\bf n},j} \not \equiv 0, j\leq l_1$.
In particular, this is always true for ${\mathcal{G}}_{{\bf
n},0}$.

\begin{lem}\label{lem:a}
Assume that $A, B \in {\mathbb{R}}, |A|+|B|>0,$ and
$n_{1,m_{1}}\geq 1$.  Then for all $j\in \{0,\ldots,m_1\}$ such
that $n_{1,j} \geq 2$, $\deg Aa_{{\bf n},j} + Ba_{{\bf n}^l,j}
\geq n_{1,j} - 2$ and ${\mathcal{G}}_{{\bf n},j} \not \equiv 0$.
\end{lem}

\begin{proof} Assume that there exists $j \in \{0,\ldots,m_1\}$ such
that $n_{1,j} \geq 2$ and $\deg Aa_{{\bf n},j} + Ba_{{\bf n}^l,j}
\leq n_{1,j} -3$ ($n_{1,j} -3 = -1$ means that $Aa_{{\bf n},j} +
Ba_{{\bf n}^l,j} \equiv 0$). Then ${\bf n}^{l_1}_1 - 2{\bf e}^j
\in \mathbb{Z}_+^{m_1+1}(*)$, where ${\bf e}^j$ denotes the
$m_1+1$ dimensional unit vector with all components equal to zero
except the component $j+1$ which equals $1$. According to Lemma
\ref{ulises} the linear form ${\mathcal{G}}_{{\bf n},0}$ has at
most $|{\bf n}_1| -2$ zeros on $\Delta^1_0$, but
${\mathcal{G}}_{{\bf n},0}$ satisfies the same orthogonality
relations (\ref{definition}) as $\mathcal{A}_{{\bf n},0}$ and,
therefore, it has at least $|{\bf n}_1|-1$ sign changes on this
interval.  This contradiction implies the statement.
\end{proof}

From this lemma it follows that if $n_{1,m_1} \geq 2$ then
${\mathcal{G}}_{{\bf n},j} \not \equiv 0, j\in \{0,\ldots,m_1\}.$

\begin{lem}\label{leminic2}
Assume that ${A,B}\in \mathbb{R}$ and $\gnj=A\lnj+B\lnlj\not\equiv
0$, for some $j\in\{0,\ldots,m_{1}\}$. If $j\leq l_{1}$ then
$\gnj$ has at most $N_{1,j}$ zeros, counting multiplicities, on
any interval disjoint from
$\Delta_{j+1}^{1}(\Delta_{m_{1}+1}^{1}=\emptyset)$. If $j>l_{1}$
then $\gnj$ has at most $N_{1,j}-1$ zeros, counting
multiplicities, on any interval disjoint from $\Delta_{j+1}^{1}$.
\end{lem}
\begin{proof} We have
\[
\gnj(z)=\sum_{k=j}^{m_{1}}(Aa_{{\bf n},k}(z)+B a_{{\bf
n}^l,k}(z))\widehat{s}_{j+1,k}^{1}(z),
\]
where $\deg a_{{\bf n},k}=n_{1,k}-1$ and $\deg a_{{\bf
n}^{l},k}=n_{1,k}^{l_{1}}-1$. By Lemma \ref{ulises},
$(1,\widehat{s}_{j+1,j+1}^{1},\ldots,\widehat{s}_{j+1,m_{1}}^{1})$
forms an AT-system with respect to
$(n^{l_1}_{1,j},\ldots,n^{l_1}_{1,m_1})$ on any interval disjoint
from $\Delta_{j+1}^{1}$, and the result follows immediately.
\end{proof}

Notice that for each $j\in\{0,\ldots,m_{1}\}$, $\gnj$ is a real
function when it is restricted to the real line.

\begin{prop}\label{theoimp} Let $n_{1,m_{1}}\geq 1$. Assume that
${A,B}\in\mathbb{R}, |A| +|B|
>0,$ and let $k = \max\{k': \mathcal{G}_{{\bf n},k'} \not\equiv 0\} \leq m_1$. Then, $k \geq
l_1$ and $\gnj \equiv 0, k < j \leq m_1$. If $j\leq l_{1}$ then
$\gnj$ has at most $N_{1,j}$ zeros in
$\mathbb{C}\setminus\Delta_{j+1}^{1}$, counting multiplicities,
and at least $N_{1,j}-1$ sign changes in the interior of
$\Delta_{j}^{1}$. If $l_1 < j \leq k$ then $\gnj$ has at most
$N_{1,j}-1$ zeros in $\mathbb{C}\setminus\Delta_{j+1}^{1}$ and at
least $N_{1,j}-2$ sign changes in the interior of
$\Delta_{j}^{1}$. Therefore, all the zeros of $\gnj$ in
$\mathbb{C}\setminus\Delta_{j+1}^{1}$ are real and simple.
\end{prop}
\begin{proof} If $j\leq l_{1}$, then $\deg a_{{\bf n}^l,l_1} > \deg
a_{{\bf n},l_1}$ and $\gnj\not\equiv 0$. Consequently, $k \geq
l_1$. Obviously, from the definition of $k$, $\gnj \equiv 0, k < j
\leq m_1$.

Assume that $\gnj, j \leq l_1,$ has at least $N_{1,j}+1$ zeros in
$\mathbb{C}\setminus\Delta_{j+1}^{1}$, counting multiplicities.
Select $N_{1,j}+1$ or $N_{1,j}+2$ zeros of $\gnj$ which are
symmetric with respect to the real axis, and let $Q_{{\bf
n},j}^{*}$ be the monic polynomial whose zeros are those
prescribed. If $j<l_{1}$ then
\[
\frac{z^{\nu}\gnj}{Q_{{\bf
n},j}^{*}}=\mathcal{O}\Big(\frac{1}{z^2}\Big)\,, \qquad
\nu=0,\ldots,N_{1,j+1}\,.
\]
From (\ref{eq:f}),  it follows that
\[
0=\int x^{\nu}{\mathcal G}_{{\bf
n},j+1}(x)\frac{d\sigma_{j+1}^{1}(x)}{Q_{{\bf
n},j}^{*}(x)}\,,\qquad \nu=0,\ldots,N_{1,j+1}\,.
\]
These orthogonality relations imply that $\gnju$ has at least
$N_{1,j+1}+1$ zeros on $\Delta_{j+1}^{1}$. Since $\gnju\not\equiv
0$ we obtain a contradiction with Lemma \ref{leminic2}.

If $j=l_{1}$ and $j<k$, then
\[
\frac{z^{\nu}{\mathcal{G}}_{{\bf n},l_{1}}}{Q_{{\bf
n},l_{1}}^{*}}=\mathcal{O}\Big(\frac{1}{z^2}\Big)\,, \qquad
\nu=0,\ldots,N_{1,l_{1}+1}-1\,.
\]
Arguing as before, it follows that ${\mathcal{G}}_{{\bf
n},l_{1}+1}$ has at least $N_{1,l_{1}+1}$ zeros on
$\Delta_{l_{1}+1}^{1}$, contradicting Lemma \ref{leminic2}. If
$j=l_{1} = k$ then ${\mathcal{G}}_{{\bf n},l_{1}+1}\equiv 0$ and
${\mathcal{G}}_{{\bf n},l_{1}} = Aa_{{\bf n},l_1} +Ba_{{\bf
n}^l,l_1}$ is a polynomial of degree at most ${ n}_{1,l_1} <
N_{1,l_1}+1$ and thus it is identically equal to zero which is
impossible. Consequently, when $j \leq l_1$, $\gnj$ has at most
$N_{1,j}$ zeros in $\mathbb{C}\setminus\Delta_{j+1}^{1}$ counting
multiplicities.

Let $l_{1} < j \leq k$ and assume that $\gnj$ has at least
$N_{1,j}$ zeros in $\mathbb{C}\setminus\Delta_{j+1}^{1}$, counting
multiplicities. If $j=m_{1}$ we get immediately a contradiction
because in this case ${\mathcal{G}}_{{\bf n},m_{1}}$ is a
polynomial of degree at most $N_{1,m_{1}}-1$. If $l_{1}<j<m_{1}$,
then there exists a polynomial $Q_{{\bf n},j}^{*}$ with real
coefficients and degree at least $N_{1,j}$ such that
\[
\frac{z^{\nu}{\mathcal{G}}_{{\bf n},j}}{Q_{{\bf
n},j}^{*}}=\mathcal{O}\Big(\frac{1}{z^2}\Big)\,, \qquad
\nu=0,\ldots,N_{1,j+1}-1\,.
\]
This implies that $\gnju$ has at least $N_{1,j+1}$ zeros on
$\Delta_{j+1}^{1}$ contradicting Lemma \ref{leminic2}.

Now, let us analyze the sign changes.  Notice that
${\mathcal{G}}_{{\bf n},0} \not\equiv 0$. Assume that
${\mathcal{G}}_{{\bf n},0}$ has $N<N_{1,0}-1=|{\bf n}_{1}|-1$ sign
changes on $\Delta_{0}^{1}=\Delta_{0}^{2}$, choose a nonzero
linear form
\[
{\mathcal{B}}_{{\bf n}_{2}}(z)=\sum_{k=0}^{m_{2}} b_{{\bf
n}_{2},k}(z)\widehat{s}_{1,k}^{2}(z)\,,\quad \deg b_{{\bf
n}_{2},k}\leq n_{2,k}-1\,,\quad k=0,\ldots,m_{2}\,,
\]
such that ${\mathcal{B}}_{{\bf n}_{2}}$ has a zero at each point
where ${\mathcal{G}}_{{\bf n},0}$ has a sign change, and a zero of
order $|{\bf n}_{2}|-1-N$ at one of the extreme points of
$\Delta_{0}^{2}$. By Lemma \ref{ulises}, ${\mathcal{B}}_{{\bf
n}_{2}}$ has at most $|{\bf n}_{2}|-1$ zeros on $\Delta_{0}^{2}$.
Thus, ${\mathcal{B}}_{{\bf n}_{2}}$ has exactly those zeros
prescribed. By definition,
\[
\int {\mathcal{B}}_{{\bf n}_{2}}(x){\mathcal{G}}_{{\bf
n},0}(x)d\sigma_{0}^{2}(x)=0\,,
\]
which contradicts the fact that ${\mathcal{B}}_{{\bf
n}_{2}}(x){\mathcal{G}}_{{\bf n},0}(x)$ has constant sign on
$\Delta_{0}^{2}$.

Let us prove by induction that for all $j \leq l_1$,
${\mathcal{G}}_{{\bf n},j}$ has at least $N_{1,j}-1$ sign changes
in the interior of $\Delta^1_j$. For $j=0$ this was proved above
and if $l_1=0$ we are done. Let us assume that for some $j < l_1$,
${\mathcal{G}}_{{\bf n},j}$ has at least $N_{1,j}-1$ sign changes
on $\Delta^1_j$, and let us show that ${\mathcal{G}}_{{\bf
n},j+1}$ has at least $N_{1,j+1}-1$ sign changes on
$\Delta^1_{j+1}$.

Let $Q_{{\bf n},j}^{*}$ be a monic polynomial whose zeros are
$N_{1,j}-1$ points where ${\mathcal{G}}_{{\bf n},j}$ has a sign
change. Then
\[
\frac{z^{\nu}{\mathcal{G}}_{{\bf n},j}}{Q_{{\bf
n},j}^{*}}=\mathcal{O}\Big(\frac{1}{z^2}\Big)\,, \qquad
\nu=0,\ldots,N_{1,j+1}-2\,.
\]
Using (\ref{eq:f}), this implies that
\[
0=\int x^{\nu}\gnju(x)\frac{d\sigma_{j+1}^{1}(x)}{Q_{{\bf
n},j}^{*}(x)}\,,\qquad \nu=0,\ldots,N_{1,j+1}-2\,.
\]
Thus, $\mathcal{G}_{{\bf n},j+1}$ has at least $N_{1,j+1}-1$ sign
changes in the interior of $\Delta_{j+1}^{1}$ as claimed.

Finally, we prove that ${\mathcal{G}}_{{\bf n},j}, l_1 < j  \leq
k, $ has at least $N_{1,j}-2$ sign changes in the interior of
$\Delta_{j}^{1}$. Let $Q_{{\bf n},l_{1}}^{*}$ be a monic
polynomial of degree $N_{1,l_{1}}-1$ whose zeros are points where
${\mathcal{G}}_{{\bf n},l_{1}}$ changes sign in the interior of
$\Delta_{l_{1}}^{1}$, then
\[
\frac{z^{\nu}{\mathcal{G}}_{{\bf n},l_{1}}}{Q_{{\bf
n},l_{1}}^{*}}=\mathcal{O}\Big(\frac{1}{z^2}\Big)\,, \qquad
\nu=0,\ldots,N_{1,l_{1}+1}-3\,.
\]
From here we get orthogonality conditions that imply that
${\mathcal{G}}_{{\bf n},l_{1}+1}$ has at least $N_{1,l_{1}+1}-2$
sign changes in the interior of $\Delta_{l_{1}+1}^{1}$. One
proceeds the same way until we arrive to $j=k$.

From the upper bound on the number of zeros and the lower bound on
the number of sign changes it follows that all the zeros are
simple and lie on the real line. \end{proof}

Let $j\in\{0,\ldots,m_{2}+1\}$. Given two real constants $A,B,$ we
define
\[\gnmj:=A\lnmj+B\lnlmj. \]
Thus, by (\ref{defnlnjnegative}),
\begin{equation}\label{eq:11}
\gnmju(z) = \int \frac{\gnmj(x)}{z-x} d\sigma^2_j(x), \qquad j
=0,\ldots,m_2.
\end{equation}
If $|A|+|B|>0$ then ${\mathcal G}_{{\bf n},0}\not\equiv 0$ and
from (\ref{eq:11}) it follows that $\gnmj\not\equiv 0$ for all
$j\in\{1,\ldots,m_{2}+1\}$.

\begin{prop}\label{theoimp2}
Let $A,B\in {\mathbb{R}}, |A|+|B|>0$. For every
$j\in\{1,\ldots,m_{2}\}$, $\gnmj$ has at most $N_{2,j}+1$ zeros on
$\mathbb{C}\setminus\Delta_{j-1}^{2}$, counting multiplicities,
and at least $N_{2,j}$ sign changes in the interior of
$\Delta_{j}^{2}$. Hence, all the zeros of $\gnmj$ on
$\mathbb{C}\setminus\Delta_{j-1}^{2}$ are real and simple.
\end{prop}
\begin{proof} Let $j\in\{0,\ldots,m_{2}\}$. By (\ref{eq:10}) we
know that
\[
\int x^{\nu} \mathcal{A}_{{\bf n}^{l},-j}(x)d s^2_{j,k}(x)=0\,,
\qquad k=j,\ldots,m_2\,, \qquad \nu =0,\ldots,n_{2,k}^{l_{2}}-1\,.
\]
Since $n_{2,k}\leq n_{2,k}^{l_{2}}$, it follows that
\begin{equation}\label{eq:12}
\int x^{\nu} \gnmj(x) d s^2_{j,k}(x)=0\,, \qquad k=j,\ldots,m_2\,,
\qquad \nu =0,\ldots,n_{2,k}-1\,.
\end{equation}
Using the same arguments employed in the previous section to show
that $\lnmj$ has at least $N_{2,j}$ sign changes in the interior
of $\Delta_{j}^{2}$, one obtains the same conclusion for $\gnmj$.

If $q$ is a polynomial with $\deg q\leq n_{2,j}$, then from
(\ref{eq:12}) we have
\[
\int \frac{q(z)-q(x)}{z-x}\gnmj(x)d\sigma_{j}^{2}(x)=0\,.
\]
Hence, for every $j\in\{0,\ldots,m_{2}\}$,
\[
\gnmju(z)=\frac{1}{q(z)}\int\frac{q(x)}{z-x}\gnmj(x)d\sigma_{j}^{2}(x)=
\mathcal{O}\Big(\frac{1}{z^{n_{2,j}+1}}\Big)\,,\quad
z\rightarrow\infty\,.
\]
Assume that for some $j\in\{0,\ldots,m_{2}-1\}$, $\gnmju$ has at
least $N_{2,j+1}+2$ zeros, counting multiplicities, on
$\mathbb{C}\setminus\Delta_{j}^{2}$. Select at least $N_{2,j+1}+2$
zeros of $\gnmju$, symmetric with respect to the real axis, and
denote by $Q_{{\bf n},-j-1}^{*}$ the monic polynomial whose zeros
are the points selected. Then,
\[
\frac{z^{\nu}\gnmju}{Q_{{\bf
n},-j-1}^{*}}=\mathcal{O}\Big(\frac{1}{z^2}\Big)\,, \qquad
\nu=0,\ldots,N_{2,j}+1\,.
\]
As before, this implies that $\gnmj$ has at least $N_{2,j}+2$
zeros in the interior of $\Delta_{j}^{2}$. Going downwards on the
index $j,$ we obtain that ${\mathcal{G}}_{{\bf n},0}$ has at least
$N_{2,0}+2=N_{1,0}+1$ zeros, which is impossible by Proposition
\ref{theoimp}. Therefore, for all $j \in \{1,\ldots,m_2+1\},$
${\mathcal{G}}_{{\bf n},j}$ has at most $N_{2,j}+1$ zeros in
$\mathbb{C}\setminus \Delta^2_{j-1}$ and, therefore, they must be
real and simple.\end{proof}

\begin{thm}\label{theointerlace}
Let ${\bf n}, {\bf n}^{l}\in
\mathbb{Z}_{+}^{m_{1}+1}(\bullet)\times\mathbb{Z}_{+}^{m_{2}+1}(\bullet),
n_{1,m_1} \geq 2$. Then, for all $j \in \{-m_2, \ldots,m_1\}$  the
zeros of $\mathcal{A}_{{\bf n},j}$ and $\mathcal{A}_{{\bf
n}^{l},j}$ interlace; that is, between two consecutive zeros of
$\mathcal{A}_{{\bf n},j}$ there is one zero of $\mathcal{A}_{{\bf
n}^{l},j}$ and viceversa.
\end{thm}

\begin{proof} Since $n_{1,m_1} \geq 2$, from Lemma \ref{lem:a} we
know that for all $j\in \{0,\ldots,m_1\}$ and for all $A,B$ real
such that $|A| +|B| > 0,$ the linear form ${\mathcal{G}}_{{\bf
n},j}$ is not identically equal to zero. This is always true for
$j\in \{-m_2,\ldots,-1\}$. Therefore, from Propositions
\ref{theoimp} and \ref{theoimp2} we know that for all real $A,B,$
such that $|A|+|B|>0$ the zeros of ${\mathcal{G}}_{{\bf n},j},
j\in \{-m_2,\ldots,m_1\},$ are real an simple. This is the basic
fact we will use in the proof.

Fix $y\in \mathbb{R} \setminus \Delta_{j+1}^{1}$. It cannot occur
that $\mathcal{A}_{{\bf n},j}(y)=\mathcal{A}_{{\bf n}^l,j}(y)=0$.
If so, $y$ would be a simple zero of $\mathcal{A}_{{\bf n},j}$ and
$\mathcal{A}_{{\bf n}^l,j}$. Thus, $\mathcal{A}_{{\bf
n},j}'(y)\neq 0$ and $\mathcal{A}_{{\bf n}^l,j}'(y)\neq 0$. Take
$A=1$ and $B=-\mathcal{A}_{{\bf n},j}'(y)/\mathcal{A}_{{\bf
n}^l,j}'(y)$ and consider $\gnj=A\lnj+B\lnlj$. With this choice of
$A$ and $B$, we have
\[
\gnj(y) = \gnj'(y)=0\,,
\]
and we obtain a contradiction because the zeros of
${\mathcal{G}}_{{\bf n},j}$ are simple.

Now, taking $A=\mathcal{A}_{{\bf n}^l,j}(y)$ and
$B=-\mathcal{A}_{{\bf n},j}(y)$, we have that $|A|+|B|>0$. Since
\[
\mathcal{A}_{{\bf n}^l,j}(y)\lnj(y)-\lnj(y)\lnlj(y)=0\,,
\]
and the zeros on $\mathbb{R}\setminus\Delta_{j+1}^{1}$ of
$\mathcal{A}_{{\bf n}^l,j}(y)\lnj(x)-\lnj(y)\lnlj(x)$ with respect
to $x$ are simple, it follows that
\[
\lnlj(y)\lnj'(y)-\lnj(y)\lnlj'(y)\neq 0\,.
\]
But $\lnlj(y)\lnj'(y)-\lnj(y)\lnlj'(y)$ is a continuous real
function on $\mathbb{R}\setminus\Delta_{j+1}^{1}$ in $y$ so it
must have constant sign on each one of the connected components of
$\mathbb{R}\setminus\Delta_{j+1}^{1}$. In particular, its sign on
$\Delta_{j}^{1}$ is constant.

Evaluating $\lnlj(y)\lnj'(y)-\lnj(y)\lnlj'(y)$ at two consecutive
zeros of $\lnlj$, since the sign of $\lnlj'$ at these two points
changes, the sign of $\lnj$ must also change. Using Bolzano's
theorem we find that there must be an intermediate zero of $\lnj$.
Analogously, one proves that between two consecutive zeros of
$\lnj$ on $\Delta_{j}^{1}$ there is one of $\lnlj$. Thus, the
interlacing property has been proved. \end{proof}

\section{Asymptotic distribution of zeros}\label{ext-pro}

Let $\{\mu_l\} \subset \mathcal{M}(E)$ be a sequence of positive
measures, where $E$ is a compact subset of the complex plane and
$\mu \in \mathcal{M}(E)$. We write
\[ *\lim_l \mu_l = \mu\,,
\]
if for every continuous function $f \in \mathcal{C}(E)$
\[ \lim_l \int f d \mu_l
 = \int f d\mu\,;\]
that is, when the sequence of measures converges to $\mu$ in the
weak star topology. Given a polynomial $q_l$ of degree $l \geq 1$,
we denote the associated normalized zero counting measure by
\[ \mu_{q_l} = \frac{1}{l} \sum_{q_l(x) = 0} \delta_x \,,\]
where $\delta_x$ is the Dirac measure with mass $1$ at $x$ (in the
sum the zeros are repeated according to their multiplicity).

\begin{lem} \label{lemextremal}
Let $E \subset \mathbb{C}$ be a compact set which is regular with
respect to the Dirichlet problem and $\phi$ a continuous function
on $E$. Then there exists a unique $\overline{\mu} \in
\mathcal{M}_1(E)$ and a constant $w$ such that
\[
V^{\overline{\mu}}(z)+\phi(z) \left\{ \begin{array}{l} \leq
w,\quad z
\in \supp{\overline{\mu}} \,, \\
\geq w, \quad z \in  E\,.
\end{array} \right.
\]
\end{lem}

If the compact set $E$ is not regular with respect to the
Dirichlet problem then  the second part of the statement is true
except on a set $e$ such that $\mbox{cap}(e) =0.$ Theorem I.1.3 in
\cite{ST} contains a proof of this lemma in this context. When $E$
is regular, it is well known that this inequality except for a set
of capacity zero implies the inequality for all points in the set.
$\overline{\mu}$ is called the equilibrium measure in the presence
of the external field $\phi$ on $E$ and $w$ is  the equilibrium
constant.

In order to determine the asymptotic zero distribution of the
polynomials $Q_{{\bf n},j}$ we use the following lemma. Different
versions of it appear in \cite{calo}, \cite{gora},  and
\cite{stto}. In \cite{gora}, it was proved assuming that
$\sop(\sigma)$ is an interval on which $\sigma'
> 0$ a.e. We wish to preserve this more restrictive condition for
stronger results in section \ref{ratio}. Theorem 3.3.3 in
\cite{stto} and Theorem 1 in \cite{calo}, do not cover the type of
external field we consider here. So, we will sketch a proof.

\begin{lem}\label{gonchar-rakhmanov}
Let $\sigma \in \mbox{\bf Reg}$, $\supp{\sigma} \subset
\mathbb{R}$, where $\supp{\sigma}$ is regular with respect to the
Dirichlet problem. Let $\{\phi_l\}, l \in \Lambda \subset
\mathbb{Z}_+,$ be a sequence of positive continuous functions on
$\supp{\sigma}$ such that
\begin{equation} \label{eq:phi}
\lim_{l\in \Lambda}\frac{1}{2l}\log\frac{1}{|\phi_l(x)|}= \phi(x)
> -\infty ,
\end{equation}
uniformly on $\sop(\sigma)$. By  $\{q_l\}, l \in \Lambda,$ denote
a sequence of monic polynomials such that $\deg q_l = l$ and
\begin{equation} \label{eq:ortogo}
\int x^k q_l(x)\phi_l(x)d\sigma(x)=0,\qquad k=0,\ldots, l-1.
\end{equation}
Then
\begin{equation} \label{eq:18}
*\lim_{l \in \Lambda}\mu_{q_l} = \overline{\mu},
\end{equation}
and
\begin{equation} \label{eq:19}
\lim_{l\in \Lambda}\left(\int |q_l(x)|^2\phi_l(x)
d\sigma(x)\right)^{1/{2l}}= e^{-w},
\end{equation}
where $\overline{\mu}$ and $w$ are the  equilibrium measure and
equilibrium constant in the presence of the external field $\phi$
on $\supp{\sigma}$. We also have that
\begin{equation} \label{eq:H}
\lim_{l \in \Lambda} \left(\frac{|q_l(z)|}{\|q_l
\phi_l^{1/2}\|_E}\right)^{1/l} = \exp{(w -
V^{\overline{\mu}}(z))}, \qquad  \mathcal{K} \subset {\mathbb{C}}
\setminus \mbox{\rm Co}(\sop(\sigma)).
 \end{equation}
\end{lem}

\begin{proof} On account of (\ref{eq:phi}) and Lemma
\ref{lemextremal}, it follows that for any $\varepsilon > 0$ there
exists $l_0$ such that for all $l \geq l_0 , l \in \Lambda,$ and
$z \in \sop (\overline{\mu}) \subset \sop(\sigma) =: E$
\[ \frac{1}{l}\log \frac{|p_l(z)|}{\|p_l \phi_l^{1/2}\|_E} \leq
\frac{1}{2l}\log \frac{1}{|\phi_l(z)|} \leq \phi(z) + \varepsilon
\leq w - V^{\overline{\mu}}(z) + \varepsilon,
\]
where $\{p_l\}, l \in \Lambda,$ is any sequence of monic
polynomials such that $\deg p_l =l$ and $\|p_l \phi_l^{1/2}\|_E =
\max_{z \in E} |(p_l\phi_l^{1/2})(z)|.$ Hence,
\[ u_l(z) := V^{\overline{\mu}}(z) + \frac{1}{l}\log \frac{|p_l(z)|}{\|p_l
\phi_l^{1/2}\|_E} \leq w + \varepsilon, \quad z \in \sop
(\overline{\mu}), \quad l \geq l_0.
\]
Since $u_l$ is subharmonic in $\overline{\mathbb{C}} \setminus
\sop(\overline{\mu})$, by the continuity and maximum principles,
we have
\[ u_l(z) \leq w + \varepsilon, \quad z \in \overline{\mathbb{C}}, \quad l \geq l_0.
\]
In particular,
\[ u_l(\infty) = \frac{1}{l} \log \frac{1}{\|p_l
\phi_l^{1/2}\|_E} \leq w + \varepsilon.
\]
The last two relations imply
\begin{equation} \label{eq:20}
\limsup_{l \in \Lambda} \left(\frac{|p_l(z)|}{\|p_l
\phi_l^{1/2}\|_E}\right)^{1/l} \leq \exp{(w -
V^{\overline{\mu}}(z))}, \qquad  \mathcal{K} \subset {\mathbb{C}},
\end{equation}
and
\begin{equation} \label{eq:21}
\liminf_{l \in \Lambda} \|p_l \phi_l^{1/2}\|_E^{1/l} \geq
\exp{(-w)}.
\end{equation}
In particular, these relations hold for the sequence of
polynomials $\{q_l\}, l \in \Lambda$.

Let $t_l$ be the weighted Fekete polynomial of degree $l$ for the
weight $e^{-\phi}$ on $\sop(\sigma)$ and $|\sigma|$ be the total
variation of $\sigma$. From the minimality property in the $L_2$
norm of $q_l$, we have
\[ \|q_l\phi_l^{1/2}\|_2 :=  \left(\int |q_l(x)|^2\phi_l(x) d\sigma(x)\right)^{1/2} \leq
\|t_l\phi_l^{1/2}\|_2 \leq |\sigma|^{1/2}\|t_l\phi_l^{1/2}\|_E
\leq
\]
\[ |\sigma|^{1/2}  \|t_le^{-l\phi}\|_E \|\phi_l^{1/2}e^{l\phi}\|_E.
\]
Then, using (\ref{eq:phi}) and Theorem III.1.9 in \cite{ST}, we
obtain that
\begin{equation} \label{eq:22} \limsup_{l \in \Lambda}
\|q_l\phi_l^{1/2}\|_2^{1/l} \leq e^{-w}.
\end{equation}

Since $\sop(\sigma)$ is regular with respect to the Dirichlet
problem, Theorem 3.2.3 vi) in \cite{stto} yields
\[
\limsup_{l \in \Lambda} \left(\frac{\|q_l
\phi_l^{1/2}\|_E}{\|q_l\phi_l^{1/2}\|_2}\right)^{1/l} \leq 1,
\]
which combined with (\ref{eq:21}) (with $p_l = q_l$) and
(\ref{eq:22}) implies
\begin{equation} \label{eq:23}
\lim_{l \in \Lambda} \left(\frac{\|q_l
\phi_l^{1/2}\|_E}{\|q_l\phi_l^{1/2}\|_2}\right)^{1/l} = 1.
\end{equation}
Thus, we obtain (\ref{eq:19}) since (\ref{eq:21}), (\ref{eq:22}),
and (\ref{eq:23}) give
\begin{equation} \label{eq:24}
\limsup_{l \in \Lambda} \|q_l\phi_l^{1/2}\|_E^{1/l} = \limsup_{l
\in \Lambda} \|q_l\phi_l^{1/2}\|_2^{1/l} = e^{-w}.
\end{equation}

All the zeros of $q_l$ lie in $\mbox{Co}(\sop(\sigma)) \subset
\mathbb{R}$. The unit ball in the weak star topology of measures
is compact. Take any subsequence of indices $\Lambda' \subset
\Lambda$ such that
\[ *\lim_{l \in \Lambda'}\mu_{q_l} = \mu_{\Lambda'}.
\]
Then,
\[ \lim_{l \in \Lambda'} \frac{1}{l} \log  |q_l(z)| = -\lim_{n \in \Lambda'}
 \int \log \frac{1}{|z - x|} \mu_{q_l}(x) =  -V^{\mu_{\Lambda'}}(z),
 \quad {\mathcal{K}} \subset \mathbb{C} \setminus
 \mbox{Co}(\sop(\sigma)).
\]
This, together with (\ref{eq:19}) and (\ref{eq:20}) (applied to
$\{q_l\}, l \in \Lambda'$), implies
\[ (V^{\overline{\mu}} - V^{\mu_{\Lambda'}})(z) \leq 0, \qquad z \in
\overline{\mathbb{C}} \setminus \mbox{Co}(\sop(\sigma)).
\]
Since $V^{\overline{\mu}} - V^{\mu_{\Lambda'}}$ is subharmonic in
$\overline{\mathbb{C}} \setminus  \sop(\overline{\mu})$ and
$(V^{\overline{\mu}} - V^{\mu_{\Lambda'}})(\infty) =0$, from the
maximum principle, it follows that $V^{\overline{\mu}} \equiv
V^{\mu_{\Lambda'}}$ in $\mathbb{C} \setminus
 \mbox{Co}(\sop(\sigma))$ and thus $\mu_{\Lambda'} =
 \overline{\mu}$. Consequently, (\ref{eq:18}) holds. (\ref{eq:18}) and (\ref{eq:19}) imply (\ref{eq:H}).   \end{proof}

Using Lemma \ref{gonchar-rakhmanov}, we can obtain the asymptotic
limit distribution of the zeros of the polynomials $Q_{{\bf
n},j}$, $j=-m_2,\ldots,m_1$. At this point, let us make a slight
change of notation. In the sequel,
\[\Delta_j = \Delta^1_{j}, \qquad \sigma_j = \sigma^1_{j}, \qquad j=
0,1,\ldots,m_1,\]
\[\Delta_j = \Delta^2_{-j}, \qquad \sigma_j = \sigma^2_{-j}, \qquad j=
0,-1,\ldots,-m_2,
\]
and
\[
 N_{{\bf n},j} = \left\{
\begin{array}{ll}
N_{1,j}({\bf n}) -1, & j=0,1\ldots,m_1, \\ N_{2,-j}({\bf n}), & j=
0,-1,\ldots,-m_2.
\end{array}
\right.
\]

According to Propositions \ref{ortogonalidad} and \ref{ortogb},
for all $j=-m_2,\ldots,m_1$ the zeros of $Q_{{\bf n},j}$ are all
simple, lie in the interior of $\Delta_j$, and total $N_{{\bf
n},j}$ points.

\begin{thm}\label{mea-equil}
Let $\Lambda =
\Lambda(p_{1,0},\ldots,p_{1,m_1};p_{2,0},\ldots,p_{2,m_2}) \subset
\mathbb{Z}_+^{m_1+1}(\bullet) \times
\mathbb{Z}_+^{m_2+1}(\bullet), (S^1,S^2)\in\mbox{\bf {Reg}},$ $
S^1= \mathcal{N}(\sigma^1_0,\ldots,\sigma^1_{m_1}),$ and $S^2=
\mathcal{N}(\sigma^2_0,\ldots,\sigma^2_{m_2})$ be given. Then
\begin{equation}\label{conv-Qnj*}
*\lim_{{\bf n}\in\Lambda}\mu_{Q_{{\bf
n},j}}=\overline{\mu}_j,\quad j=-m_2,\ldots,m_1,
\end{equation}
where $\overline{\mu}=\overline{\mu}(\mathcal{C})\in
\mathcal{M}_1$ is the vector equilibrium measure determined by the
matrix $\mathcal{C}$ in $(\ref{matriz})$ on the system of compact
sets $E_j = \sop (\sigma^1_j), j=0,\ldots,m_1, E_j = \sop
(\sigma^2_{-j}), j=-m_2,\ldots,0$. Moreover,
\begin{equation} \label{eq:4*}
\lim_{{\bf n} \in \Lambda} \left(\int \frac{Q_{{\bf
n},j}^2(x)}{|Q_{{\bf n},j-1}(x)|}\frac{|\mathcal{A}_{{\bf
n},j}(x)|}{|Q_{{\bf n},j}(x)|}d|\sigma_j|(x)\right)^{1/2|{\bf
n}_1|} = \exp\left(-\sum_{k=j}^{m_1}
{\omega_k^{\overline{\mu}}}/{P_k}\right)\,,
\end{equation}
where the $\omega_k^{\overline{\mu}}$ denote the corresponding
equilibrium constants.
\end{thm}
\begin{proof} The unit ball in the cone of positive Borel measures
is weak star compact; therefore, it is sufficient to show that
each one of the sequences of measures $\{\mu_{Q_{{\bf n},j}}\}$,
${\bf n}\in\Lambda$, $j=-m_2,\ldots,m_1,$ has only one
accumulation point which coincides with the corresponding
component of the vector measure $\overline{\mu}(\mathcal{C})$. Let
$\Lambda'\subset \Lambda$ be a subsequence of multi-indices such
that for each $j=-m_2,\ldots,m_1$
\[
*\lim_{{\bf n}\in\Lambda'}\mu_{Q_{{\bf n},j}}=\mu_j.
\]
Notice that $\mu_j\in\mathcal{M}_1(E_j)$, $j=-m_2,\ldots,m_1$.
Therefore,
\begin{equation}\label{conv-Qnj}
\lim_{{\bf n}\in\Lambda'}|Q_{{\bf n},j}(z)|^{1/|{\bf
n}_1|}=\exp(-P_j V^{\mu_j}(z)),
\end{equation}
uniformly on compact subsets of $\C\setminus\Delta_j$, where $P_j
= \lim_{{\bf n} \in \Lambda^{\prime}} N_{{\bf n},j}/|{\bf n}_1|$.

Because of the normalization adopted on $a_{{\bf n},m_1}$,
$\mathcal{A}_{{\bf n},m_1} = Q_{{\bf n},m_1}$; consequently, when
$j=m_1$, (\ref{orto1}) takes the form
\[ \int x^{\nu} Q_{{\bf n},m_1}(x) \frac{d|\sigma_{m_1}|(x)}{|Q_{{\bf
n},m_1-1}(x)|} = 0\,, \qquad \nu=0,\ldots,N_{{\bf n},m_1}-1\,.
\]
(By $|\sigma|$ we denote the total variation of the measure
$\sigma$.) According to (\ref{conv-Qnj})
\[ \lim_{{\bf n} \in \Lambda'}\frac{1}{2N_{{\bf n},m_1}}\log|Q_{{\bf n},m_1-1}(x)| =
-\frac{P_{m_1-1}}{2P_{m_1}} V^{\mu_{m_1-1}}(x)\,,
\]
uniformly on $\Delta_{m_1}$. Using Lemma \ref{gonchar-rakhmanov},
it follows that $\mu_{m_1}$ is the unique solution of the extremal
problem
\begin{equation} \label{eq:1}
V^{\mu_{m_1}}(x) - \frac{P_{m_1-1}}{2P_{m_1}}V^{\mu_{m_1-1}}(x)
\left\{
\begin{array}{l} = \omega_{m_1},\quad x
\in \sop( {\mu_{m_1}}) \,, \\
\geq \omega_{m_1}, \quad x \in E_{m_1} \,,
\end{array} \right.
\end{equation}
and
\begin{equation} \label{eq:2}
\lim_{{\bf n} \in \Lambda'} \left(\int \frac{Q_{{\bf
n},m_1}^2(x)}{|Q_{{\bf
n},m_1-1}(x)|}d|\sigma_{m_1}|(x)\right)^{1/2N_{{\bf n},m_1}} =
e^{-\omega_{m_1}}\,.
\end{equation}

Let us show by induction on decreasing values of $j$, that for all
$j\in \{-m_2,\ldots,m_1\}$
\begin{equation} \label{eq:3}
V^{\mu_j}(x) - \frac{P_{j-1}}{2P_j}V^{\mu_{j-1}}(x) -
\frac{P_{j+1}}{2P_j}V^{\mu_{j+1}}(x) +
\frac{P_{j+1}}{P_j}\omega_{j+1} \left\{
\begin{array}{l} = \omega_j,\quad x
\in \sop( {\mu_j}) \,, \\
\geq \omega_j, \quad x \in E_j \,,
\end{array} \right.
\end{equation}
where $P_{-m_2-1} = P_{m_1+1} = 0$, and
\begin{equation} \label{eq:4}
\lim_{{\bf n} \in \Lambda'} \left(\int \frac{Q_{{\bf
n},j}^2(x)}{|Q_{{\bf n},j-1}(x)|}\frac{|\mathcal{A}_{{\bf
n},j}(x)|}{|Q_{{\bf n},j}(x)|}d|\sigma_j|(x)\right)^{1/2N_{{\bf
n},j}} = e^{-\omega_j}\,,
\end{equation}
where $Q_{{\bf n},-m_2-1} \equiv 1$. For $j = m_1$ these relations
are non other than (\ref{eq:1})-(\ref{eq:2}) and the initial
induction step is settled. Let us assume that the statement is
true for $j+1 \in \{-m_2+1,\ldots,m_1\}$ and let us prove it for
$j$.

It is easy to see that the orthogonality relations (\ref{orto1})
and (\ref{orto2}) can be expressed as
\[ \int x^{\nu} Q_{{\bf n},j}(x) \frac{|Q_{{\bf
n},j+1}(x)\mathcal{A}_{{\bf n},j}(x)|}{|Q_{{\bf
n},j}(x)|}\frac{d|\sigma_j|(x)}{|Q_{{\bf n},j-1}(x)Q_{{\bf
n},j+1}(x)|} = 0\,, \qquad \nu=0,\ldots,N_{{\bf n},j}-1\,.
\]
On account of (\ref{formint}) and (\ref{formint2}) taking $q =
Q_{{\bf n},j+1}$, this can be further transformed into
\[ \int x^{\nu} Q_{{\bf n},j}(x) \left( \int \frac{Q_{{\bf
n},j+1}^2(t)}{|Q_{{\bf n},j}(t)|}\frac{|\mathcal{A}_{{\bf
n},j+1}(t)|}{|Q_{{\bf
n},j+1}(t)|}\frac{d|\sigma_{j+1}|(t)}{|x-t|}\right)\frac{d|\sigma_j|(x)}{|Q_{{\bf
n},j-1}(x)Q_{{\bf n},j+1}(x)|} = 0\,,
\]
for $\nu=0,\ldots,N_{{\bf n},j}-1\,.$

Relation (\ref{conv-Qnj}) implies that
\begin{equation}\label{eq:5}
\lim_{{\bf n} \in \Lambda'} \frac{1}{2N_{{\bf n},j}}\log|Q_{{\bf
n},j-1}(x)Q_{{\bf n},j+1}(x)| = -
\frac{P_{j-1}}{2P_j}V^{\mu_{j-1}}(x) -
\frac{P_{j+1}}{2P_j}V^{\mu_{j+1}}(x)\,,
\end{equation}
uniformly on $\Delta_j.$ (Since $Q_{{\bf n},-m_2-1}\equiv 1$, when
$j=-m_2$ we only get the second term on the right hand side of
this limit; that is, $P_{-m_2-1} = 0$.)

Set
\[ K_{{\bf n},j+1} = \left(\int \frac{Q_{{\bf
n},j+1}^2(t)}{|Q_{{\bf n},j}(t)|} \frac{ |{\mathcal{A}}_{{\bf
n},j+1}(t)|}{|Q_{{\bf n},j+1}(t)|}
d|\sigma_{j+1}|(t)\right)^{-1/2}.
\]
It follows that for all $x \in \Delta_j$
\[ \frac{1}{\delta_{j+1}^*K_{{\bf n},j+1}^2} \leq \int \frac{Q_{{\bf
n},j+1}^2(t)}{|Q_{{\bf n},j}(t)|}\frac{|{\mathcal{A}}_{{\bf
n},j+1}(t)|}{|Q_{{\bf n},j+1}(t)|}\frac{d|\sigma_{j+1}|(t)}{|x-t|}
\leq \frac{1}{\delta_{j+1}K_{{\bf n},j+1}^2},
\]
where $0 < \delta_{j+1} = \inf\{|x-t|: t \in \Delta_{j+1}, x \in
\Delta_j\} \leq \max\{|x-t|: t \in \Delta_{j+1}, x \in \Delta_j\}
= \delta_{j+1}^* < \infty.$ Taking into consideration these
inequalities, from the induction hypothesis, we obtain that
\begin{equation} \label{eq:6}
\lim_{{\bf n} \in \Lambda'} \left(\int \frac{Q_{{\bf
n},j+1}^2(t)}{|Q_{{\bf n},j}(t)|}\frac{|{\mathcal{A}}_{{\bf
n},j+1}(t)|}{|Q_{{\bf
n},j+1}(t)|}\frac{d|\sigma_{j+1}|(t)}{|x-t|}\right)^{1/2N_{{\bf
n},j}} = e^{-P_{j+1}\omega_{j+1}/P_j}.
\end{equation}

Taking (\ref{eq:5}) and (\ref{eq:6}) into account, Lemma
\ref{gonchar-rakhmanov} yields that $\mu_j$ is the unique solution
of the extremal problem (\ref{eq:3}) and
\[  \lim_{{\bf n} \in \Lambda'} \left(\int   \int \frac{Q_{{\bf
n},j+1}^2(t)}{|Q_{{\bf n},j}(t)|}\frac{|{\mathcal{A}}_{{\bf
n},j+1}(t)|}{|Q_{{\bf n},j+1}(t)|}\frac{d|\sigma_{j+1}|(t)}{|x-t|}
\frac{Q_{{\bf n},j}^2(x)d|\sigma_j|(x)}{|Q_{{\bf n},j-1}(x)Q_{{\bf
n},j+1}(x)|}\right)^{1/2N_{{\bf n},j}} = e^{-\omega_j}.
\]
According to (\ref{formint}) and (\ref{formint2}) with $q =
Q_{{\bf n},j+1}$
\[ \frac{1}{|Q_{{\bf n},j+1}(x)|} \int \frac{Q_{{\bf n},j+1}^2(t)}{|Q_{{\bf n},j}(t)|}
\frac{|{\mathcal{A}}_{{\bf n},j+1}(t)|}{|Q_{{\bf
n},j+1}(t)|}\frac{d|\sigma_{j+1}|(t)}{|x-t|}=
\frac{|{\mathcal{A}}_{{\bf n},j}(x)|}{|Q_{{\bf n},j}(x)|}\,, \quad
x \in \Delta_j\,,
\]
which allows to reduce the previous formula to (\ref{eq:4}) thus
concluding the induction.

Now, we can rewrite (\ref{eq:3}) multiplying through by $P_j^2$
and taking the constant term on the left to the right to obtain
the system of boundary value equations
\begin{equation} \label{eq:a1}
P_j^2 V^{\mu_j}(x) - \frac{P_{j-1}P_j}{2}V^{\mu_{j-1}}(x) -
\frac{P_jP_{j+1}}{2 }V^{\mu_{j+1}}(x)  \left\{
\begin{array}{l} = \omega_j',\quad x
\in \sop( {\mu_j}) \,, \\
\geq \omega_j', \quad x \in E_j \,,
\end{array} \right.
\end{equation}
for $j=-m_2,\ldots,m_1$, where
\begin{equation}\label{eq:c}
\omega_j' = P_j^2\omega_j - P_jP_{j+1}\omega_{j+1}.
\end{equation}
The terms with $P_{-m_2-1}$ and $P_{m_1+1}$ do not appear when
$j=-m_2$ and $j=m_1$, respectively. By Lemma \ref{niksor},
$(\mu_{-m_2},\ldots,\mu_{m_1})=
(\overline{\mu}_{-m_2},\ldots,\overline{\mu}_{m_1})$ and
$(\omega_{-m_2}',\ldots,\omega_{m_1}')=(\omega_{-m_2}^{\overline{\mu}},\ldots,\omega_{m_1}^{\overline{\mu}})$
for any convergent subsequence showing the existence of the limits
in (\ref{conv-Qnj*}) as stated.

Notice that (\ref{eq:4}) implies that
\[ \lim_{{\bf n} \in \Lambda'}
\left(\int \frac{Q_{{\bf n},j}^2(x)}{|Q_{{\bf
n},j-1}(x)|}\frac{|\mathcal{A}_{{\bf n},j}(x)|}{|Q_{{\bf
n},j}(x)|}d|\sigma_j|(x)\right)^{1/2|{\bf n}_1|} =
e^{-P_j\omega_j}.
\]
On the other hand, from (\ref{eq:c}) it follows that
$P_{m_1}\omega_{m_1} = \omega_{m_1}^{\overline{\mu}}/P_{m_1}$ when
$j=m_1.$  Suppose that $P_{j+1}\omega_{j+1} = \sum_{k=j+1}^{m_1}
\frac{\omega_{k}^{\overline{\mu}}}{P_{k}}, j+1 \in
\{-m_2+1,\ldots,m_1\}$. Then, according to (\ref{eq:c})
\[ P_j\omega_j = \frac{\omega_j^{\overline{\mu}}}{P_j} +
P_{j+1}\omega_{j+1} = \sum_{k=j}^{m_1}
\frac{\omega_{k}^{\overline{\mu}}}{P_{k}}
\]
and (\ref{eq:4*}) immediately follows.
\end{proof}

\section{Proof of Theorem 1}\label{thm1}

Here, we maintain the change of notation introduced in the
previous section. Theorem \ref{principal} is a consequence of the
following more general result.

\begin{thm}\label{general}
Let $\Lambda =
\Lambda(p_{1,0},\ldots,p_{1,m_1};p_{2,0},\ldots,p_{2,m_2}) \subset
\mathbb{Z}_+^{m_1+1}(\bullet) \times
\mathbb{Z}_+^{m_2+1}(\bullet), (S^1,S^2)\in\mbox{\bf {Reg}},$
$S^1= \mathcal{N}(\sigma^1_0,\ldots,\sigma^1_{m_1}),$ and $S^2=
\mathcal{N}(\sigma^2_0,\ldots,\sigma^2_{m_2})$ be given. Let
$\{\mathcal{A}_{{\bf n},j}\}, {\bf n} \in \Lambda, j =
-m_2-1,\ldots,m_1,$ be the sequences of ``monic'' linear forms
associated with the corresponding mixed type orthogonal
polynomials. Then, for each $j=-m_2-1,\ldots,m_1$
\begin{equation}\label{conver1}
\lim_{{\bf n} \in \Lambda}|\mathcal{A}_{{\bf n},j}(z)|^{1/|{\bf
n}_1|}=G_j(z), \qquad {\mathcal{K}} \subset {\C}\setminus
(\Delta_j \cup \Delta_{j+1})
\end{equation}
$( \Delta_{-m_2-1} = \Delta_{m_1+1} = \emptyset )$, where
\begin{equation}\label{eqGj1}
G_j(z)=\exp\left(P_{j+1} V^{\overline{\mu}_{j+1}}(z)-P_j
V^{\overline{\mu}_j}(z) - 2\sum_{k=j+1}^{m_1}
\frac{\omega_k^{\overline{\mu}}}{P_k}\right),\quad
j=-m_2-1,\ldots,m_1-1,
\end{equation}
$(P_{-m_2-1} = P_{m_1+1} =0)$ and
\begin{equation}\label{eqGj2}
G_{m_1}(z)=\exp\left(-P_{m_1} V^{\overline{\mu}_{m_1}}(z)\right).
\end{equation}
$\overline{\mu}= \overline{\mu}(\mathcal{C})=
(\overline{\mu}_{-m_2},\ldots,\overline{\mu}_{m_1})$ is  the
equilibrium vector measure   and
$(\omega_{-m_2}^{\overline{\mu}},\ldots,\omega_{m_1}^{\overline{\mu}})$
is the system of equilibrium constants for the vector potential
problem determined by the interaction matrix $\mathcal{C}$ defined
in $(\ref{matriz})$  on the system of compact sets $E_j = \sop
(\sigma^1_j), j=0,\ldots,m_1, E_j = \sop (\sigma^2_{-j}),
j=-m_2,\ldots,0$.
\end{thm}
\begin{proof} If $j=m_1, \mathcal{A}_{{\bf n},m_1} = Q_{{\bf n},m_1}$ and
(\ref{conv-Qnj*}) directly implies that
\[ \lim_{{\bf n}\in\Lambda} |\mathcal{A}_{{\bf n},m_1}(z)|^{1/|{\bf n}_1|}
=  \exp\left(-P_{m_1} V^{\overline{\mu}_{m_1}}(z)\right), \qquad
\mathcal{K} \subset \mathbb{C} \setminus \Delta_{m_1}.
\]
For $j \in \{-m_2-1,\ldots,m_1-1\}$, using (\ref{formint}) and
(\ref{formint2}) with $q= Q_{{\bf n},j+1}$, we obtain
\begin{equation}\label{eq:7} \mathcal{A}_{{\bf n},j}(z) =
\frac{Q_{{\bf n},j}(z)}{Q_{{\bf n},j+1}(z)}\int \frac{Q_{{\bf
n},j+1}^2(x)}{Q_{{\bf n},j}(x)} \frac{\mathcal{A}_{{\bf
n},j+1}(x)}{Q_{{\bf n},j+1}(x)} \frac{d\sigma_{j+1}(x)}{z-x},
\end{equation}
$(Q_{{\bf n},-m_2-1} \equiv 1$.) From (\ref{conv-Qnj*}), it
follows that
\begin{equation}\label{eq:8}
\lim_{{\bf n}\in\Lambda}\left|\frac{Q_{{\bf n},j}(z)}{Q_{{\bf
n},j+1}(z)}\right|^{ {1}/{|{\bf n}_1|}}=\exp\left(P_{j+1}
V^{\overline{\mu}_{j+1}}(z)-P_j V^{\overline{\mu}_j}(z)\right),
\qquad \mathcal{K} \subset \mathbb{C} \setminus (\Delta_j \cup
\Delta_{j+1})
\end{equation}
(we also use that the zeros of $Q_{{\bf n},j}$ and $Q_{{\bf
n},j+1}$ lie in $\Delta_j$ and $\Delta_{j+1}$, respectively). It
remains to find the $|{\bf n}_1|$-th root asymptotic behavior of
the integral.

Fix a compact set $\mathcal{K} \subset \mathbb{C} \setminus
\Delta_{j+1}.$ It is easy to verify that (for the definition of
$K^2_{{\bf n},j+1}$ see proof of Theorem \ref{mea-equil} above)
\[ \frac{C_1}{K_{{\bf n},j+1}^2} \leq \left|\int \frac{Q_{{\bf n},j+1}^2(x)}{Q_{{\bf
n},j}(x)} \frac{\mathcal{A}_{{\bf n},j+1}(x)}{Q_{{\bf n},j+1}(x)}
\frac{d\sigma_{j+1}(x)}{z-x}\right| \leq \frac{C_2}{K_{{\bf
n},j+1}^2} \,,
\]
where
\[ C_1 = \frac{\min \{ \max\{|u-x|,|v|: z = u+iv\}: z \in \mathcal{K}, x
\in \Delta_{j+1}\}}{ \max\{|z-x|^2: z \in \mathcal{K}, x \in
\Delta_{j+1}\}} > 0
\]
and
\[ C_2 = \frac{1}{\min\{|z-x|: z \in \mathcal{K}, x \in
\Delta_{j+1}\}} < \infty.
\]
Taking into account (\ref{eq:4*})
\begin{equation}\label{eq:9}
\lim_{{\bf n} \in \Lambda} \left|\int \frac{Q_{{\bf
n},j+1}^2(x)}{Q_{{\bf n},j}(x)} \frac{\mathcal{A}_{{\bf
n},j+1}(x)}{Q_{{\bf n},j+1}(x)}
\frac{d\sigma_{j+1}(x)}{z-x}\right|^{1/|{\bf n}_1|} = \exp\left(-2
\sum_{k=j+1}^{m_1} {\omega_k^{\overline{\mu}}}/{P_k}\right)\,.
\end{equation}
From (\ref{eq:7})-(\ref{eq:9}), we obtain (\ref{conver1}) and we
are done.
\end{proof}

\begin{rem} Taking into consideration that the polynomials $Q_{{\bf
n},j}$ (see Propositions \ref{ortogonalidad} and \ref{ortogb}) and
the functions
\[ \int \frac{Q_{{\bf n},j}^2(x)}{Q_{{\bf
n},j-1}(x)} \frac{\mathcal{A}_{{\bf n},j}(x)}{Q_{{\bf n},j}(x)}
\frac{d\sigma_{j}(x)}{z-x},
\]
may have at most one zero in each of the connected components of
$\Delta_j \setminus  E_j$, in place of (\ref{conver1}) one can
prove convergence in capacity on each compact subset $\mathcal{K}
\subset \mathbb{C} \setminus (E_j \cup E_{j+1})$. More precisely,
for any such compact set $\mathcal{K}$ and each $\varepsilon > 0$
\[ \lim_{{\bf n} \in \Lambda} \mbox{cap} \left\{z \in \mathcal{K}: \left||\mathcal{A}_{{\bf n},j}(z)|^{1/|{\bf
n}_1|}-G_j(z)\right| > \varepsilon \right\} = 0.
\] \fp
\end{rem}

Set
\[ U_j^{\overline{\mu}}(z) = P_j V^{{\overline{\mu}_j}}(z) -P_{j+1}
V^{\overline{\mu}_{j+1}}(z) + 2\sum_{k=j+1}^{m_1}
\frac{\omega_k^{\overline{\mu}}}{P_k}, \qquad
j=-m_2-1,\ldots,m_1-1,
\]
and
\[   U_{m_1}^{\overline{\mu}}(z) = P_{m_1} V^{{\overline{\mu}_{m_1}}}(z).
\]
Hence, $G_j(z) = \exp(-U_j^{\overline{\mu}}(z) ),
j=-m_2-1,\ldots,m_1.$

We have that for $j=-m_2,\ldots,m_1 (P_{-m_2-1} = P_{m_1+1} = 0)$
\[ \frac{P_j}{2}(U_j^{\overline{\mu}}(z) - U_{j-1}^{\overline{\mu}}(z)) =
-\frac{P_{j+1}P_j}{2}
 V^{{\overline{\mu}_{j+1}}}(z)+ P_j^2 V^{{\overline{\mu}_j}}(z) -
\frac{P_jP_{j-1}}{2} V^{\overline{\mu}_{j-1}}(z) -
{\omega_j^{\overline{\mu}}}\,.
\]
From the equilibrium property (see Lemma \ref{niksor} and
(\ref{eq:a1})), it follows that
\[ U_j^{\overline{\mu}}(x) - U_{j-1}^{\overline{\mu}}(x) \left\{
\begin{array}{l} = 0,\quad x
\in \sop( \overline{\mu}_j) \,, \\
\geq 0, \quad x \in E_j \,.
\end{array} \right.
\]
Define
\[ p_j =  \left\{
\begin{array}{ll}
p_{1,j}, & j=0,\ldots,m_1, \\
-p_{2,-j-1}, & j=-m_2-1,\ldots,-1.
\end{array}
\right.
\]
It is easy to verify that for $j= -m_2,\ldots,m_1$
\begin{equation} \label{eq:a2} U_j^{\overline{\mu}}(z) - U_{j-1}^{\overline{\mu}}(z) =
\mathcal{O}((p_j -p_{j-1})\log1/|z|), \qquad z \to \infty.
\end{equation}
In particular, $U_j^{\overline{\mu}}(z) -
U_{j-1}^{\overline{\mu}}(z) = \mathcal{O}(1), z \to \infty,$
whenever $p_j = p_{j-1}$. By assumption, $p_j -p_{j-1} \leq 0,
j=-m_2,\ldots,m_1$ except for  $p_0 - p_{-1} = p_{1,0} + p_{2,0} >
0$.

For all $j$, the function $U_j^{\overline{\mu}} -
U_{j-1}^{\overline{\mu}} $ is subharmonic in ${\mathbb{C}}
\setminus \sop (\overline{\mu}_j).$ If $p_j \geq p_{j-1},$  then
it is subharmonic in all $\overline{\mathbb{C}} \setminus \sop
(\overline{\mu}_j).$ According to what was said above, when $j=0$
or $p_j = p_{j-1}$, from the equilibrium condition and the maximum
principle, we have that $U_j^{\overline{\mu}}  -
U_{j-1}^{\overline{\mu}}  \equiv 0$ on $\sop ({\sigma}_j) = E_j$
and $U_j^{\overline{\mu}}  < U_{j-1}^{\overline{\mu}} $ on
$\overline{\mathbb{C}} \setminus \sop ( {\sigma}_j)$. In
particular, in this case we have that $\sop ( \overline{\mu}_j) =
\sop ( {\sigma}_j).$

When $p_{j-1} > p_j$, (\ref{eq:a2}) implies that in a neighborhood
of $z=\infty, U_j^{\overline{\mu}} > U_{j-1}^{\overline{\mu}}$.
Let $\gamma_j = \{z \in {\mathbb{C}}: U_j^{\overline{\mu}}(z) -
U_{j-1}^{\overline{\mu}}(z) =  0\}$. The equilibrium condition
entails that $\gamma_j \supset \sop(\overline{\mu}_j)$ and the
initial remark indicates that $\gamma_j$ is bounded. Consider any
bounded component of the complement of $\gamma_j$. On it,
$U_j^{\overline{\mu}} - U_{j-1}^{\overline{\mu}}$ is subharmonic
and on its boundary $U_j^{\overline{\mu}} -
U_{j-1}^{\overline{\mu}}=0$. Thus, on any bounded component of the
complement of $\gamma_j$ we have that $U_j^{\overline{\mu}} <
U_{j-1}^{\overline{\mu}}$. From the initial remark  it follows
that on the unbounded component of the complement of $\gamma_j,
U_j^{\overline{\mu}}
> U_{j-1}^{\overline{\mu}}.$

Fix $j \in \{0,\ldots,m_1\}$. For each $k \in \{j,\ldots,m_1\}$
define
\[ D^j_k := \{z \in {\mathbb{C}} \setminus \cup_{i=j}^{m_1} \Delta_i: U_k^{\overline{\mu}}(z) <
U_i^{\overline{\mu}}(z), i = j,\ldots,m_1, i\neq k\}, \qquad
D^{m_1}_{m_1} := \mathbb{C} \setminus \Delta_{m_1}.
\]
Let
\[\zeta_j(z) = \min\{U_k^{\overline{\mu}}(z): k = j,\ldots,m_1\}.
\]

\begin{cor}\label{polinom}
Let $\Lambda =
\Lambda(p_{1,0},\ldots,p_{1,m_1};p_{2,0},\ldots,p_{2,m_2}) \subset
\mathbb{Z}_+^{m_1+1}(\bullet) \times
\mathbb{Z}_+^{m_2+1}(\bullet), (S^1,S^2)\in\mbox{\bf {Reg}},$ $
S^1= \mathcal{N}(\sigma^1_0,\ldots,\sigma^1_{m_1}),$ $S^2=
\mathcal{N}(\sigma^2_0,\ldots,\sigma^2_{m_2})$ be given. Let
$(a_{{\bf n},0},a_{{\bf n},1},\ldots,a_{{\bf n},m_1}), {\bf n} \in
\Lambda, $ be the associated sequence of ``monic'' mixed type
multiple orthogonal polynomials. Then, for $j=0,\ldots,m_1$
\begin{equation}\label{eq:a3}
\lim_{{\bf n} \in \Lambda}|a_{{\bf n},j}(z)|^{1/|{\bf n}_1|}=
\exp(-\zeta_j(z)), \qquad \mathcal{K} \subset \cup_{k=j}^{m_1}
D_k^j,
\end{equation}
and
\begin{equation}\label{eq:a4}
\limsup_{{\bf n} \in \Lambda}|a_{{\bf n},j}(z)|^{1/|{\bf n}_1|}
\leq \exp(-\zeta_j(z)), \qquad {\mathcal{K}} \subset {\mathbb{C}}
\setminus \cup_{k=j}^{m_1} \Delta_k .
\end{equation}
In particular, if $p_{1,0} = \cdots = p_{1,m_1} = 1/(m_1+1)$, then
\begin{equation}\label{eq:a5}
\lim_{{\bf n} \in \Lambda}|a_{{\bf n},j}(z)|^{1/|{\bf n}_{1}|}=
\exp(-U^{\overline{\mu}}_{m_1}(z)), \qquad {\mathcal{K}} \subset
\mathbb{C} \setminus \cup_{k=j}^{m_1} \Delta_k.
\end{equation}
\end{cor}

\begin{proof}
For $j=m_1$, $\mathcal{A}_{{\bf n},m_1} = a_{{\bf n},m_1},
D^{m_1}_{m_1} = \mathbb{C} \setminus \Delta_{m_1}$ and
$\zeta_{m_1} = U^{\overline{\mu}}_{m_1}$. Therefore, (\ref{eq:a3})
reduces to (\ref{conver1}) and implies (\ref{eq:a4}). Let us prove
these relations for $j=0,\ldots,m_1-1$.

The $\mathcal{A}_{{\bf n},j}$ are expressed in terms of the
$a_{{\bf n},k}, k = j,\ldots,m_1,$ through a linear triangular
scheme of equations with function coefficients which do not depend
on ${\bf n}.$  Using this system, we can solve for $a_{{\bf
n},j},$ in terms of $\mathcal{A}_{{\bf n},k}, k = j,\ldots,m_1.$

Given $j \in \{1,\ldots,m_1\}$ and $0 \leq i < j,$ we have
\[  (-1)^{j-i} \langle\sigma^1_i,\ldots,\sigma^1_j\rangle^{\widehat{}}(z) = \int\cdots\int \frac{ d\sigma^1_i(x_i)\cdots
d\sigma^1_j(x_j)}{(z-x_i)(x_{i+1}-x_i)\cdots (x_j-x_{j-1})},
\]
where $\langle \cdot \rangle^{\widehat{}}(z)$ denotes the Cauchy
transform of the indicated measure, and
\[ \langle\sigma^1_j,\ldots,\sigma^1_i\rangle^{\widehat{}}(z) = \int\cdots\int \frac{ d\sigma^1_i(x_i)\cdots
d\sigma^1_j(x_j)}{(x_{i+1}-x_i)\cdots (x_j-x_{j-1})(z-x_j)}.
\]
Consequently,
\[ (-1)^{j-i} \langle\sigma^1_i,\ldots,\sigma^1_j\rangle^{\widehat{}}(z)-
\langle\sigma^1_j,\ldots,\sigma^1_i\rangle^{\widehat{}}(z)  =
\int\cdots\int \frac{-(x_j-x_i) d\sigma^1_i(x_i)\cdots
d\sigma^1_j(x_j)}{(z-x_i)(x_{i+1}-x_i)\cdots
(x_j-x_{j-1})(z-x_j)}.
\]
Since $x_j -x_i = x_j - x_{j-1} + x_{j-1}- \cdots - x_{i+1} +
x_{i+1} -x_i$, substituting this in the previous formula, we
obtain
\begin{equation} \label{formula}
\langle\sigma^1_j,\ldots,\sigma^1_i\rangle^{\widehat{}}(z) =
\sum_{k=i}^{j-1}(-1)^{k-i}
\langle\sigma^1_i,\ldots,\sigma^1_k\rangle^{\widehat{}}(z) \langle
\sigma^1_j, \ldots,\sigma^1_{k+1}\rangle^{\widehat{}}(z) +
(-1)^{j-i}
\langle\sigma^1_i,\ldots,\sigma^1_j\rangle^{\widehat{}}(z).
\end{equation}
(This formula is applicable to any Nikishin system. We will use it
on $S^2$ in the last section.)

Using formula (\ref{formula}) it is easy to deduce that (the sum
is empty when $j=m_1$)
\[ a_{{\bf n},j}(z) = \mathcal{A}_{{\bf n},j}(z) + \sum_{k=j+1}^{m_1}
(-1)^{k-j} \langle \sigma^1_k, \ldots,\sigma^1_{j+1}
\rangle^{\widehat{}}(z) \mathcal{A}_{{\bf n},k}(z).
\]
Taking (\ref{conver1})  into consideration, on $D^j_k$ the term
containing   $\mathcal{A}_{{\bf n},k}$ dominates the sum (notice
that $\langle \sigma^1_k, \ldots,\sigma^1_{j+1}
\rangle^{\widehat{}}(z) \neq 0, z \in \mathbb{C} \setminus
\Delta_{k}$) and (\ref{eq:a3}) immediately follows. On the
complement of $\cup_{k=j}^{m_1} D^j_k$ there is no dominating term
and all we can conclude from the previous equality is
(\ref{eq:a4}).

Let $p_{1,0}=\cdots=p_{1,m_1} = 1/(m_1+1)$. In this case, on
$\mathbb{C} \setminus \cup_{k=j}^{m_1}\Delta_k$ we have that
$U^{\overline{\mu}}_{m_1}(z) < U^{\overline{\mu}}_{m_1-1}(z) <
\cdots < U^{\overline{\mu}}_j(z)$ and (\ref{eq:a5}) follows from
(\ref{eq:a3}).
\end{proof}

\section{Ratio asymptotic}\label{ratio}

Here, we study the convergence of the sequences  $\{Q_{{\bf
n}^l,j}/Q_{{\bf n},j}\}, {\bf n} \in \Lambda \subset
\mathbb{Z}_+^{m_1+1}(\bullet)\times \mathbb{Z}_+^{m_2+1}(\bullet)$
and of the ratio of the corresponding linear forms. We maintain
the notation introduced in Section  \ref{ext-pro}, namely
\[\Delta_j = \Delta^1_{j}, \qquad \sigma_j = \sigma^1_{j}, \qquad j=
0,1,\ldots,m_1,\]
\[\Delta_j = \Delta^2_{-j}, \qquad \sigma_j = \sigma^2_{-j}, \qquad j=
0,-1,\ldots,-m_2,\] and
\begin{equation} \label{eq:N}
N_{{\bf n},j} = \left\{
\begin{array}{ll}
N_{1,j}({\bf n}) -1, & j=0,1\ldots,m_1, \\ N_{2,-j}({\bf n}), & j=
0,-1,\ldots,-m_2.
\end{array}
\right.
\end{equation}
Set
\[ \mathcal{H}_{{\bf n},j} = \frac{Q_{{\bf n},j+1}\mathcal{A}_{{\bf n},j}}{Q_{{\bf
n},j}}, \qquad j = -m_2-1,\ldots,m_1,
\]
($Q_{{\bf n},-m_2-1} \equiv Q_{{\bf n},m_1+1}\equiv 1$ and
$\mathcal{H}_{{\bf n},m_1} \equiv 1$). With these notations,
relations (\ref{orto1}), (\ref{orto2}), (\ref{formint}), and
(\ref{formint2}) (replacing general $q$ by $Q_{{\bf n},j+1}$ and
shifting the index $j$ by $-1$) can be rewritten as follows
\begin{equation} \label{orto3}
\int x^{\nu} Q_{{\bf n},j}(x) \frac{|\mathcal{H}_{{\bf
n},j}(x)|\,d |\sigma_{j}|(x)}{|Q_{{\bf n},j-1}(x)Q_{{\bf
n},j+1}(x)|}=0,  \quad \nu=0,\ldots,N_{{\bf n},j}-1, \quad
j=-m_2,\ldots,m_1,
\end{equation}
and
\begin{equation} \label{formint3}
\mathcal{H}_{{\bf n},j-1}(z) = \int \frac{Q_{{\bf n},j}^2(x)}{z-x}
\frac{\mathcal{H}_{{\bf n},j}(x)d\sigma_{j}(x)}{Q_{{\bf
n},j-1}(x)Q_{{\bf n},j+1}(x)}, \qquad j=-m_2,\ldots,m_1.
\end{equation}
Since on the interval $\Delta_j$ the measure $\sigma_j$ and the
functions $\mathcal{H}_{{\bf n},j},Q_{{\bf n},j-1}Q_{{\bf
n},j+1},$ preserve a constant sign, we can take their absolute
values in (\ref{orto3}) without altering the orthogonality
relations.

For each $j=-m_2,\ldots,m_1,$ define
\begin{equation} \label{eq:K} K_{{\bf n},j} =   \left(\int_{\Delta_j}
Q_{{\bf n},j}^2(x) \frac{|\mathcal{H}_{{\bf n},j}(x)|
\,d|\sigma_j|(x)}{|Q_{{\bf n},j-1}(x)Q_{{\bf n},j+1}(x)|}
\right)^{-1/2}.
\end{equation}
Take
\[ K_{{\bf n},m_1+1} = 1 \;, \quad \kappa_{{\bf n},j} = \frac{K_{{\bf n},j}}{K_{{\bf
n},j+1}} \;, \quad j=-m_2,\ldots,m_1 \;.
\]
Define
\begin{equation} \label{eq:orton} q_{{\bf n},j} =
\kappa_{{\bf n},j}Q_{{\bf n},j} \;, \quad h_{{\bf n},j-1}  =
K_{{\bf n},j}^2 \mathcal{H}_{{\bf n},j-1} \;,
\end{equation}
and
\begin{equation} \label{eq:var1}
d\rho_{{\bf n},j}(x) = \frac{h_{{\bf n},j}(x)d\sigma_j(x)}{Q_{{\bf
n},j-1}(x)Q_{{\bf n},j+1}(x)} \,.
\end{equation}
From (\ref{orto3}) and the notation introduced above, we obtain
\begin{equation} \label{orton1}
\int_{\Delta_j} x^{\nu} Q_{{\bf n},j}(x)d|\rho_{{\bf
n},j}|(x)=0,\quad \nu=0,\ldots,N_{{\bf n},j} -1,\quad
j=-m_2,\ldots,m_1\,,
\end{equation}
and $q_{n,j}$ is orthonormal with respect to the varying measure
$|\rho_{{\bf n},j}|.$ On the other hand, using (\ref{formint3}) it
follows that
\begin{equation} \label{eq:5.7*}
h_{{\bf n},j-1}(z)= \varepsilon_{{\bf n},j}\int_{\Delta_{j}}
\frac{ q^2_{{\bf n},j}(x)}{z-x}d|\rho_{{\bf n},j}|(x)\,, \quad
j=-m_2,\ldots,m_1,
\end{equation}
where $\varepsilon_{{\bf n},j}$ denotes the sign of the varying
measure $\rho_{{\bf n},j}$.

The proof of Theorem \ref{teofund} below has three steps. First,
we show that for each $j \in \{-m_2,\ldots,m_1\}$ the sequence
$\{Q_{{\bf n}^l,j}/Q_{{\bf n},j}\}$ is uniformly bounded on each
compact subset contained in $\mathbb{C} \setminus \sop(\sigma_j)$
(for all sufficiently large $|{\bf n}_1|$). Taking a subsequence
of multi-indices such that all the sequences of ratios of
polynomials have limit, we show that the limit functions must
satisfy a system of boundary value problems. This system happens
to have a unique solution from which we derive that all convergent
subsequences have the same limit. Finally, we show that the limit
functions can be expressed in terms of the branches of certain
conformal representations of a related compact Riemann surface
onto the extended complex plane.

In this section, we assume that $\sop (\sigma_j) =
\widetilde{\Delta}_j \cup e_j, j=-m_2,\ldots,m_1$, where
$\widetilde{\Delta}_j = [a_j,b_j]$ is a bounded interval of the
real line, $|\sigma_j^{\prime}| > 0$ a.e. on
$\widetilde{\Delta}_j$, and $e_j$ is a set without accumulation
points in $\overline{\mathbb{R}} \setminus \widetilde{\Delta}_j$.
We denote this writing $S^1=
\mathcal{N}'(\sigma^1_0,\ldots,\sigma^1_{m_1}), S^2=
\mathcal{N}'(\sigma^2_0,\ldots,\sigma^2_{m_2})$.

In order to fulfill the first step, Theorem \ref{theointerlace}
would be sufficient if $\Delta_j = \widetilde{\Delta}_j,
j=-m_2,\ldots,m_1.$ In order to allow the compact sets to enter
the connected components of $\Delta_j \setminus \sop(\sigma_j)$,
we need to show that the zeros falling in the intervals $I$ (see
Propositions \ref{ortogonalidad} and \ref{ortogb}) are attracted
to points in $\sop(\sigma_j) \setminus \widetilde{\Delta}_j$. In
our aid comes the next result.

\begin{lem} \label{lm:3}
Let $S^1= \mathcal{N}'(\sigma^1_0,\ldots,\sigma^1_{m_1}), S^2=
\mathcal{N}'(\sigma^2_0,\ldots,\sigma^2_{m_2})$ be given, and let
$\Lambda \subset \mathbb{Z}_+^{m_1+1}(\bullet)\times
\mathbb{Z}_+^{m_2+1}(\bullet)$ be an infinite sequence of distinct
multi-indices such that
\begin{equation}\label{eq:multcond}
\displaystyle{\sup_{{\bf n} \in \Lambda }((m_2 +1) n_{2,0} - |{\bf
n}_2| ) }< \infty\,, \qquad \displaystyle{\sup_{{\bf n} \in
\Lambda }((m_1+1)n_{1,0} -  |{\bf n}_1| ) }<
\infty\,.\end{equation} For any continuous function $f$ on
$\supp{\sigma_{j}}$
\begin{equation} \label{eq:5.7**}
\lim_{{\bf n} \in \Lambda}  \int f(x)q^2_{{\bf n},j}(x)
d|\rho_{{\bf n},j}|(x) = \frac{1}{\pi} \int_{a_j}^{b_j} f(x)
\frac{d x}{\sqrt{(b_j -x)(x - a_j)}} \,,
\end{equation}
where $\widetilde{\Delta}_j =[a_j,b_j]$, $-m_{2}\leq j\leq m_{1}$.
In particular,
\begin{equation} \label{eq:h}
\lim_{{\bf n} \in \Lambda} \varepsilon_{{\bf n},j}h_{{\bf
n},j-1}(z) = \frac{1}{\sqrt{(z - b_j)(z - a_j)}}\,,\quad
\mathcal{K} \subset \mathbb{C} \setminus \supp{\sigma_{j}}\,,
\end{equation}
where $\sqrt{(z - b_j)(z - a_j)} > 0$ if $z > b_{j}$.
Consequently, for $j = -m_2,\ldots,m_1,$ each point of
$\supp{\sigma_j} \setminus \widetilde{\Delta}_j$ is a limit of
zeros of $\{Q_{{\bf n},j}\}, {\bf n} \in \Lambda.$
\end{lem}

\begin{proof} We will prove this by induction on $j$. For $j=m_{1}$,
using Corollary 3 in \cite{BCL} and the second condition in
(\ref{eq:multcond}), it follows that
\[ \lim_{{\bf n} \in \Lambda}  \int_{\Delta_{m_{1}}} f(x)q^2_{{\bf
n},m_{1}}(x) \frac{ d|\sigma_{m_{1}}|(x)}{| Q_{{\bf
n},m_{1}-1}(x)|} = \frac{1}{\pi} \int_{\widetilde{\Delta}_{m_{1}}}
f(x) \frac{d x}{\sqrt{(b_{m_{1}} -x)(x - a_{m_{1}})}}\,,
\]
where $f$ is continuous on $\sop ({\sigma_{m_{1}}})$. Take $f(x) =
(z-x)^{-1}$ where $z \in \mathbb{C} \setminus \sop
({\sigma_{m_{1}}})$. According to (\ref{eq:5.7*}) and the previous
limit one obtains that
\[ \lim_{{\bf n} \in \Lambda} \varepsilon_{{\bf n},m_{1}}h_{{\bf n},m_{1}-1}(z) =
\frac{1}{\sqrt{(z - b_{m_{1}})(z - a_{m_{1}})}} =: h_{m_{1}}(z)\,,
\]
pointwise on $\mathbb{C} \setminus \sop ({\sigma_{m_{1}}})$. Since
\[\left| \int_{\Delta_{m_{1}}} \frac{q^2_{{\bf
n},m_{1}}(x)}{z-x} \frac{ d|\sigma_{m_{1}}|(x)}{|Q_{{\bf
n},m_{1}-1}(x)|}\right| \leq \frac{1}{d(\mathcal{K},\sop
(\sigma_{m_{1}}))}\,, \quad z \in \mathcal{K} \subset \mathbb{C}
\setminus \sop ({\sigma_{m_{1}}})\,,
\]
where $d(\mathcal{K},\sop ({\sigma_{m_{1}}}))$ denotes the
distance between the two compact sets, the sequence $\{h_{{\bf
n},m_{1}-1}\},$ ${\bf n} \in \Lambda,$ is uniformly bounded on
compact subsets of $\mathbb{C} \setminus \sop ({\sigma_{m_{1}}})$
and (\ref{eq:h}) follows for $j=m_{1}$.

Let $\zeta \in \sop ({\sigma_{m_{1}}}) \setminus
\widetilde{\Delta}_{m_{1}}$. Take $r > 0$ sufficiently small so
that the circle $C_{r} = \{z: |z - \zeta| = r \}$ surrounds no
other point of $\sop ({\sigma_{m_{1}}}) \setminus
\widetilde{\Delta}_{m_{1}}$ and contains no zero of $q_{{\bf
n},m_{1}}, {\bf n} \in \Lambda$. From (\ref{eq:h}) for $j=m_{1}$
\begin{equation} \label{argument}
\lim_{{\bf n} \in \Lambda} \frac{1}{2\pi i}\int_{C_r}
\frac{\varepsilon_{{\bf n},m_{1}}h_{{\bf
n},m_{1}-1}^{\prime}(z)}{\varepsilon_{{\bf n},m_{1}}h_{{\bf
n},m_{1}-1}(z)} \,d z= \frac{1}{2\pi i}\int_{C_r}
\frac{h_{m_{1}}^{\prime}(z)}{h_{m_{1}}(z)}\,dz = 0\,.
\end{equation}

Since $\zeta$ is a mass point of $\sigma_{m_1}$, formula
(\ref{eq:5.7*}) indicates that either $h_{{\bf n},m_1-1}$  has a
simple pole at $\zeta$ or $Q_{{\bf n},m_1}(\zeta)=0$.  In any
case, from (\ref{argument}) and the argument principle, it follows
that for all sufficiently large $|{\bf n}|, {\bf n} \in \Lambda$,
$Q_{{\bf n},m_{1}}$ must have a simple zero inside $C_r$. The
parameter $r$ can be taken arbitrarily small; therefore, the last
statement of the lemma readily follows and the basis of induction
is fulfilled.

Let us assume that the lemma is satisfied for $j \in
\{k+1,\ldots,m_{1}\}, -m_{2}\leq k \leq m_{1}-1,$ and let us prove
that it is also true for $k$. From (\ref{eq:h}) applied to
$j=k+1$, we have that
\[ \lim_{{\bf n} \in \Lambda} |h_{{\bf n},k}(x)| =
\frac{1}{\sqrt{|(x - b_{k+1})(x - a_{k+1})|}}\,,
\]
uniformly on $\Delta_{k} \subset \mathbb{C} \setminus \sop
(\sigma_{k+1})$. It follows that $\{|h_{{\bf
n},k}|d|\sigma_{k}|\},{\bf n} \in \Lambda,$ is a sequence of
Denisov type measures according to Definition 3 in \cite{BCL}.
Additionally, $(\{|h_{{\bf n},k}|d|\sigma_{k}|\},\{|Q_{{\bf
n},k-1} Q_{{\bf n},k+1}|\},l),{\bf n} \in \Lambda,$ is strongly
admissible as in Definition 2 of \cite{BCL} for each $l \in
\mathbb{Z}$ (see paragraph just after both definitions in the
referred paper). Therefore, we can apply Corollary 3 in \cite{BCL}
of which (\ref{eq:5.7**}) for $j=k$ is a particular case. In the
proof of Corollary 3 of \cite{BCL} (see also Theorem 9 in
\cite{BernardoGuillermo1}) it is required that the inequality
$\deg (Q_{{\bf n},j-1}Q_{{\bf n},j+1}) - 2\deg(Q_{{\bf n},j}) \leq
C$ holds for every ${\bf n}\in \Lambda$, where $C \geq 0$ is a
constant. It is straightforward to check that this condition is
satisfied under the assumption (\ref{eq:multcond}).

Now we return to the induction argument. From (\ref{eq:5.7**}) for
$j=k$, (\ref{eq:h}) and the rest of the statements of the lemma
immediately follow just as in the case when $j=m_{1}$. With this
we conclude the proof.
\end{proof}

Now, we are ready to  prove normality.

\begin{lem} \label{lm:3a}
Let $S^1= \mathcal{N}'(\sigma^1_0,\ldots,\sigma^1_{m_1}), S^2=
\mathcal{N}'(\sigma^2_0,\ldots,\sigma^2_{m_2})$ be given, and let
$\Lambda \subset \mathbb{Z}_+^{m_1+1}(\bullet)\times
\mathbb{Z}_+^{m_2+1}(\bullet)$ be an infinite sequence of distinct
multi-indices such that
\[
\displaystyle{\sup_{{\bf n} \in \Lambda }((m_2 +1) n_{2,0} - |{\bf
n}_2| ) }< \infty\,, \qquad \displaystyle{\sup_{{\bf n} \in
\Lambda }((m_1+1)n_{1,0} -  |{\bf n}_1| ) }< \infty\,.
\]
Let us assume that there exists $l=(l_{1};l_{2}), 0\leq l_{1}\leq
m_{1}, 0\leq l_{2}\leq m_{2},$ such that for all ${\bf n} \in
\Lambda$ we have that ${\bf n}^l \in
\mathbb{Z}_{+}^{m_1+1}(\bullet)\times
\mathbb{Z}_{+}^{m_2+1}(\bullet).$ Then, for each
$j=-m_2,\ldots,m_1,$ and each compact set $\mathcal{K} \subset
\mathbb{C}\setminus\sop({\sigma_j})$ there exist positive
constants $C_{j,1}(\mathcal{K}), C_{j,2}(\mathcal{K})$ such that
\[C_{j,1}(\mathcal{K})\leq \inf_{z\in \mathcal{K}}\left| \frac{Q_{{\bf
n}^l,j}(z)}{Q_{{\bf n},j}(z)}\right| \leq \sup_{z\in
\mathcal{K}}\left| \frac{Q_{{\bf n}^l,j}(z)}{Q_{{\bf
n},j}(z)}\right| \leq C_{j,2}(\mathcal{K}),
\]
for all sufficiently large $|{\bf n}_1|, {\bf n}\in\Lambda$.
\end{lem}

\begin{proof} The uniform bound from above and below on each fixed compact
subset $\mathcal{K} \subset \mathbb{C} \setminus \Delta_j$ (for
all ${\bf n} \in \Lambda$) is a direct consequence of the
interlacing property of the zeros of $Q_{{\bf n}^l,j}$ and
$Q_{{\bf n},j}$. In fact, comparing distances to $z \in
\mathcal{K}$ of consecutive interlacing zeros, it is easy to
verify that
\[\min\Big\{d_1,\frac{d_{1}}{d_{2}}\Big\}
\leq \inf_{z\in \mathcal{K}}\left| \frac{Q_{{\bf
n}^l,j}(z)}{Q_{{\bf n},j}(z)}\right| \leq \sup_{z\in
\mathcal{K}}\left| \frac{Q_{{\bf n}^l,j}(z)}{Q_{{\bf
n},j}(z)}\right| \leq \frac{\max\{d_{2},d_{2}^{2}\}}{d_{1}}\,,
\]
where $d_2$ denotes the diameter of $ \mathcal{K}\cup \Delta_j$
and $d_1$ denotes the distance between $\mathcal{K}$ and
$\Delta_j$. So, for such compact sets the assertion has been
proved.

The additional restrictions made in the lemma guarantee that the
zeros of the polynomials $Q_{{\bf n}^l,j}$ and $Q_{{\bf n},j}$
lying in $\Delta_j \setminus \sop{(\sigma_{j})}$ converge to the
mass points as results from Lemma \ref{lm:3}.   Let $\mathcal{K}
\subset \mathbb{C} \setminus \sop{(\sigma_j)}$ and suppose that
$\mathcal{K} \cap \Delta_j \neq \emptyset$. Notice that $
\mathcal{K}$ can intersect at most a finite number of open
intervals $I_1,\ldots,I_M $ forming the connected components of
$\Delta_j \setminus \sop{(\sigma_{j})}$. The polynomials $Q_{{\bf
n}^l,j}$ and $Q_{{\bf n},j}$ can have at most one zero in each of
those intervals. Consequently, for all $|{\bf n}_1|, {\bf n} \in
\Lambda,$ sufficiently large, the zeros of $Q_{{\bf n}^l,j}$ and
$Q_{{\bf n},j}$ lie at a positive distance $\varepsilon$ from
$\mathcal{K}$. Now, it is easy to show that for all sufficiently
large $|{\bf n}_1|$
\[\min\Big\{\varepsilon,\frac{\varepsilon}{d_{2}}\Big\} \leq \inf_{z\in \mathcal{K}}
\left| \frac{Q_{{\bf n}^l,j}(z)}{Q_{{\bf n},j}(z)}\right| \leq
\sup_{z\in \mathcal{K}}\left| \frac{Q_{{\bf n}^l,j}(z)}{Q_{{\bf
n},j}(z)}\right| \leq
\frac{\max\{d_{2},d_{2}^{2}\}}{\varepsilon}\,.
\]
This concludes the proof.
\end{proof}

From Lemma \ref{lm:3a} we know  that the sequences
\[\left\{ {Q_{{\bf n}^l,j}}/{Q_{{\bf n},j}}\right\}_{{\bf n}\in\Lambda},\qquad
j=-m_2,\ldots,m_1,
\]
are uniformly bounded on each compact subset of
$\mathbb{C}\setminus\sop{(\sigma_{j})}$ for all sufficiently large
$|{\bf n}_1|$. By Montel's theorem, there exists a subsequence of
multi-indices $\Lambda' \subset \Lambda$ and a collection of
functions $\widetilde{F}^{(l)}_{j}$, holomorphic in
$\mathbb{C}\setminus\sop{(\sigma_{j})}$, such that
\begin{equation}\label{eq26}
\lim_{{\bf n}\in\Lambda'}\,\frac{Q_{{\bf n}^l,j}(z)}{Q_{{\bf
n},j}(z)}=\widetilde{F}^{(l)}_{j}(z),\quad \mathcal{K}
\subset\mathbb{C}\setminus\sop{(\sigma_j)},\quad
j=-m_2,\ldots,m_1.
\end{equation}

In principle, the functions $\widetilde{F}^{(l)}_{j}$ may depend
on $\Lambda'$. We shall see that this is not the case and,
therefore, the limit in (\ref{eq26}) holds for ${\bf n} \in
\Lambda$. First, let us obtain some general information on the
functions $\widetilde{F}^{(l)}_{j}$.

The points in $\sop{(\sigma_j)} \setminus \widetilde{\Delta}_j$
are isolated singularities of $\widetilde{F}^{(l)}_{j}$. Let
$\zeta \in \sop{(\sigma_j)} \setminus \widetilde{\Delta}_j$. By
Lemma \ref{lm:3}, $\zeta~$ is a limit of zeros of $Q_{{\bf n},j}$
and $Q_{{\bf n}^l,j}$ as $|{\bf n}_1| \to \infty, {\bf n} \in
\Lambda,$ and in a sufficiently small neighborhood of $\zeta$, for
large $|{\bf n}_1|, {\bf n} \in \Lambda$, there can be at most one
zero of these polynomials (so there is exactly one, for all
sufficiently large $|{\bf n}_1|$). Let $\lim_{{\bf n} \in \Lambda}
\zeta_{\bf n} = \zeta$ where $Q_{{\bf n},j}(\zeta_{\bf n}) =0$.
From formula (\ref{eq26})
\[ \lim_{{\bf n}\in\Lambda'}\,\frac{(z - \zeta_{\bf n})Q_{{\bf n}^l,j}(z)}{Q_{{\bf
n},j}(z)}=(z - \zeta)\widetilde{F}^{(l)}_{j}(z),\quad \mathcal{K}
\subset (\mathbb{C}\setminus\sop{(\sigma_j)}) \cup \{\zeta\}\,,
\]
and  $(z - \zeta)\widetilde{F}^{(l)}_{j}(z)$ is analytic in a
neighborhood of $\zeta$. Hence $\zeta$ is not an essential
singularity of $\widetilde{F}^{(l)}_{j}$. Taking into
consideration that $Q_{{\bf n}^l,j}, {\bf n}\in \Lambda,$ also has
a sequence of zeros converging to $\zeta$, from the argument
principle it follows that $\zeta$ is a removable singularity of
$\widetilde{F}^{(l)}_{j}$ which is not a zero. By Lemma
\ref{lm:3a} we also know that the sequence of functions  $|Q_{{\bf
n}^l,j}/Q_{{\bf n},j}|, {\bf n} \in \Lambda,$ is uniformly bounded
from below by a positive constant for all sufficiently large
$|{\bf n}_1|$. Therefore, in $\mathbb{C}\setminus\sop{(\sigma_j)}$
the function $\widetilde{F}^{(l)}_{j}$ is also different from
zero. According to the definitions of $Q_{{\bf n},j}, Q_{{\bf
n}^l,j},$ and Propositions \ref{ortogonalidad} and \ref{ortogb}
(see also (\ref{eq:N})), when $-l_2 \leq j \leq l_1,$ we have that
$\deg{Q_{{\bf n}^l,j}} = N_{{\bf n}^l,j} = N_{{\bf n},j}+1=
\deg{Q_{{\bf n},j}} +1$ whereas, for $j \in \{-m_2,\ldots,-l_2-1\}
\cup \{l_1+1,\ldots,m_1\} $, we obtain that $\deg{Q_{{\bf n}^l,j}}
= N_{{\bf n}^l,j} = N_{{\bf n},j} = \deg{Q_{{\bf n},j}}$.
Consequently, when $-l_2 \leq j \leq l_1,$ the function
$\widetilde{F}^{(l)}_{j}$ has a simple pole at infinity and
$(\widetilde{F}^{(l)}_{j})^{\prime}(\infty) =1$, whereas, for $j
\in \{-m_2,\ldots,-l_2-1\} \cup \{l_1+1,\ldots,m_1\} $, it is
analytic at infinity and $\widetilde{F}^{(l)}_{j}(\infty) =1$.

Now, let us prove that the functions $\widetilde{F}^{(l)}_{j},
j=-m_2,\ldots,m_1,$ satisfy a system of boundary value problems.

\begin{lem} \label{sistema}
Let $S^1= \mathcal{N}'(\sigma^1_0,\ldots,\sigma^1_{m_1}), S^2=
\mathcal{N}'(\sigma^2_0,\ldots,\sigma^2_{m_2})$ be given, and let
$\Lambda \subset \mathbb{Z}_+^{m_1+1}(\bullet)\times
\mathbb{Z}_+^{m_2+1}(\bullet)$ be an infinite sequence of distinct
multi-indices such that
\[
\displaystyle{\sup_{{\bf n} \in \Lambda }((m_2 +1) n_{2,0} - |{\bf
n}_2| ) }< \infty, \qquad \displaystyle{\sup_{{\bf n} \in \Lambda
}((m_1+1)n_{1,0} -  |{\bf n}_1| ) }< \infty.
\]
Let us assume that there exists $l=(l_{1};l_{2}), 0\leq l_{1}\leq
m_{1}, 0\leq l_{2}\leq m_{2},$ such that for all ${\bf n} \in
\Lambda$ we have that ${\bf n}^l \in
\mathbb{Z}_{+}^{m_1+1}(\bullet)\times
\mathbb{Z}_{+}^{m_2+1}(\bullet).$ Then, there exists a
normalization $F^{(l)}_{j}$, $j=-m_2,\ldots,m_1,$ by positive
constants, of the functions $\widetilde{F}^{(l)}_{j}$  given in
$(\ref{eq26})$, which verifies the system of boundary value
problems
\begin{equation} \label{sisec}
\aligned &1) \qquad
 F^{(l)}_{j},\, {1}/{F^{(l)}_{j}}\in
H(\mathbb{C}\setminus \widetilde{\Delta}_{j})\,,\\
&2) \qquad (F^{(l)}_{j})'(\infty)>0\,,\quad j \in \{-l_2,\ldots,l_1\}\,,\\
&2') \qquad F^{(l)}_{j}(\infty)>0\,,
\quad j \in \{-m_2,\ldots,-l_2-1\} \cup \{l_1+1,\ldots,m_1\}\,,\\
&3) \qquad
|F^{(l)}_{j}(x)|^{2}\frac{1}{|(F^{(l)}_{j-1}\,F^{(l)}_{j+1})(x)|}=1,\,\,
x\in \widetilde{\Delta}_{j}\,,\\
\endaligned
\end{equation}
 where $F^{(l)}_{-m_2-1}\equiv\,F^{(l)}_{m_1+1}\equiv\,1$.
\end{lem}
\begin{proof}
The assertions 1), 2), and 2') were proved above for the functions
$\widetilde{F}_j^{(l)}$. Consequently, they are satisfied by any
normalization of these functions by means of positive constants.

From (\ref{orton1}) applied to ${\bf n}$ and ${\bf n}^l$, for each
$j=-m_2,\ldots,m_1$, we have
\[
\int x^{\nu} Q_{{\bf n},j}(x)d|\rho_{{\bf n},j}|(x)=0,\qquad
\nu=0,\ldots,N_{{\bf n},j} -1\,,
\]
and
\[
\int x^{\nu} Q_{{\bf n}^l,j}(x) g_{{\bf n},j}(x)d|\rho_{{\bf
n},j}|(x)=0\,, \qquad \nu=0,\ldots,N_{{\bf n}^l,j} -1\,,
\]
where
\[ g_{{\bf n},j}(x) = \frac{|Q_{{\bf
n},j-1}(x)Q_{{\bf n},j+1}(x)|}{|Q_{{\bf n}^l,j-1}(x)Q_{{\bf
n}^l,j+1}(x)|}\frac{|h_{{\bf n}^l,j}(x)|}{|h_{{\bf n},j}(x)|}\,,
\quad d\rho_{{\bf n},j}(x) = \frac{h_{{\bf
n},j}(x)d\sigma_{j}(x)}{Q_{{\bf n},j-1}(x)Q_{{\bf n},j+1}(x)}\,.
\]
From (\ref{eq:h}) and (\ref{eq26})
\begin{equation} \label{eq:var}
\lim_{{\bf n} \in \Lambda'} g_{{\bf n},j}(x) =
 {|(\widetilde{F}_{j-1}^{(l)}\widetilde{F}_{j+1}^{(l)})(x)|^{-1}}
\end{equation}
uniformly on $\Delta_j$.

Fix $j \in \{-m_2,\ldots,-l_2-1\}\cup\{l_1+1,\ldots,m_1\}$. As
mentioned above, for this selection of $j$ we have that $\deg
Q_{{\bf n}^l,j} = \deg Q_{{\bf n},j} = N_{{\bf n},j}$. Due to
(\ref{eq:var}) and (\ref{eq26}), from Theorems 1 and  2 of
\cite{BCL}, it follows that
\begin{equation} \label{eq:otraa}
\lim_{{\bf n}\in \Lambda^{\prime}} \frac{Q_{{\bf
n}^l,j}(z)}{Q_{{\bf n},j}(z)}=\frac{S_j(z)}{S_j(\infty)} =
\widetilde{S}_j(z) = \widetilde{F}_j^{(l)}(z)\,, \qquad
\mathcal{K} \subset \overline{\mathbb{C}}\setminus
\sop{(\sigma_j)} \,,
\end{equation}
where $S_j$ is  the Szeg\H{o} function on $\overline{\mathbb{C}}
\setminus \widetilde{\Delta}_j$ with respect to
${|\widetilde{F}_{j-1}^{(l)}(x)\widetilde{F}_{j+1}^{(l)}(x)|^{-1}},
x \in \widetilde{\Delta}_j.$  The function $S_j$ is uniquely
determined by
\begin{equation}       \label{eq:Se}
\aligned &1)\qquad S_j,1/S_j \in H(\overline{\mathbb{C}} \setminus
\widetilde{\Delta}_j)\,,
\\
&2)\qquad  S_j(\infty)>0\,,
\\
&3 )\qquad |S_j(x)|^2
\frac1{\bigl|(\widetilde{F}_{j-1}^{(l)}\widetilde{F}_{j+1}^{(l)})(x)\bigr|}=1,
\qquad x \in \widetilde{\Delta}_j \,.
\endaligned
\end{equation}

Now, fix $j \in \{-l_2, \ldots , l_1\}$. In this situation $\deg
Q_{{\bf n}^l,j} = \deg Q_{{\bf n},j} +1 = N_{{\bf n},j} +1.$ Let
$Q_{{\bf n},j}^*(x)$ be the monic polynomial of degree $N_{{\bf
n},j}$ orthogonal with respect to the varying measure $g_{{\bf
n},j}\,d|\rho_{{\bf n},j}|$. Using the same arguments as above, we
have
\begin{equation}          \label{eq:comparativa}
\lim_{{\bf n}\in \Lambda^{\prime}} \frac{Q_{{\bf
n},j}^*(z)}{Q_{{\bf n},j}(z)}=\frac{S_j(z)}{S_j(\infty)} =
\widetilde{S}_j(z)\,, \qquad \mathcal{K} \subset
\overline{\mathbb{C}}\setminus \sop{(\sigma_{j})} \,.
\end{equation}
On the other hand, since $\deg Q_{{\bf n}^l,j}=\deg Q_{{\bf
n},j}^* +1$ and both of these polynomials are orthogonal with
respect to the same varying weight,  by  Theorem 1 of \cite{BCL}
and (\ref{eq:h}), it follows that
\begin{equation}                \label{eq:ratio}
\lim_{{\bf n}\in {\Lambda^{\prime}}} \frac{Q_{{\bf
n}^{l},j}(z)}{Q_{{\bf n},j}^*(z)}=
\frac{\varphi_j(z)}{{\varphi_j^{\prime}}(\infty)} =
\widetilde{\varphi}_j(z)\,, \qquad \mathcal{K} \subset
\mathbb{C}\setminus \sop{(\sigma_{j})}\,,
\end{equation}
where $\varphi_j$ denotes the conformal representation of
$\overline{\mathbb{C}}\setminus \widetilde{\Delta}_j$ onto $\{w :
|w|
> 1\}$ such that $\varphi_j(\infty) = \infty$ and
$\varphi_j^{\prime}(\infty)
> 0$. The function $\varphi_j$ is uniquely determined by
\begin{equation}    \label{eq:Zh}
   \aligned &1 )\qquad \varphi_j,1/\varphi_j \in
H(\mathbb{C} \setminus \widetilde{\Delta}_j)\,,
\\
&2)\qquad  \varphi_j^{\prime}(\infty)>0\,,
\\
&3 )\qquad |\varphi_j(x)|=1, \quad  x \in \widetilde{\Delta}_j \,.
\endaligned
\end{equation}
From (\ref{eq26}), (\ref{eq:comparativa}), and (\ref{eq:ratio}),
we obtain
\begin{equation} \label{eq:otra}
\lim_{{\bf n}\in \Lambda^{\prime}} \frac{Q_{{\bf
n}^l,j}(z)}{Q_{{\bf n},j}(z)}=
(\widetilde{S}_j\widetilde{\varphi}_j)(z) =
\widetilde{F}_j^{(l)}(z)\,, \quad \mathcal{K} \subset
\mathbb{C}\setminus \sop{(\sigma_j)}  \,.
\end{equation}

Thus,
\begin{equation}    \label{eq:FS}
\widetilde{F}_j^{(l)}= \left\{
\begin{array}{ll} \widetilde{S}_j \widetilde{\varphi}_j\,, &  j \in \{ -l_2,
\ldots , l_1\}\,, \\
\widetilde{S}_j \,, &  j \in
\{-m_2,\ldots,-l_2-1\}\cup\{l_1+1,\ldots,m_1\}\,,
\end{array}
\right.
\end{equation}
and  from (\ref{eq:Se}) and (\ref{eq:FS}) it follows that
\begin{equation} \label{eq:sysinter}
|\widetilde{F}_j^{(l)}(x)|^2
\frac1{\bigl|(\widetilde{F}_{j-1}^{(l)}\widetilde{F}_{j+1}^{(l)})(x)\bigr|}=
\frac{1}{\omega_j}
 \,, \qquad x \in \widetilde{\Delta}_j\,, \qquad j=-m_2,\ldots,m_1\,,
\end{equation}
where
\begin{equation} \label{eq:xd} \omega_j = \left\{
\begin{array}{ll}
(S_j\,\varphi_j^{\prime})^2(\infty)\,, & j \in
\{-l_2,\ldots,l_1\}\,, \\
S_j^2(\infty)\,, & j \in
\{-m_2,\ldots,-l_2-1\}\cup\{l_1+1,\ldots,m_1\} \,.
\end{array}
\right.
\end{equation}

Now, let us show that there exist positive constants $c_j,
j=-m_2,\ldots,m_1,$ such that the functions ${F_j^{(l)}} = c_j
\widetilde{F}_j^{(l)}$ satisfy (\ref{sisec}). In fact, according
to (\ref{eq:sysinter}) for any such constants $c_j$ we have that
\[ |{F_j^{(l)}}(x)|^2
\frac1{\bigl|({F_{j-1}^{(l)}}{F_{j+1}^{(l)}})(x)\bigr|}=
\frac{c_j^2}{c_{j-1}c_{j+1}\omega_j}
 \,, \qquad x \in \widetilde{\Delta}_j\,,\qquad j=-m_2,\ldots,m_1\,,
\]
where $c_{-m_2-1} = c_{m_1+1} = 1$. The problem reduces to finding
appropriate constants $c_j$ such that
\begin{equation} \label{eq:aa}
\frac{c_j^2}{c_{j-1}c_{j+1}\omega_j} = 1\,, \qquad j
=-m_2,\ldots,m_1\,.
\end{equation}
Taking logarithm, we obtain the linear system of equations
\begin{equation} \label{eq:bb} 2\log c_j - \log c_{j-1} -
\log c_{j+1} = \log \omega_j \,, \qquad j=-m_2,\ldots,m_1
\end{equation}
$(c_{-m_2-1} = c_{m_1+1} = 1)$ on the unknowns $\log c_j\,.$ This
system has a unique solution with which we conclude the proof.
\end{proof}

Consider the $(m_1+m_2+2)$-sheeted Riemann surface
$$
\mathcal R=\overline{\bigcup_{k=-m_2-1}^{m_1} \mathcal R_k} ,
$$
formed by the consecutively ``glued'' sheets
$$
\mathcal R_{-m_2-1}:=\overline {\mathbb{C}} \setminus
\widetilde{\Delta}_{-m_2},\,\,\, \mathcal R_k:=\overline
{\mathbb{C}} \setminus (\widetilde{\Delta}_k \cup
\widetilde{\Delta}_{k+1}),\,\, k=-m_2,\dots,m_1-1,\,\,\, \mathcal
R_{m_1}:=\overline {\mathbb{C}} \setminus
\widetilde{\Delta}_{m_1},
$$
where the upper and lower banks of the slits of two neighboring
sheets are identified. Fix $l = (l_1,l_2), 0 \leq l_1 \leq m_1, 0
\leq l_2 \leq m_2$. Let $\psi^{(l)}$ be a singled valued function
defined on $\mathcal{R}$ onto the extended complex plane
satisfying
\[ \psi^{(l)}(z)=\frac{C_{1}}{z}+\mathcal{O}(\frac{1}{z^{2}}),\quad
z\rightarrow\infty ^{(-l_2-1)},\]
\[ \psi^{(l)}(z)=C_{2}\,z+\mathcal{O}(1),\quad z\rightarrow\infty ^{(l_1)},\]
where $C_{1}$ and $C_{2}$ are nonzero constants. Since the genus
of $\mathcal{R}$ is zero, $\psi^{(l)}$ exists and is uniquely
determined up to a multiplicative constant. Consider the branches
of $\psi^{(l)}$, corresponding to the different sheets
$k=-m_2-1,\ldots,m_1$ of $\mathcal{R}$
\[\psi^{(l)}:=\{\psi^{(l)}_{k}\}_{k=-m_2-1}^{m_1}\,. \]
We normalize $\psi^{(l)}$ so that
\begin{equation} \label{mormaliz}
\prod_{k=-m_2-1}^{m_1}\,|\psi^{(l)}_{k}(\infty)|=1, \qquad C_1 \in
\mathbb{R}\setminus \{0\}.
\end{equation}
Certainly, there are two $\psi^{(l)}$ verifying this
normalization. Since the product of all the branches
$\prod_{k=-m_2-1}^{m_1} \psi^{(l)}_{k}$ is a single valued
analytic function in $\overline{\mathbb{C}}$ without
singularities, by Liouville's theorem it is constant and because
of the normalization introduced above this constant is either $1$
or $-1$.

If $\psi^{(l)}$ is such that $C_1 \in \mathbb{R}\setminus \{0\},$
then
\[ \psi^{(l)}(z) = \overline{\psi^{(l)}(\overline{z})}, \qquad z \in
\mathcal{R}.
\]
In fact, let $\phi(z) := \overline{\psi^{(l)}(\overline{z})}$.
$\phi$ and $\psi^{(l)}$ have the same divisor; consequently, there
exists a constant $C$ such that $\phi= C\psi^{(l)}$. Comparing the
leading coefficients of the Laurent expansion of these functions
at $\infty^{(-l_2-1)}$,  we conclude that $C=1$ since $C_1 \in
\mathbb{R} \setminus \{0\}$.

In terms of the branches of $\psi^{(l)}$, the symmetry formula
above means that  for each  $k= -m_2-1,\ldots,m_1$,
\[
\psi^{(l)}_k: \overline{\mathbb{R}} \setminus
(\widetilde{\Delta}_k\cup \widetilde{\Delta}_{k+1})
\longrightarrow \overline{\mathbb{R}}
\]
$(\widetilde{\Delta}_{-m_2-1} =
\widetilde{\Delta}_{m_1+1}=\emptyset)$; therefore, the
coefficients (in particular, the leading one) of the Laurent
expansion at $\infty$ of these branches are real numbers, and
\begin{equation}  \label{contacto}\psi^{(l)}_k(x_{\pm}) =
\overline{\psi^{(l)}_k(x_{\mp})} =
\overline{\psi^{(l)}_{k+1}(x_{\pm})}, \qquad x \in
\widetilde{\Delta}_{k+1}. \end{equation}

Given an arbitrary function $F(z)$  which has in a neighborhood of
infinity a Laurent expansion of the form $F(z)= Cz^k +
{\mathcal{O}}(z^{k-1}), C \neq 0,$ and $ k \in {\mathbb{Z}},$ we
denote
\[
\widetilde{F}:= {F}/{C}\,.
\]
$C$ is called the leading coefficient of $F$. When $C \in
\mathbb{R}$, $\mbox{sg}(F(\infty))$ will represent the sign of
$C$.

We are ready to state and prove one of the main results of this
section.

\begin{thm} \label{teofund}
Let $S^1= \mathcal{N}'(\sigma^1_0,\ldots,\sigma^1_{m_1}), S^2=
\mathcal{N}'(\sigma^2_0,\ldots,\sigma^2_{m_2})$ be given, and let
${\bf n} \in \Lambda \subset \mathbb{Z}_+^{m_1+1}(\bullet)\times
\mathbb{Z}_+^{m_2+1}(\bullet)$ be an infinite sequence of distinct
multi-indices such that
\[
\displaystyle{\sup_{{\bf n} \in \Lambda }((m_2 +1) n_{2,0} - |{\bf
n}_2| ) }< \infty, \qquad \displaystyle{\sup_{{\bf n} \in \Lambda
}((m_1+1)n_{1,0} -  |{\bf n}_1| ) }< \infty.
\]
Let us assume that there exists $l=(l_{1};l_{2}), 0\leq l_{1}\leq
m_{1}, 0\leq l_{2}\leq m_{2},$ such that for all ${\bf n} \in
\Lambda$ we have that ${\bf n}^l \in
\mathbb{Z}_{+}^{m_1+1}(\bullet)\times
\mathbb{Z}_{+}^{m_2+1}(\bullet)$. Let $\{Q_{{\bf
n},j}\}_{j=-m_2}^{m_1}$, ${\bf n}\in\Lambda$, be the corresponding
sequences of polynomials defined in section \ref{normal}. Then,
for each fixed $j\in\{-m_2,\ldots,m_1\}$, we have
\begin{equation}\label{eq30}
\lim_{{\bf n}\in\Lambda}\,\frac{Q_{{\bf n}^l,j}(z)}{Q_{{\bf
n},j}(z)}=\widetilde{F}^{(l)}_{j}(z),\qquad z\in \mathcal{K}
\subset\mathbb{C}\setminus\supp{\sigma_j}
\end{equation}
where the functions satisfying $(\ref{sisec})$ are
\begin{equation}\label{eq31}
F^{(l)}_{j}:=\mbox{sg}\left(\prod_{k=j}^{m_1}\,\psi^{(l)}_{k}(\infty)\right)
\prod_{k=j}^{m_1}\,\psi^{(l)}_{k}.
\end{equation}
\end{thm}
\begin{proof} Since the families of functions
\[ \left\{ {Q_{{\bf n}^l,j}}/{Q_{{\bf n},j}}\right\}_{{\bf n}\in \Lambda}
\,,\qquad j=-m_2,\ldots,m_1,
\]
are uniformly bounded on each compact subset $\mathcal{K} \subset
\mathbb{C} \setminus \sop{(\sigma_j)}$ for all sufficiently large
$|{\bf n}_{1}|, {\bf n} \in \Lambda$,  uniform convergence on
compact subsets of the indicated region follows from proving that
any convergent subsequence  has the same limit. According to Lemma
\ref{sistema} the limit functions, appropriately normalized, of a
convergent subsequence satisfy the same system of boundary value
problems (\ref{sisec}). According to Lemma 4.2 in \cite{AptLopRoc}
this system has a unique solution.

It remains to show that the functions defined in (\ref{eq31})
satisfy (\ref{sisec}). When multiplying two consecutive branches,
the singularities on the common slit cancel out by the Schwarz
reflection principle; therefore, $1)$ takes place since only the
singularities of $\psi_j^{(l)}$ on $\widetilde{\Delta}_j$ remain.
From the definition of $\psi^{(l)}$ it also follows that for
$j=-l_2,\ldots,l_1$, $F_j^{(l)}$ has at infinity a simple pole,
whereas it is regular and different from zero at infinity when $j
\in \{-m_2,\ldots,-l_2-1\}\cup\{l_1+1,\ldots,m_1\}$. The factor
sign in front of (\ref{eq31}) guarantees the positivity claimed in
$2)$ and $2^{\prime})$.

In order to verify $3)$, notice that $F_j^{(l)}/F_{j-1}^{(l)} =
\mbox{sg} (\psi_{j-1}^{(l)}(\infty))/\psi_{j-1}^{(l)}$. Therefore,
if $j=-m_2+1,\ldots,m_1,$
\[ \frac{|F_{j}^{(l)}(x)|^2}{|F_{j-1}^{(l)}(x)F_{j+1}^{(l)}(x)|} =
\frac{|\psi_j^{(l)}(x)|}{|\psi_{j-1}^{(l)}(x)|} = 1, \qquad x \in
\widetilde{\Delta}_j,
\]
on account of (\ref{contacto}). For $j=-m_2$, from the definition
and (\ref{contacto})
\[ \frac{|F_{-m_2}^{(l)}(x) |^2}{| F_{-m_2+1}^{(l)}(x) |} =
|\psi_{-m_2}^{(l)}(x)||\prod_{k=-m_2}^{m_1}\,\psi^{(l)}_{k}(x)| =
|\prod_{k=-m_2-1}^{m_1}\,\psi^{(l)}_{k}(x)| =1, \qquad x \in
\widetilde{\Delta}_{-m_2},
\]
since $\prod_{k=-m_2-1}^{m_1}\,\psi^{(l)}_{k}$ is constantly equal
to $1$ or $-1$ on all $\overline{\mathbb{C}}.$ The proof is
complete.
\end{proof}

The following corollary complements Theorem \ref{teofund}. The
proof is similar to that of Corollary 4.1 in \cite{AptLopRoc}.

\begin{cor} \label{cor3}
Let $S^1= \mathcal{N}'(\sigma^1_0,\ldots,\sigma^1_{m_1}), S^2=
\mathcal{N}'(\sigma^2_0,\ldots,\sigma^2_{m_2})$ be given, and let
$\Lambda \subset \mathbb{Z}_+^{m_1+1}(\bullet)\times
\mathbb{Z}_+^{m_2+1}(\bullet)$ be an infinite sequence of distinct
multi-indices such that
\[
\displaystyle{\sup_{{\bf n} \in \Lambda }((m_2 +1) n_{2,0} - |{\bf
n}_2| ) }< \infty, \qquad \displaystyle{\sup_{{\bf n} \in \Lambda
}((m_1+1)n_{1,0} -  |{\bf n}_1| ) }< \infty.
\]
Let us assume that there exists $l=(l_{1};l_{2}), 0\leq l_{1}\leq
m_{1}, 0\leq l_{2}\leq m_{2},$ such that for all ${\bf n} \in
\Lambda$ we have that ${\bf n}^l \in
\mathbb{Z}_{+}^{m_1+1}(\bullet)\times
\mathbb{Z}_{+}^{m_2+1}(\bullet)$. Let $\{q_{{\bf n},j} =
\kappa_{{\bf n},j}\,Q_{{\bf n},j} \}_{j=-m_2}^{m_1} ,{\bf n}\in
{\Lambda},$ be the system of orthonormal polynomials as defined in
$(\ref{eq:orton})$ and $\{K_{{\bf n},j}\}_{j=-m_2}^{m_1} ,{\bf
n}\in {\Lambda},$ the values given by $(\ref{eq:K})$. Then, for
each fixed $j = -m_2,\ldots,m_1,$ we have
\begin{equation} \label{eq:xa}\lim_{{\bf n}\in
{\Lambda}}\frac{\kappa_{{\bf n}^l,j}}{\kappa_{{\bf n},j}}=
\kappa^{(l)}_j\,,
\end{equation}
\begin{equation} \label{eq:xk}
\lim_{{\bf n}\in {\Lambda}}\frac{K_{{\bf n}^l,j}}{K_{{\bf n},j}}=
\kappa^{(l)}_{j}\cdots\kappa^{(l)}_{m_1}\,,
\end{equation}
and
\begin{equation} \label{eq:xb}
\lim_{{\bf n}\in {\Lambda}}\frac{q_{{\bf n}^l,j}(z)}{q_{{\bf
n},j}(z)}= \kappa^{(l)}_j \widetilde{F}_j^{(l)}(z), \qquad
\mathcal{K} \subset {\mathbb{C}} \setminus \supp{\sigma_j} \,,
\end{equation}
where
\begin{equation} \label{eq:xc} \kappa^{(l)}_j =
\frac{c_{j}^{(l)}}{\sqrt{c_{j-1}^{(l)}c_{j+1}^{(l)}}}\,, \qquad
c_{j}^{(l)} = \left\{
\begin{array}{ll}
(F^{(l)}_j)^{\prime}(\infty)\,, & j \in \{-l_2,\ldots,l_1\} \,, \\
F^{(l)}_j(\infty)\,, & j  \notin \{-l_{2},\ldots,l_{1}\}\,,
\end{array}
\right.
\end{equation}
$(c_{-m_{2}-1}^{(l)}=c_{m_{1}+1}^{(l)}=1)$ and the functions
$F^{(l)}_j$ are defined by $(\ref{eq31})$.
\end{cor}

\begin{proof} By  Theorem  \ref{teofund}, we have limit in
(\ref{eq:var}) along the whole sequence $\Lambda$. Reasoning as in
the deduction of formulas (\ref{eq:otraa}) and (\ref{eq:otra}),
but now in connection with orthonormal polynomials (see Theorems 1
and 2 of \cite{BCL}),   it follows that
\[
\lim_{{\bf n}\in \Lambda} \frac{q_{{\bf n}^l,j}(z)}{q_{{\bf
n},j}(z)}= \left\{
\begin{array}{ll}
(S_j\,\varphi_j)(z) \,, & j \in \{-l_2,\ldots,l_1\}\,, \\
S_j(z) \,, & j \in \{-m_2,\ldots,-l_2-1\} \cup
\{l_1+1,\ldots,m_1\}\,,
\end{array}
\right.
\]
uniformly on compact subsets of $\mathbb{C}\setminus
\sop{(\sigma_j)}$, where $S_j$ is defined in (\ref{eq:Se}). This
formula, divided by (\ref{eq:otraa}) or (\ref{eq:otra}) according
to the value of $j$ gives
\[ \lim_{{\bf n}\in {\Lambda}}\frac{\kappa_{{\bf n}^l,j}}{\kappa_{{\bf n},j}}=
\sqrt{\omega_j} = \frac{c_j}{\sqrt{c_{j-1}c_{j+1}}}\,,
\]
where $\omega_j$ is defined in (\ref{eq:xd}), and the $c_j$ are
the normalizing constants found in Lemma \ref{sistema} solving the
linear system of equations (\ref{eq:bb}) which ensure that
\[ F_j^{(l)} \equiv c_j \widetilde{F}_j^{(l)}\,, \qquad j=-m_2,\ldots,m_1\,,
\]
with $F_j^{(l)}$ satisfying (\ref{sisec}) and thus given by
(\ref{eq31}). Since $(\widetilde{F}_j^{(l)})^{\prime}(\infty) = 1,
j\in \{-l_2,\ldots,l_1\},$ and $(\widetilde{F}_j^{(l)})(\infty) =
1, j \in \{-m_2,\ldots,-l_2-1\} \cup \{l_1+1,\ldots,m_1\}$ formula
(\ref{eq:xa}) immediately follows with $\kappa_j^{(l)}$ as in
(\ref{eq:xc}).

From the definition of $\kappa_{{\bf n},j}\,,$ we have that
\[ K_{{\bf n},j} = \kappa_{{\bf n},j}\cdots\kappa_{{\bf n},m_1} \,.
\]
Taking the ratio of these constants for the multi-indices ${\bf
n}$ and ${\bf n}^l$ and using (\ref{eq:xa}), we get (\ref{eq:xk}).
Formula (\ref{eq:xb}) is an immediate consequence of (\ref{eq:xa})
and (\ref{eq30}).
\end{proof}

Let $\lcm{(a,b)}$ denote the least common multiple of two integers
$a$ and $b$, and define $m:=\lcm{(m_{1}+1,m_{2}+1)}$,
$d_{1}:=m/(m_{1}+1)$, $d_{2}:=m/(m_{2}+1)$. Within the class of
pairs $l=(l_{1};l_{2})$ with $0\leq l_{1}\leq m_{1}$, $0\leq
l_{2}\leq m_{2}$, we distinguish the subclass
\[L:=\{(l_{1};l_{2}): l_{1}\equiv r\,\text{mod}\,(m_{1}+1),\,\,
l_{2}\equiv r\,\text{mod}\,(m_{2}+1)\,\,\mbox{for some}\,\,0\leq
r\leq m-1\}\,.
\]
It is easy to check that for different $r, 0\leq r\leq m-1,$ the
pairs $(l_1,l_2)$ in $L$ are different. Let
$\mathbf{p}:=(\mathbf{p}_{1};\mathbf{p}_{2})$, where
$\mathbf{p}_{1}=(d_{1},\ldots,d_{1})$ and
$\mathbf{p}_{2}=(d_{2},\ldots,d_{2})$ have $m_{1}+1$ and $m_{2}+1$
components, respectively. By ${\bf n}+{\bf p}$ we denote the
multi-index
$(\mathbf{n}_{1}+\mathbf{p}_{1};\mathbf{n}_{2}+\mathbf{p}_{2})$.

\begin{cor} \label{corperiod1}
Let $S^1= \mathcal{N}'(\sigma^1_0,\ldots,\sigma^1_{m_1}), S^2=
\mathcal{N}'(\sigma^2_0,\ldots,\sigma^2_{m_2})$ be given, and let
$\Lambda \subset \mathbb{Z}_+^{m_1+1}(\bullet)\times
\mathbb{Z}_+^{m_2+1}(\bullet)$ be an infinite sequence of distinct
multi-indices such that
\[
\displaystyle{\sup_{{\bf n} \in \Lambda }((m_2 +1) n_{2,0} - |{\bf
n}_2| ) }< \infty, \qquad \displaystyle{\sup_{{\bf n} \in \Lambda
}((m_1+1)n_{1,0} -  |{\bf n}_1| ) }< \infty.
\]
Then, for each fixed $j\in\{-m_2,\ldots,m_1\}$, we have
\begin{equation}\label{eqperiod1}
\lim_{{\bf n}\in\Lambda}\,\frac{Q_{{\bf n}+{\bf p},j}(z)}{Q_{{\bf
n},j}(z)}=\prod_{l\in L}\widetilde{F}^{(l)}_{j}(z),\qquad
\mathcal{K} \subset\mathbb{C}\setminus\supp{\sigma_{j}}\,,
\end{equation}
\begin{equation}\label{eqperiod12}
\lim_{{\bf n}\in\Lambda}\,\frac{\kappa_{{\bf n}+{\bf
p},j}}{\kappa_{{\bf n},j}}=\prod_{l\in L}\kappa^{(l)}_{j}\,,
\end{equation}
and
\begin{equation}\label{eqperiod13}
\lim_{{\bf n}\in\Lambda}\,\frac{q_{{\bf n}+{\bf p},j}(z)}{q_{{\bf
n},j}(z)}=\prod_{l\in
L}\kappa^{(l)}_{j}\widetilde{F}^{(l)}_{j}(z),\qquad \mathcal{K}
\subset\mathbb{C}\setminus\supp{\sigma_{j}}\,.
\end{equation}
\end{cor}
\begin{proof}
Given ${\bf n}\in\Lambda$ and $0\leq r\leq m$, let
$\mathbf{n}(r):=\mathbf{n}+\mathbf{q}(r)$ where
$\mathbf{q}(r)=(\mathbf{q}_{1}(r);\mathbf{q}_{2}(r))$ is the
multi-index satisfying
\[
\mathbf{q}_{i}(r)=(\underbrace{k+1,\ldots,k+1}_{s},k,\ldots,k)\,,\quad
r=k(m_{i}+1)+s\,,\,\,0\leq s\leq m_{i}\,.
\]
Hence, $\mathbf{n}(0)=\mathbf{n}$, $\mathbf{n}(m)={\bf n}+{\bf p}$
and $\mathbf{n}(r)\in\mathbb{Z}_+^{m_1+1}(\bullet)\times
\mathbb{Z}_+^{m_2+1}(\bullet)$ for every $r$.

We have
\[
\frac{Q_{{\bf n}+{\bf p},j}(z)}{Q_{{\bf
n},j}(z)}=\prod_{r=0}^{m-1}\frac{Q_{{\bf n}(r+1),j}(z)}{Q_{{\bf
n}(r),j}(z)}\,.
\]
In addition, by (\ref{eq30}) we know that for each fixed $0\leq
r\leq m-1$,
\[
\lim_{{\bf n}\in\Lambda}\,\frac{Q_{{\bf n}(r+1),j}(z)}{Q_{{\bf
n}(r),j}(z)}=\widetilde{F}^{(l)}_{j}(z),\qquad z\in \mathcal{K}
\subset\mathbb{C}\setminus\supp{\sigma_j}\,,
\]
where $l=(l_{1};l_{2})$ is precisely the multi-index satisfying
$l_{1}\equiv r\,\text{mod}\,(m_{1}+1)$, $l_{2}\equiv
r\,\text{mod}\,(m_{2}+1)$. Therefore (\ref{eqperiod1}) follows.
Relations (\ref{eqperiod12}) and (\ref{eqperiod13}) are proved
analogously in view of (\ref{eq:xa}) and  (\ref{eq:xb}).
\end{proof}

Now, we need to introduce some notations. For
$j\in\{-m_{2},\ldots,m_{1}-1\}$, set
\[
\delta_{j}:=\left\{
\begin{array}{ccccc}
1, & \mbox{if} & \Delta_{j}
& \mbox{is to the left of} & \Delta_{j+1} \,,\\
-1, & \mbox{if} & \Delta_{j}
& \mbox{is to the right of} & \Delta_{j+1} \,.\\
\end{array}
\right.
\]
For multi-indices $l=(l_{1};l_{2})$ such that $l_{1}+l_{2}\geq 2$
we define
\[
\Delta_{j,\,l}:=\left\{
\begin{array}{ccc}
1, & \mbox{if} & j\geq l_{1}+2 \,,\\
\delta_{j-1}, & \mbox{if} & j\in\{l_{1},l_{1}+1\} \,,\\
-\delta_{j-1}\delta_{j}, & \mbox{if} & j\in\{-l_{2}+1,\ldots,l_{1}-1\}\,,\\
-\delta_{j}, & \mbox{if} & j\in\{-l_{2}-1, -l_{2}\} \,,\\
1, & \mbox{if} & j\leq -l_{2}-2 \,.
\end{array}
\right.
\]
If $l_{1}+l_{2}=1$ then
\[
\Delta_{j,\,l}:=\left\{
\begin{array}{ccc}
1, & \mbox{if} & j\geq l_{1}+2 \,,\\
\delta_{j-1}, & \mbox{if} & j\in\{l_{1},l_{1}+1\} \,,\\
-\delta_{j}, & \mbox{if} & j\in\{-l_{2}-1,-l_{2}\}\,,\\
1, & \mbox{if} & j\leq -l_{2}-2 \,,
\end{array}
\right.
\]
and for $l_{1}=l_{2}=0$
\[
\Delta_{j,(0;0)}:=\left\{
\begin{array}{ccc}
1, & \mbox{if} & j\geq 2 \,,\\
\delta_{0}, & \mbox{if} & j=1 \,,\\
1, & \mbox{if} & j=0\,,\\
-\delta_{-1}, & \mbox{if} & j=-1 \,,\\
1, & \mbox{if} & j\leq -2\,.
\end{array}
\right.
\]

Recall that $\varepsilon_{{\bf n},j}$ denotes the sign of the
varying measure
\[
d\rho_{{\bf n},j}(x)=\frac{h_{{\bf n},j}(x)d\sigma_j(x)}{Q_{{\bf
n},j-1}(x)Q_{{\bf n},j+1}(x)}\,.
\]
\begin{lem}\label{quotienteps}
For any ${\bf n}, {\bf n}^l
\in\mathbb{Z}_{+}^{m_1+1}(\bullet)\times
\mathbb{Z}_{+}^{m_2+1}(\bullet)$ and $-m_{2}\leq j\leq m_{1}$
\begin{equation}\label{ratioepsilons}
\frac{\varepsilon_{{\bf n}^{l},j}}{\varepsilon_{{\bf n},j}}
=\prod_{k=j}^{m_{1}}\Delta_{k,l}\,.
\end{equation}
\end{lem}
\begin{proof}
We will denote by $\sign{f,\Delta}$ the sign of a function $f$ on
the interval $\Delta$. Thus
\begin{equation}\label{eqepsilon1}
\frac{\varepsilon_{\mathbf{n}^{l},j}}{\varepsilon_{\mathbf{n},j}}
=\sign{\frac{\mathcal{H}_{\mathbf{n}^{l},j}\,Q_{\mathbf{n},j-1}\,Q_{\mathbf{n},j+1}}{
{\mathcal{H}_{\mathbf{n},j}\,Q_{\mathbf{n}^{l},j-1}}\,
Q_{\mathbf{n}^{l},j+1}},\Delta_{j}}\,.
\end{equation}
If $-l_{2}\leq j-1\leq l_{1},$ then
$\deg(Q_{\mathbf{n}^{l},j-1})=1+\deg(Q_{\mathbf{n},j-1})$ and,
therefore,
\begin{equation}\label{eq:14}
\sign{\frac{Q_{\mathbf{n},j-1}}{Q_{\mathbf{n}^{l},j-1}},\Delta_{j}}=\delta_{j-1}\,.
\end{equation}
If $j-1<-l_{2}$ or $j-1>l_{1}$, then
$\deg(Q_{\mathbf{n}^{l},j-1})=\deg(Q_{\mathbf{n},j-1})$, implying
that
\begin{equation}\label{eq:15}
\sign{\frac{Q_{\mathbf{n},j-1}}{Q_{\mathbf{n}^{l},j-1}},\Delta_{j}}=1\,.
\end{equation}
Analogously, we have that for $-l_{2}\leq j+1\leq l_{1}$
\begin{equation}\label{eq:16}
\sign{\frac{Q_{\mathbf{n},j+1}}{Q_{\mathbf{n}^{l},j+1}},\Delta_{j}}=-\delta_{j}
\end{equation}
and for $j+1<-l_{2}$ or $j+1>l_{1}$
\begin{equation}\label{eq:17}
\sign{\frac{Q_{\mathbf{n},j+1}}{Q_{\mathbf{n}^{l},j+1}},\Delta_{j}}=1\,.
\end{equation}
From (\ref{eq:14})-(\ref{eq:16}) it follows that
\begin{equation}\label{eq:18*}
\sign{\frac{Q_{\mathbf{n},j-1}\,Q_{\mathbf{n},j+1}}{
{Q_{\mathbf{n}^{l},j-1}}\,
Q_{\mathbf{n}^{l},j+1}},\Delta_{j}}=\Delta_{j,\,l}\,.
\end{equation}
Now, by (\ref{formint3})
\[\frac{\mathcal{H}_{\mathbf{n}^{l},j}(x)}{\mathcal{H}_{\mathbf{n},j}(x)}=
\frac{\int\frac{Q_{\mathbf{n}^{l},j+1}^{2}(t)}{x-t}\frac{\mathcal{H}_{\mathbf{n}^{l},j+1}(t)
\,d\sigma_{j+1}(t)}{Q_{\mathbf{n}^{l},j}(t)Q_{\mathbf{n}^{l},j+2}(t)}}{\int\frac{Q_{\mathbf{n},j+1}^{2}(t)}{x-t}
\frac{\mathcal{H}_{\mathbf{n},j+1}(t)
\,d\sigma_{j+1}(t)}{Q_{\mathbf{n},j}(t)Q_{\mathbf{n},j+2}(t)}}\,.\]
Therefore,
\begin{equation}\label{eqepsilon2}
\sign{
{\mathcal{H}_{\mathbf{n}^{l},j}}/{\mathcal{H}_{\mathbf{n},j}},
\Delta_{j}}=\frac{\varepsilon_{\mathbf{n}^{l},j+1}}
{\varepsilon_{\mathbf{n},j+1}}\,.
\end{equation}
Since $\mathcal{H}_{{\bf n}^{l},m_{1}}\equiv\mathcal{H}_{{\bf
n},m_{1}}\equiv 1$, the right hand side of (\ref{eqepsilon2})
equals $1$ for $j=m_{1}$. Hence (\ref{ratioepsilons}) follows from
(\ref{eqepsilon1}), (\ref{eq:18*}) and (\ref{eqepsilon2}).
\end{proof}

This lemma shows that $
{\varepsilon_{\mathbf{n}^{l},j}}/{\varepsilon_{\mathbf{n},j}}$
depends on $j,l,$ and the relative positions of the intervals
$\Delta_j$ but not on ${\bf n}$. Define the functions
\[
{\mathcal{A}}^{(l)}_{j}:=\widetilde{\psi}^{(l)}_{j}
\prod_{k=j+1}^{m_{1}}\frac{\Delta_{k,\,l}}{(\kappa^{(l)}_{k})^{2}}
\]
(the product should be understood to be equal to $1$ when
$j=m_1$).
\begin{thm} \label{teoLratio}
Let $S^1= \mathcal{N}'(\sigma^1_0,\ldots,\sigma^1_{m_1}), S^2=
\mathcal{N}'(\sigma^2_0,\ldots,\sigma^2_{m_2})$ be given, and let
$\Lambda \subset \mathbb{Z}_+^{m_1+1}(\bullet)\times
\mathbb{Z}_+^{m_2+1}(\bullet)$ be an infinite sequence of distinct
multi-indices such that
\[
\displaystyle{\sup_{{\bf n} \in \Lambda }((m_2 +1) n_{2,0} - |{\bf
n}_2| ) }< \infty\,, \qquad \displaystyle{\sup_{{\bf n} \in
\Lambda }((m_1+1)n_{1,0} -  |{\bf n}_1| ) }< \infty\,.
\]
Let us assume that there exists $l=(l_{1};l_{2}), 0\leq l_{1}\leq
m_{1}, 0\leq l_{2}\leq m_{2},$ such that for all ${\bf n} \in
\Lambda$ we have that ${\bf n}^l \in
\mathbb{Z}_{+}^{m_1+1}(\bullet)\times
\mathbb{Z}_{+}^{m_2+1}(\bullet)$. Let $\{\mathcal{A}_{{\bf
n},j}\}_{j=-m_2-1}^{m_1} ,{\bf n}\in {\Lambda},$ be the associated
sequences of ``monic'' linear forms. Then, for each fixed $j =
-m_2-1,\ldots,m_1,$
\begin{equation} \label{asympforms}
\lim_{{\bf n}\in {\Lambda}}\frac{\mathcal{A}_{{\bf
n}^{l},j}(z)}{\mathcal{A}_{{\bf n},j}(z)}=
{\mathcal{A}}^{(l)}_{j}\,, \qquad \mathcal{K} \subset {\mathbb{C}}
\setminus (\supp{\sigma_j}\cup\supp{\sigma_{j+1}})
\end{equation}
$(\supp{\sigma_{-m_{2}-1}}=\supp{\sigma_{m_{1}+1}}=\emptyset)$.
\end{thm}
\begin{proof}
It follows from definition of $\mathcal{H}_{{\bf n},j}$ and
$\mathcal{H}_{{\bf n}^{l},j}$ that
\[
\frac{\mathcal{A}_{{\bf n}^{l},j}(z)}{\mathcal{A}_{{\bf n},j}(z)}=
\frac{\varepsilon_{{\bf n}^{l},j+1}h_{{\bf
n}^{l},j}(z)}{\varepsilon_{{\bf n},j+1}h_{{\bf
n},j}(z)}\frac{\varepsilon_{{\bf n},j+1}}{\varepsilon_{{\bf
n}^{l},j+1}}\frac{K_{{\bf n},j+1}^{2}}{K_{{\bf n}^{l},j+1}^{2}}
\frac{Q_{{\bf n}^{l},j}(z)}{Q_{{\bf n},j}(z)}\frac{Q_{{\bf
n},j+1}(z)}{Q_{{\bf n}^{l},j+1}(z)}\,.
\]
By Lemma \ref{lm:3},
\[
\lim_{{\bf n}\in {\Lambda}}\frac{\varepsilon_{{\bf
n}^{l},j+1}h_{{\bf n}^{l},j}(z)}{\varepsilon_{{\bf n},j+1}h_{{\bf
n},j}(z)}=1\,,\qquad \mathcal{K} \subset {\mathbb{C}} \setminus
\sop{(\sigma_{j+1})}\,.
\]
Using Lemma \ref{quotienteps} and Corollary \ref{cor3}, we have
\[
\lim_{{\bf n}\in {\Lambda}}\frac{\varepsilon_{{\bf
n},j+1}}{\varepsilon_{{\bf n}^{l},j+1}}\frac{K_{{\bf
n},j+1}^{2}}{K_{{\bf
n}^{l},j+1}^{2}}=\prod_{k=j+1}^{m_{1}}\frac{\Delta_{k,\,l}}{(\kappa^{(l)}_{k})^{2}}\,.
\]
Finally, applying (\ref{eq30}) and (\ref{eq31}) one obtains
\[
\lim_{{\bf n}\in {\Lambda}}\frac{Q_{{\bf n}^{l},j}(z)}{Q_{{\bf
n},j}(z)}\frac{Q_{{\bf n},j+1}(z)}{Q_{{\bf
n}^{l},j+1}(z)}=\widetilde{\psi}_{j}^{(l)}(z)\,,\qquad \mathcal{K}
\subset {\mathbb{C}} \setminus
(\sop{(\sigma_j)}\cup\sop{(\sigma_{j+1})})\,.
\]
Putting these relations together we get (\ref{asympforms}).
\end{proof}

\begin{cor} \label{corperiod2}
Let $S^1= \mathcal{N}'(\sigma^1_0,\ldots,\sigma^1_{m_1}), S^2=
\mathcal{N}'(\sigma^2_0,\ldots,\sigma^2_{m_2})$ be given, and let
$\Lambda\subset \mathbb{Z}_+^{m_1+1}(\bullet)\times
\mathbb{Z}_+^{m_2+1}(\bullet)$ be an infinite sequence of distinct
multi-indices such that
\begin{equation} \label{eq:m}
\displaystyle{\sup_{{\bf n} \in \Lambda }((m_2 +1) n_{2,0} - |{\bf
n}_2| ) }< \infty\,, \qquad \displaystyle{\sup_{{\bf n} \in
\Lambda }((m_1+1)n_{1,0} -  |{\bf n}_1| ) }< \infty\,.
\end{equation}
Then, for each fixed $j\in\{-m_2,\ldots,m_1\}$, we have
\begin{equation}\label{eqperiod2}
\lim_{{\bf n}\in\Lambda}\,\frac{\mathcal{A}_{{\bf n}+{\bf
p},j}(z)}{\mathcal{A}_{{\bf n},j}(z)}=\prod_{l\in
L}\mathcal{A}^{(l)}_{j}(z)\,,\qquad \mathcal{K} \subset
{\mathbb{C}} \setminus (\supp{\sigma_j}\cup\supp{\sigma_{j+1}})
\end{equation}
$(\supp{\sigma_{-m_{2}-1}}=\supp{\sigma_{m_{1}+1}}=\emptyset)$.
Consequently,
\begin{equation} \label{eq:98}
\lim_{{\bf n}\in\Lambda} \left|\mathcal{A}_{{\bf
n},j}(z)\right|^{1/|{\bf n}_1|}  = \prod_{l\in
L}|\mathcal{A}^{(l)}_{j}(z) |^{1/m}\,,\qquad \mathcal{K} \subset
{\mathbb{C}} \setminus (\supp{\sigma_j}\cup\supp{\sigma_{j+1}}),
\end{equation}
where $m = \mbox{\rm lcm}(m_1+1,m_2+1).$
\end{cor}

\begin{proof} Using the same arguments employed to prove Corollary
\ref{corperiod1}, we obtain (\ref{eqperiod2}). From
(\ref{eqperiod2}) it is easy to deduce  the $|{\bf n}_1|$-th root
asymptotic of the linear forms.

In fact, it is easy to see that for each ${\bf n} \in \Lambda$
there exists ${\bf n}_0 \in \mathbb{Z}_+^{m_1+1}(\bullet)\times
\mathbb{Z}_+^{m_2+1}(\bullet)$ (which may depend on ${\bf n}$),
whose entries are uniformly bounded by a constant $C$ independent
of ${\bf n}$ (condition (\ref{eq:m}) is used), such that ${\bf n}
= r{\bf p} + {\bf n}_0$ for some $r \in {\mathbb{Z}}_+$. Write
\[ \mathcal{A}_{{\bf n},j}(z) = \frac{\mathcal{A}_{{\bf n},j}(z)}{\mathcal{A}_{{\bf n}-{\bf
p},j}(z)} \frac{\mathcal{A}_{{\bf n}-{\bf
p},j}(z)}{\mathcal{A}_{{\bf n}-2{\bf p},j}(z)}\cdots
\frac{\mathcal{A}_{{\bf n}_0+{\bf p},j}(z)}{\mathcal{A}_{{\bf
n}_0,j}(z)} \mathcal{A}_{{\bf n}_0,j}(z).
\]
Then
\[ \frac{1}{|{\bf n}_1|}\log |\mathcal{A}_{{\bf n},j}(z)| = \frac{1}{|{\bf n}_1|} \log |\mathcal{A}_{{\bf
n}_0,j}(z)| + \frac{1}{|{\bf n}_1|}\sum_{k=0}^{r-1} \log \left|
\frac{\mathcal{A}_{{\bf n}_0 + (k+1){\bf
p},j}(z)}{\mathcal{A}_{{\bf n}_0 + k{\bf p},j}(z)}\right|.
\]
Obviously,
\[ \lim_{{\bf n} \in \Lambda} \frac{1}{|{\bf n}_1|}\log |\mathcal{A}_{{\bf
n}_0,j}(z)| = 0,\qquad \mathcal{K} \subset {\mathbb{C}} \setminus
(\supp{\sigma_j}\cup\supp{\sigma_{j+1}}),
\]
and because of (\ref{eqperiod2})
\[\lim_{{\bf n} \in \Lambda} \frac{1}{|{\bf n}_1|}\sum_{k=0}^{r-1} \log \left|
\frac{\mathcal{A}_{{\bf n}_0 + (k+1){\bf
p},j}(z)}{\mathcal{A}_{{\bf n}_0 + k{\bf p},j}(z)}\right| =
\frac{1}{m} \log \left|\prod_{l\in
L}\mathcal{A}^{(l)}_{j}(z)\right| \,,\quad \mathcal{K} \subset
{\mathbb{C}} \setminus (\supp{\sigma_j}\cup\supp{\sigma_{j+1}}),
\]
since $|{\bf n}_{1}| = r |{\bf p}_{1}| + {\mathcal{O}}(1) = rm +
{\mathcal{O}}(1), |{\bf n}_{1}| \to \infty.$ \end{proof}

The function appearing on the right hand side of (\ref{eq:98})
corresponds with the one on the right hand side of (\ref{conver1})
associated to the vector equilibrium problem with interaction
matrix $\mathcal{C}$ constructed taking $p_{1,k}=1/(m_{1}+1)$,
$0\leq k\leq m_{1}$, and $p_{2,k}=1/(m_{2}+1)$, $0\leq k\leq
m_{2}$. In that case, for each $j=-m_{2}-1,\ldots,m_{1}$, we have
\[
G_{j}(z)=\prod_{l\in L}|\mathcal{A}^{(l)}_{j}(z)|^{1/m}\,,\qquad
z\in{\mathbb{C}}\setminus
(\widetilde{\Delta}_{j}\cup\widetilde{\Delta}_{j+1})
\]
$(\Delta_{-m_{2}-1}=\Delta_{m_{1}+1}=\emptyset),$ where
$m=\lcm{(m_{1}+1,m_{2}+1)}$.

\section{Application to mixed type Hermite Pad\'{e} approximation}

Let $S^1= \mathcal{N}(\sigma_0^1,\ldots,\,\sigma^1_{m_1}), S^2
=\mathcal{N}(\sigma^2_0,\ldots,\,\sigma^2_{m_2}), \sigma_0^1 =
\sigma_0^2$ be given. Let us  introduce the row vectors
\[ {\mathbb{U}} =  (1,\widehat{s}^{2}_{1,1},\ldots,\widehat{s}^{2}_{1,m_2}),
\qquad {\mathbb{V }} =
(1,\widehat{s}^{1}_{1,1},\ldots,\widehat{s}^{1}_{1,m_1})
\]
and the $(m_2+1) \times (m_1 +1)$ dimensional matrix
\[ {\mathbb{W}} = {\mathbb{U}}^t {\mathbb{V}},
\]
where the super-index $t$ means taking transpose. Define the
matrix Markov type function
\[ \widehat{\mathbb{S}}(z) = \int \frac{{\mathbb{W}}(x)d\sigma_0^2 (x)}{z-x}
\]
understanding that integration is carried out entry by entry on
the matrix $\mathbb{W}$.

Fix ${\bf
n}_1=(n_{1,0},\,n_{1,1},\ldots,\,n_{1,m_1})\in{{\mathbb{Z}}}_+^{m_1+1}$
and ${\bf
n}_2=(n_{2,0},\,n_{2,1},\ldots,\,n_{2,m_2})\in{{\mathbb{Z}}}_+^{m_2+1},
|{\bf n}_2| = |{\bf n}_1| -1$. It is easy to see that there exists
a non zero vector polynomial
\[{\mathbb{A}}_{\bf n} = (a_{{\bf
n},0},\ldots,a_{{\bf n},m_1}), \qquad \deg (a_{{\bf n},k}) \leq
n_{1,k} -1, \qquad k=0,\ldots,m_1,
\]
such that
\begin{equation} \label{eq:26}
\widehat{\mathbb{S}}(z){\mathbb{A}}^t_{{\bf n}}(z) -
{\mathbb{D}}^t_{{\bf n}}(z) =
({\mathcal{O}}(1/z^{n_{2,0}+1}),\ldots,{\mathcal{O}}(1/z^{n_{2,m_2}+1}))^t
=: {\mathcal{O}}(1/z^{{\bf n}_{2}+1}), \qquad z \to \infty,
\end{equation}
where ${\mathbb{D}}_{{\bf n}} = (d_{{\bf n},0},\ldots,d_{{\bf
n},m_{2}})$ is some vector polynomial. When $m_2 = 0$, this
construction is called type I Hermite-Pad\'{e} approximation. If
$m_1 =0$ it is called of type II. When $m_1 = m_2 =0$ it reduces
to diagonal Pad\'{e} approximation. This definition is  of mixed
type.

\begin{lem} For $j=0,\ldots,m_2, (\widehat{s}^1_{1,0} \equiv 1)$
\begin{equation}\label{definition**}
 \int x^{\nu} \sum_{k=0}^{m_1} a_{{\bf n},k}(x)\widehat{s}^1_{1,k}(x) d s^2_{j}(x)  =0, \qquad \nu =0,\ldots,n_{2,j}
 -1.
\end{equation}

\end{lem}
\begin{proof} In fact, notice that according to (\ref{eq:26}), for each $\nu,
0\leq \nu \leq n_{2,j} -1, j=0,\ldots,m_2,$
\[ z^{\nu}\left(\sum_{k=0}^{m_1} a_{{\bf n},k}(z) \int \frac{\widehat{s}^2_{1,j}(x)\widehat{s}^1_{1,k}(x)d\sigma^2_0(x)}{z-x}  - d_{{\bf
n},j}(z)\right) = {\mathcal{O}}\left(1/z^2\right), \qquad z \to
\infty,
\]
$(\widehat{s}^2_{1,0}  \equiv 1)$ and the function on the left
hand side is holomorphic in $\overline{\mathbb{C}} \setminus
\mbox{Co}(\sop(\sigma^2_{0}))$. Using Lemma \ref{reduc}, we obtain
(\ref{definition**}).
\end{proof}

Because of this Lemma, we see that ${\mathbb{A}}_{\bf n}$ is an
${\bf n}$-th mixed type multiple orthogonal polynomial with
respect to the pair $(S^1,S^2)$ and in the sequel we assume that
it is ``monic''. If
\[{\mathbb{B}}_{\bf n} = (b_{{\bf n},0}, \ldots,b_{{\bf n},m_2}),
\qquad \deg(b_{{\bf n},j}) \leq n_{2,j} -1, \qquad j=0,\ldots,m_2,
\]
denotes  a generic vector polynomial with the indicated degrees,
(\ref{definition**}) may be rewritten in matrix form as
\begin{equation} \label{eq:25}
\int {\mathbb{B}}_{{\bf n}}(x) {\mathbb{W}}(x) {\mathbb{A}}_{{\bf
n}}^t(x) d\sigma_0^2(x) = 0, \qquad \mbox{for all}\,\,
{\mathbb{B}}_{\bf n}.
\end{equation}

Fix $j \in \{0,\ldots,m_2\}$. For each $k \in \{-1,\ldots,-j-1\}$
define
\[ \Omega^j_k = \{z \in {\mathbb{C}} \setminus \cup_{i=0}^{-j-1} \Delta_i: U_{k}^{\overline{\mu}}(z) <
U_{i}^{\overline{\mu}}(z), i = -1,\ldots,-j-1, i\neq k\}, \quad
\Omega^{0}_{-1} = \mathbb{C} \setminus (\Delta_{0} \cup
\Delta_{-1}).
\]
Set
\[\chi_j(z) := \min\{U_k^{\overline{\mu}}(z): k = -1,\ldots,-j-1\}\]
and
\[ (\mathcal{R}_{{\bf n},0},\ldots,\mathcal{R}_{{\bf n},m_2})^t:= \widehat{\mathbb{S}}(z){\mathbb{A}}^t_{{\bf n}}(z) -
{\mathbb{D}}^t_{{\bf n}}(z).
\]

\begin{thm}\label{markov}
Let $\Lambda =
\Lambda(p_{1,0},\ldots,p_{1,m_1};p_{2,0},\ldots,p_{2,m_2}) \subset
\mathbb{Z}_+^{m_1+1}(\bullet) \times
\mathbb{Z}_+^{m_2+1}(\bullet), (S^1,S^2)\in\mbox{\bf {Reg}},$
$S^1= \mathcal{N}(\sigma^1_0,\ldots,\sigma^1_{m_1}),$ and $S^2=
\mathcal{N}(\sigma^2_0,\ldots,\sigma^2_{m_2})$ be given. Then for
each $j \in \{0,\ldots,m_2\}$
\begin{equation} \label{eq:m1}
\lim_{{\bf n} \in \Lambda}|\mathcal{R}_{{\bf n},j}(z)|^{1/|{\bf
n}_1|}= \exp(-\chi_j(z)), \qquad \mathcal{K} \subset
\cup_{k=-1}^{-j-1} \Omega_k^j,
\end{equation}
and
\begin{equation} \label{eq:m2}
\lim_{{\bf n} \in \Lambda}|\mathcal{R}_{{\bf n},j}(z)|^{1/|{\bf
n}_1|} \leq \exp(-\chi_j(z)), \qquad {\mathcal{K}} \subset
{\mathbb{C}} \setminus ( \cup_{k=0}^{-j-1} \Delta_k).
\end{equation}
In particular, if $p_{2,0} = \cdots = p_{2,m_2} = 1/(m_2+1)$, then
\begin{equation}\label{eq:m3}
\lim_{{\bf n} \in \Lambda}|\mathcal{R}_{{\bf n},j}(z)|^{1/|{\bf
n}_1|}= \exp(-U_{-1}^{\overline{\mu}}(z)), \qquad \mathcal{K}
\subset {\mathbb{C}} \setminus ( \cup_{k=0}^{-j-1} \Delta_k).
\end{equation}
$\overline{\mu}= \overline{\mu}(\mathcal{C})=
(\overline{\mu}_{-m_2},\ldots,\overline{\mu}_{m_1})$ is  the
equilibrium vector measure   and
$(\omega_{-m_2}^{\overline{\mu}},\ldots,\omega_{m_1}^{\overline{\mu}})$
is the system of equilibrium constants for the vector potential
problem determined by the interaction matrix $\mathcal{C}$ defined
in $(\ref{matriz})$  on the system of compact sets $E_k = \sop
(\sigma^1_k), k=0,\ldots,m_1, E_k = \sop (\sigma^2_{-k}),
k=-m_2,\ldots,0$.
\end{thm}

\begin{proof} Notice that (\ref{eq:25}) implies that
\[ \widehat{\mathbb{S}}(z){\mathbb{A}}^t_{{\bf n}}(z) - \int \frac{\mathbb{W}(x)({\mathbb{A}}^t_{{\bf
n}}(z) - {\mathbb{A}}^t_{{\bf n}}(x))d \sigma_0^2(x)}{z-x} = \int
\frac{\mathbb{W}(x) {\mathbb{A}}^t_{{\bf n}}(x)d
\sigma_0^2(x)}{z-x} = \mathcal{O}(1/z^{{\bf n}_2 +1}), \quad z \to
\infty,
\]
and taking
\[ {\mathbb{D}}^t_{{\bf n}}(z) = \int \frac{\mathbb{W}(x)({\mathbb{A}}^t_{{\bf
n}}(z) - {\mathbb{A}}^t_{{\bf n}}(x))d \sigma_0^2(x)}{z-x}
\]
we obtain an integral expression for the remainder in
(\ref{eq:26}).

Then
\[ (\mathcal{R}_{{\bf n},0}(z),\ldots,\mathcal{R}_{{\bf n},m_2}(z))^t = \int
\frac{\mathbb{W}(x) {\mathbb{A}}^t_{{\bf n}}(x)d
\sigma_0^2(x)}{z-x}.
\]
In scalar form this says that
\[ \mathcal{R}_{{\bf n},j}(z) = \int \frac{\mathcal{A}_{{\bf
n},0}(x)}{z-x} ds^2_j(x), \qquad j =0,\ldots,m_2.
\]
Notice that (see (\ref{defnlnjnegative}))
\[ \mathcal{R}_{{\bf n},0}(z) = \mathcal{A}_{{\bf n},-1}(z).
\]
Let us establish a connection between the remainders
$\mathcal{R}_{{\bf n},j}(z)$ and the forms $\mathcal{A}_{{\bf
n},k}(z)$ with negative indices $k \in \{ -1,\ldots,-j-1\}.$

Fix $j \in \{1,\ldots,m_2\}.$ We have
\[  (-1)^j\mathcal{R}_{{\bf n},j}(z) = \int\cdots\int \frac{\mathcal{A}_{{\bf
n},0}(x_0)\,d\sigma^2_0(x_0)\cdots
d\sigma^2_j(x_j)}{(z-x_0)(x_1-x_0)\cdots (x_j-x_{j-1})},
\]
and
\[ \mathcal{A}_{{\bf n},-j-1}(z) = \int\cdots\int \frac{\mathcal{A}_{{\bf
n},0}(x_0)\,d\sigma^2_0(x_0)\cdots
d\sigma^2_j(x_j)}{(x_1-x_0)\cdots (x_j-x_{j-1})(z-x_j)}.
\]
Consequently,
\[ (-1)^j\mathcal{R}_{{\bf n},j}(z) - \mathcal{A}_{{\bf n},-j-1}(z) =
\int\cdots\int \frac{-(x_j-x_0)\mathcal{A}_{{\bf
n},0}(x_0)\,d\sigma^2_0(x_0)\cdots
d\sigma^2_j(x_j)}{(z-x_0)(x_1-x_0)\cdots (x_j-x_{j-1})(z-x_j)}.
\]
Since $x_j -x_0 = x_j - x_{j-1} + x_{j-1}- \cdots - x_1 + x_1
-x_0$, substituting this in the previous formula, we obtain
\[ \mathcal{A}_{{\bf n},-j-1}(z) =  \sum_{k=0}^{j-1}(-1)^k
\langle \sigma^2_j, \ldots,\sigma^2_{k+1}\rangle^{\widehat{}}(z)
\mathcal{R}_{{\bf n},k}(z) + (-1)^j\mathcal{R}_{{\bf n},j}(z) .
\]

We have a triangular scheme of linear equations whose coefficients
do not depend on ${\bf n}$. We can solve for $\mathcal{R}_{{\bf
n},j}$ in terms of $\mathcal{A}_{{\bf n},-1},\ldots,
\mathcal{A}_{{\bf n},-j-1}$. Using (\ref{formula}) one obtains
that for each $j \in \{0,\ldots,m_2\}$ (when $j=0$ the sum below
is empty)
\[ \mathcal{R}_{{\bf n},j}(z) = \sum_{k=1}^{j}
(-1)^{k-1} \langle \sigma^2_k,
\ldots,\sigma^2_{j}\rangle^{\widehat{}}(z)\mathcal{A}_{{\bf
n},-k}(z) + (-1)^{j}\mathcal{A}_{{\bf n},-j-1}(z).
\]
Taking (\ref{conver1})  into consideration, on $\Omega^j_{-k}$ the
term containing   $\mathcal{A}_{{\bf n},-k}$ dominates the sum
(notice that $\langle \sigma^2_k,
\ldots,\sigma^2_{j}\rangle^{\widehat{}}(z) \neq 0, z \in
\mathbb{C} \setminus \Delta_{-k}$) and (\ref{eq:m1}) immediately
follows . On the complement of $\cup_{k=-1}^{-j-1} \Omega_k^j$
there is no dominating term and all we can conclude from the
previous equality is (\ref{eq:m2}).

Let $p_{2,0}=\cdots=p_{2,m_2} = 1/(m_2+1)$. In this case, on
$\mathbb{C} \setminus \cup_{k=0}^{-j-1} \Delta_k$ we have that
$U^{\overline{\mu}}_{-1}(z) < U^{\overline{\mu}}_{-2}(z) < \cdots
< U^{\overline{\mu}}_{-j-1}(z)$ (see third sentence before
Corollary \ref{polinom}) and (\ref{eq:m3}) follows from
(\ref{eq:m1}).
\end{proof}

\begin{rem} Fix $j \in \{0,\ldots,m_2\}$. For each $k \in \{-1,\ldots,-j-1\}$
we could have defined
\[ \Omega^j_k = \{z \in {\mathbb{C}} \setminus \cup_{i=0}^{-j-1} E_i: U_{k}^{\overline{\mu}}(z) <
U_{i}^{\overline{\mu}}(z), i = -1,\ldots,-j-1, i\neq k\}, \quad
\Omega^{0}_{-1} = \mathbb{C} \setminus (E_{0} \cup E_{-1}).
\]
Taking into account that the polynomials $Q_{{\bf n},i}$ and the
forms $\mathcal{A}_{{\bf n},i}$ may have at most one zero in each
of the connected components of $\Delta_i \setminus E_i$, one can
prove in place of (\ref{eq:m1})-(\ref{eq:m3}) convergence in
capacity on each compact subset of the corresponding regions. \fp
\end{rem}

We say that $\mathcal{I}_1 \subset \mathbb{Z}_+^{m_1+1}(\bullet)$
is a complete, ordered, sequence of multi-indices if:
\begin{itemize}
\item[a)] For each $n \in {\mathbb{Z}}_+$, there exists a unique $
{\bf n}_1 \in \mathcal{I}_1$ such that $|{\bf n}_1| = n$.
\item[b)] Any two multi-indices in $\mathcal{I}_1$ are ordered in
the sense that all components of one of them are less than or
equal to the corresponding components of the other one, or they
are identical.
\end{itemize}

Fix $\mathcal{I}_1 \subset \mathbb{Z}_+^{m_1+1}(\bullet),
\mathcal{I}_2 \subset \mathbb{Z}_+^{m_2+1}(\bullet),$ two
complete, ordered sequences of multi-indices. Each $n \in
\mathbb{Z}_+$ determines a unique ${\bf n}_1 \in \mathcal{I}_1$
and ${\bf n}_2 \in \mathcal{I}_2$ such that $n = |{\bf n}_1|=|{\bf
n}_2| +1$. The corresponding ``monic'' mixed type multiple
orthogonal polynomials we denote by $\mathbb{A}_n$. We can
interchange the roles of the Nikishin systems $S^1,S^2,$ and
determine a sequence of ``monic'' mixed type multiple orthogonal
polynomials which we denote $\mathbb{B}_n$. It is easy to verify
that the sequences $\{\mathbb{A}_n\}, \{\mathbb{B}_n\}, n \in
\mathbb{Z}_+$ are bi-orthogonal. That is,
\begin{equation} \label{eq:m25}
\int {\mathbb{B}}_{n'}(x) {\mathbb{W}}(x) {\mathbb{A}}_{{n}}^t(x)
d\sigma_0^2(x)  \left\{
\begin{array}{cc}
= 0, & n\neq n', \\
\neq 0, & n = n'.
\end{array}
\right.
\end{equation}
The inequality in (\ref{eq:m25}) is a consequence of Lemma
\ref{ulises}. With the same hypothesis, all the results of this
paper hold true for the sequence $\{\mathbb{B}_n\}, n \in
\mathbb{Z}_+$.

\end{document}